\documentclass[%
 reprint,
preprintnumbers,
 amsmath,amssymb,
 aps,
]{revtex4-1}

\usepackage{amsmath}%
\usepackage{amsfonts}%
\usepackage{amssymb}%
\usepackage{graphicx}
\usepackage{verbatim}   
\usepackage[dvips]{color}      

\usepackage{subfigure}  
\usepackage{dcolumn}
\usepackage{bm}

\begin{document}

\preprint{APS/123-QED}

\title[SKEW-VARICOSE INSTABILITY IN 2D GENERALIZED SH MODEL EQUATIONS]{SKEW-VARICOSE INSTABILITY IN TWO DIMENSIONAL\\ GENERALIZED SWIFT--HOHENBERG EQUATIONS}

\author{J. A. Weliwita}
 \email{mmjaw@leeds.ac.uk}
\author{A. M. Rucklidge}%
\homepage{http://www.maths.leeds.ac.uk/~alastair/}
\author{S.M. Tobias}%
\homepage{http://www.maths.leeds.ac.uk/~smt/}
\affiliation{%
Department of Applied Mathematics, University of Leeds, Leeds LS2 9JT, UK.\\
}%

\date{\today}
\begin{abstract}
We apply analytical and numerical methods  to  study the linear stability of stripe
patterns  in two generalizations of the two-dimensional Swift--Hohenberg equation that
include coupling to a mean flow. A projection operator is included in our models  to allow  exact stripe solutions.  In the generalized models, stripes become unstable to the skew-varicose, oscillatory skew-varicose  and cross-roll instabilities, in addition to the usual Eckhaus and zigzag instabilities. We analytically derive  stability boundaries for the skew-varicose instability in various cases, including several asymptotic limits. We also use  numerical techniques to determine eigenvalues and hence stability boundaries of  other instabilities. We extend our analysis to  both stress-free and no-slip boundary conditions and we note a cross over from the behaviour characteristic of no-slip to that of stress-free boundaries as the coupling to the mean flow increases or as the Prandtl number decreases. 

Close to the critical value of the bifurcation parameter, the skew varicose instability has the same curvature as the Eckhaus instability provided the coupling to the mean flow is greater than a critical value. The region of stable stripes is completely eliminated by the  cross-roll instability  for large coupling to the mean flow.

\end{abstract}

\pacs{Valid PACS appear here}

\maketitle
\section{\label{sec:level1}Introduction\protect\\}
Spiral defect chaos (SDC) was discovered in low-viscosity convection nearly 20 years
ago~\cite{Morris:1993}, and yet much of the detail of its origin remain
unexplained. The Swift--Hohenberg equation (SHE), originally proposed as a model
equation for fluctuations near the onset of convection~\cite{swift:1977}, has
yielded considerable insight into the question of wavenumber selection in pattern
formation~\cite{Cross:1993}. However the SHE has a major drawback as a model of
low-viscosity convection: mean flows are excluded. The importance of the large scale mean flow on the stability of convection rolls was first investigated by Siggia and Zippelius~\cite{Siggia:1981}. Mean flows are known to play
an important role in the dynamics of spiral defect chaos~\cite{Chiam2003b}, and
so two related generalizations of the SHE have been developed to include the
effects of mean flows~\cite{Greenside:1988,Xi:1993}.

 SDC was found numerically
in solutions of one of these generalized Swift--Hohenberg equations~\cite{Xi:1993,cross:1995} and
of the Boussinesq equations for convection~\cite{decker1994,Chiam2003b,pesch1996,Morris:1996,Bodenschatz:2000}. 
The generalized SHE can reproduce some qualitative features of spiral defect chaos,
at least as transients, though direct comparison between the model and convection appears not to be
appropriate~\cite{pesch:2002}.

Including mean flows in these models allows an interesting long-wavelength instability, the skew-varicose instability (SVI)~\cite{Busse:1978} whose spatial dependence is neither entirely
longitudinal nor transverse to the local pattern wavevector. The SVI resembles the  Eckhaus instability but the most unstable modes are those at an angle to the original roll
axes. In the SVI,  rolls bend and become irregular in order to decrease
their effective wavenumber, and often dislocation pairs will form
\cite{Cross:1993}. Indeed, chaotic spiral patterns have been observed during the
transition from the conducting state to rolls. Mean flows are also crucial to the formation of SDC: Chiam {\it et al.}~\cite{Chiam2003b} showed that spiral defect chaos collapses to a stationary pattern when the mean flow is quenched. The onset of SDC and defect chaos has therefore been  tentatively associated with the occurrence of the SVI~\cite{gollub:1982,Bodenschatz:2000,Hu}. It is this connection that motivates this detailed investigation into the \hbox{SVI}.

The classical
problem of linear stability of convection rolls with stress-free horizontal
boundaries near the onset of convection has been studied by a number of authors~\cite{Zippelius:1983,Bolton:1985,Bernoff:1994,manneville:2006}, while the linear stability of convection with no-slip boundaries~\cite{Busse:1979} is much less well understood.
An additional problem that occurs in the analysis of the  SVI is that there is no consistent 
relative  scaling of lengths parallel and perpendicular to the roll axes, owing to the singular nature of the slow length scale expansion of the stability problem, as detailed below.

Zippelius  and Siggia~\cite{Zippelius:1983} were the first
to study the \hbox{SVI}. They used asymptotic scalings to derive modulation equations from the original equations of
convection with stress-free boundaries, and showed that mean flows have an 
important role. The stability analysis includes some restrictive assumptions about the wavenumbers of perturbations, and thus does not capture all mechanisms of instability~\cite{Busse:1984}. Bernoff~\cite{Bernoff:1994} extended the work of Zippelius and Siggia  by deriving a  set of equations for the amplitudes of rolls and the mean flow, retaining terms with  mixed asymptotic orders. The effect of mean flow on the SVI for  stress-free and no-slip convection has also been considered by  Ponty {\it et al.}~\cite{ponty:1997}, who derived the SVI boundary for stress-free and no-slip convection from the  Boussinesq equations for convection. However, the  perturbation terms that they included in their calculations do not adequately capture the mechanism of the SVI in all circumstances~\cite{mathews}.  Milke~\cite{Milke:1997} carried out a more
comprehensive study of stability of rolls by means of a Lyapunov--Schmidt
reduction. He recovered domains in Rayleigh, Prandtl and wave number space where convection rolls are unstable, the so-called Busse balloon.

In this paper, we revisit the two generalized SH models that include a
coupling to the mean flow~\cite{Green:1985,Greenside:1988,Xi:1993,manneville:1983}, and analyse in detail  how the
skew-varicose instability depends on the parameters in the models.

Although the model equations of  interest have been  well
studied, a complete analysis of the SVI in the models is not available. One difficulty in the analysis comes from 
the contribution of  terms proportional to $k^{2}l^{2}/(k^{2}+l^{2})$, where $(k,l)$ is the small
wavevector of the perturbations associated with the SVI. Terms like this are responsible for the absence of a consistent asymptotic scaling in the limit of small amplitude and small $k$ and $l$~\cite{Hoyle:2006}. This
difficulty is resolved in this paper.

The remaining part of this paper is organized as follows:
in  section (II), we present the two generalizations of the Swift--Hohenberg equation, both of which incorporate the effect of a mean flow. We discuss two operators: $\mathcal{P}_{\alpha}$ is a projection operator included in the model  to allow  exact stripe solutions, and  $\mathcal{F}_{\gamma}$  is a filtering operator that suppresses the cross-roll instability. $\mathcal{F}_{\gamma}$ is present in an original formulation of the models~\cite{Green:1985}; $\mathcal{P}_{\alpha}$  was suggested by Ian Melbourne~\cite{Ian}. In section (III), we present a detailed analysis of the  linear stability of the stripe solution, showing how the SVI can be located analytically. We explore the importance of the mean flow in long-wavelength instabilities. The structure of the  maximum eigenvalue for  instability for small $k$ and $l$ is also presented.

 In  section (IV), we  analyze the zigzag and Eckhaus long-wavelength instabilities. Section (V) describes the interesting skew-varicose instability. We also consider various asymptotic limits.
 
The work is extended in section (VI) to the  stress-free boundary condition case, for which there are skew-varicose and oscillatory skew-varicose instabilities. Section (VII) covers short-wavelength instabilities (the cross-roll instability) for stress-free and no-slip boundary conditions. Some numerical illustrations of stability boundaries with  growth rates of perturbations are illustrated in section (VIII), where we also discuss some curious behavior of the growth rates of perturbations. Finally, in section (IX), we show how the region of stable stripes is  affected by the coupling to the mean flow and how it can be completely eliminated if the coupling is strong enough, in the stress-free case. We conclude in section (X).
\section{Problem Definition}
\label{sec:2}
In this section, we set out the two models we will investigate, and discuss 
basic properties of the models. Both models are generalizations of the 
two-dimensional Swift--Hohenberg equation, where a real field~$\psi(x,y,t)$,
representing the amplitude of convection, couples to a mean flow given by a 
stream function $\zeta(x,y,t)$ and its vertical 
vorticity~$\omega(x,y,t)=-\nabla^2\zeta$.  


\subsection{Description of Models}

The standard Swift--Hohenberg equation~\cite{swift:1977} is
 \begin{equation}
 \frac{\partial \psi}{\partial t} =
   \left[\mu-(1+\nabla^2)^2\right]\psi - \psi^3, 
 \label{eq:she}
 \end{equation}
where $\psi(x,y,t)$ represents the pattern-forming field, and 
$\mu$ is the driving parameter (in convection, $\mu$ represents the temperature difference between the top and
the bottom layer), taking the value
zero at the onset of pattern formation. Both models introduce a 
$\left(\mathbf{U}\cdot\mathbf{\nabla}\right)\psi$ term to the right-hand side 
of~(\ref{eq:she}), where $\mathbf{U}(x,y,t)$ is a mean flow calculated from the 
stream function~$\zeta(x,y,t)$: 
 \[
 \mathbf{U} = \left(\frac{\partial\zeta}{\partial y},
                   -\frac{\partial\zeta}{\partial x}\right).
 \]
The mean flow has vertical vorticity $\omega(x,y,t)=-\nabla^2\zeta$. The way that
vorticity is generated by nonlinear forcing from~$\psi$ differs in the two 
models.

In the first model~\cite{Xi:1993},
the vertical vorticity~$\omega(x,y,t)$ has its own independent dynamics:
 
 \begin{align}
  \label{eq:model1psi}
 \frac{\partial \psi}{\partial t} + 
   \left(\mathbf{U}\cdot \mathbf{\nabla}\right )\psi &= 
   \left[\mu-(1+\nabla^2)^2\right]\psi - 
     \mathcal{P_{\alpha}}\left(\psi^3 \right),\\
 \left[\frac{\partial}{\partial t} - Pr(\nabla^2 - c^2)\right]\omega &=
 - g_m \mathcal{F}_{\gamma}
     \left[\nabla(\nabla^2\psi)\times \nabla \psi
     \right]\cdot\mathbf{\widehat{z}}, 
 \label{eq:model1omega}
 \end{align}
where $Pr$, $c$ and $g_m$ are parameters. The Prandtl number~$Pr$ (the ratio
between kinematic viscosity and thermal diffusivity) is effectively a viscosity
parameter for the mean flow, which plays a much greater role in low Prandtl number
convection. Indeed, in the limit of large $Pr$, the vertical vorticity is
hardly excited and the dynamics of~$\psi$ becomes purely
relaxational~\cite{Busse:1971,Cross:1993,pesch:2002}, reducing model~1 back to the
\hbox{SHE}. The coefficient~$g_m$ is a coupling parameter that
controls the strength of the mean flow effects relative to the ordinary
Swift--Hohenberg nonlinear term~$-\psi^3$. The parameter~$c$ models the effect
of top and bottom boundary conditions on the vertical vorticity, with $c=0$
corresponding to stress-free boundary conditions (where
$\omega=\mbox{constant}$ is an allowed solution of the linearised vorticity
equation), and $c\neq0$ corresponding to more realistic no-slip boundary
conditions ($\omega$~decays to zero in the absence of nonlinear forcing). The
operators $\mathcal{P_{\alpha}}$ and $\mathcal{F}_{\gamma}$ are 
explained in more detail below.

The second model~\cite{Green:1985,Greenside:1988,manneville:1983} has the vertical vorticity slaved to the
nonlinear driving term:
 \begin{align}
 \label{eq:model2psi}
 \frac{\partial \psi}{\partial t} + 
   \left(\mathbf{U}\cdot \mathbf{\nabla}\right )\psi &=
   \left[\mu-(1+\nabla^2)^2\right]\psi - 
     \mathcal{P_{\alpha}}\left(\psi^3 \right),\\
 \omega &= -g\mathcal{F}_{\gamma}\left[\mathbf{\nabla}(\nabla^2\psi) \times \mathbf{\nabla}\psi \right]\cdot\mathbf{\widehat{z}},
 \label{eq:model2omega}
 \end{align}
so the vertical vorticity responds instantly to the nonlinear driving. The
coefficient~$g$ is a coupling parameter that controls the relative strength of
mean flow effects compared to the ordinary nonlinearity.

We use two operators $\mathcal{P_{\alpha}}$ and          
$\mathcal{F}_{\gamma}$ to ease the analysis and to ensure that the SVI is not 
pre-empted by other instabilities. They both act as filters in Fourier space; 
the first is a projection:
 \[ 
 \mathcal{P_{\alpha}}(e^{i\mathbf{K}\cdot\mathbf{x}}) = 
   \left\{ \begin{array}{ll}
         e^{i\mathbf{K}\cdot\mathbf{x}} & \mbox{if $\scriptstyle\left|\mathbf{K}\right|\leq\alpha$};\\
         0 & \mbox{if $\scriptstyle\left|\mathbf{K}\right|>\alpha$}.
           \end{array} \right. 
 \]
By setting $\alpha=2.5$, we allow stripes with a single Fourier mode with 
wavenumber close to one to be exact solutions of the
PDEs~(\ref{eq:model1psi}--\ref{eq:model1omega}) and
(\ref{eq:model2psi}--\ref{eq:model2omega})~\cite{Ian}. This is shown in
section~(\ref{sec:basic}) below. The second operator, $\mathcal{F}_{\gamma}$,
reduces short-wavelength modulations of the mean flow~\cite{Green:1985}: in
Fourier space, the operator is defined by
 \[
 \mathcal{F}_{\gamma}(e^{i\mathbf{K}\cdot\mathbf{x}})
 =e^{-\gamma^{2}|\mathbf{K}|^{2}}e^{i\mathbf{K}\cdot\mathbf{x}}.
 \]
Throughout this work, we set $\gamma=2.5$. The effect of this filtering on short-wavelength instabilities, particularly the cross-roll instability, is discussed in  section (\ref{cr}).


\begin{figure}[t]
\scalebox{0.37}{\includegraphics*[viewport=0.8in 1.05in 9.25in 6.45in]{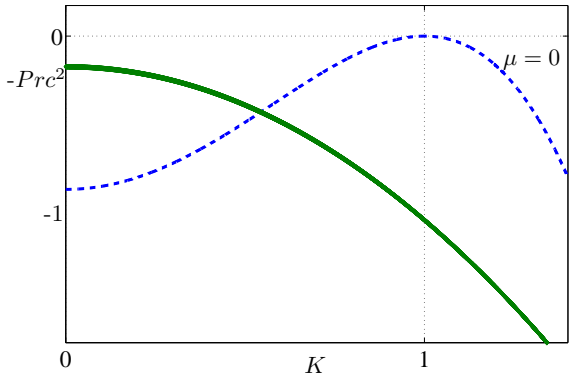}}
\caption{(Color online) Growth rates as functions of wavenumber $K=|\mathbf{K}|$.
 Dashed blue curve: $\sigma_{1}$ for $\mu=0$, which peaks at $K=K_{critical}=1$.
 Solid green curve: $\sigma_{2}$, which takes the value $-Pr\,c^{2}$ at $K=0$ and
 decreases as $K$ increases. 
 }
 \label{fig:1}
\end{figure}

\subsection{Basic Properties of the Models}
\label{sec:basic}

For model~1, linearisation yields:
 \[
 \frac{\partial \psi}{\partial t} = \left[\mu-(1+\nabla^2)^2\right]\psi
 \quad\hbox{and}\quad
 \frac{\partial\omega}{\partial t}=Pr\left(\nabla^2 - c^2\right)\omega\;;
 \]
only the first equation is relevant to model~2.
Normal mode solutions to these linear equations are given by
  $\psi=F_{1}e^{\sigma_{1} t+i\mathbf{K}\cdot\mathbf{x}}$ and
  $\omega=F_{2}e^{\sigma_{2} t+i\mathbf{K}\cdot\mathbf{x}}$, where
  $\sigma_{1}$ and $\sigma_{2}$ are growth rates,
  $\mathbf{K}$~is a wavevector, 
  and $F_{1}$ and
  $F_{2}$ are constants. Substituting these into the linearized
  equations gives
  $\sigma_{1}=\mu-\left(1-K^{2}\right)^{2}$ and
  $\sigma_{2}=-Pr\left(K^{2}+c^{2}\right)$.
These are shown in Figure~\ref{fig:1} for $\mu=0$, when the trivial
solution is marginally stable.~The most unstable wavenumber is
$K_{critical}=1$, and for $\mu>0$, a band of wavenumbers close to $K=1$ is
linearly unstable, signaling the onset of pattern formation. The
vorticity~$\omega$ (in model~1) is always linearly damped, unless $c=0$.

We define $q=K-K_{critical}$ and thus the trivial solution loses stability for
any $q$ at $\mu_{Existence}$, where 
 \[ \mu_{Existence}=(1-(1+q)^{2})^{2}.  \]
Note
that as $q\rightarrow0$, $\mu_{Existence}\rightarrow4q^{2}$.

In this paper, we are interested in the stability of the nonlinear equilibrium 
stripe solution of the PDEs. We note that 
 \begin{eqnarray}
\psi_{0}=\sqrt{\beta}\left(e^{i(1+q)x}+e^{-i(1+q)x}\right),\quad\omega_0=0
\nonumber\\
\hbox{with}\quad\beta=\frac{\mu-(1-(1+q)^{2})^{2}}{3},
\label{eq:nonlinearstripes}
\end{eqnarray} is an \textit{exact} solution of both models. This can be shown by substituting 
the expressions for $\psi_{0}$ and $\omega_0$ into the PDEs:
 \begin{eqnarray*}
 0=\mu\sqrt{\beta}\left(e^{i(1+q)x}+e^{-i(1+q)x}\right)
 - (1+\nabla^2)^2\left(\sqrt{\beta}\left(e^{i(1+q)x}\right.\right.\\
 \left.\left.+e^{-i(1+q)x}\right)\right)
 - \mathcal{P}_{\alpha}\left(\left(\sqrt{\beta}\left(e^{i(1+q)x}+e^{-i(1+q)x}\right)\right)^3\right).
 \end{eqnarray*} The nonlinear term in the vorticity equation is zero.
Next, we note that 
  \[(1+\nabla^2)^2e^{\pm{i}(1+q)x}=(1-(1+q)^2)^2e^{\pm{i}(1+q)x}\]
and
\begin{eqnarray*}
 \mathcal{P}_{2.5}\left[\left(\sqrt{\beta}\left(e^{i(1+q)x}+e^{-i(1+q)x}\right)\right)^3\right] = 3\beta\sqrt{\beta}\nonumber\\
 \times\left(e^{i(1+q)x}+e^{-i(1+q)x}\right),
 \end{eqnarray*} provided that $q$ is between $-0.167$ and $1.5$. The consequence is that $(\psi_{0},\omega_0)$ is an exact solution.
In the limit of small~$q$, we have $\beta=(\mu-4q^2)/3$, so stripe solutions 
exist when $\mu>4q^2$. 

The advantage of using the projection~$\mathcal{P}_{\alpha}$ is that it allows 
this exact stripe solution of the PDEs~\cite{Ian}. The alternative would be to consider 
the limit of small $\mu$ and~$\beta$, but by having an exact solution, which 
matches the asymptotic result we would obtain without the projection, we do 
not have to be concerned with the relative sizes of these parameters compared 
to other small parameters that will be introduced below. 

The two models can be related to each other close to onset, regardless of the
projection and filtering. By scaling $\mu=\mathcal{O}\left(\epsilon^{2}\right)$
with $\epsilon\ll1$, introducing a slow time scale $\partial/\partial
t\rightarrow \epsilon^{2}\partial/\partial t$, and assuming that the
wavenumbers that are excited in the vorticity variable~$\omega$ are of order
$\mathcal{O}(\epsilon)$, the largest term on the left-hand side of the
vorticity equation~(\ref{eq:model1psi}) in model~1 is $Pr\,c^2\omega$. Thus
model~1 reduces to model~2 in this limit, with the relation
$g=g_m/(Pr\,c^2)$.


\section{Linear stability of stripes}

The linear stability theory for stripes in the SHE is
well known~\cite{Hoyle:2006}: stripes with wavenumber $1+q$ exist provided
$\beta>4q^2$, and they are stable with respect to the Ekhaus and zigzag
instabilities provided $\beta>12q^2$ and $q>0$. Once mean flows are included, the zigzag instability needs to be modified at finite Prandtl number by the
presence of mean-flow modes with non-zero vertical vorticity
\cite{Zippelius:1983,Busse:1984}. The Eckhaus instability is unchanged. We consider the stability of the stripe
solutions with respect to  long-wavelength perturbations, deriving
three (model~1) or two (model~2) linear ODEs for the perturbation amplitudes,
and so determine the parameter regions in which the stripe configuration is
stable, as well as the boundary of the SVI.

\subsection{Linearisation}

We proceed by considering perturbations to the basic stripe solution. 
We suppose that the vorticity perturbation contains wavevectors $(k,l)$ and 
$(-k,-l)$. These interact with the wavevectors $(1+q,0)$ and $(-1-q,0)$ in the 
stripe solution to give four new wavevectors, and so the perturbed solution 
can be written as $\psi=\psi_{0}+\psi^{'}$ and $\omega=\omega_{0}+\omega^{'}$, with $\omega_0=0$ and 
 \begin{eqnarray}
 \label{eq:perturbpsi}
 \psi^{'}=A(t)e^{i(1+q+k,l)\cdot\mathbf{x}}+ B(t)e^{i(-1-q+k,l)\cdot\mathbf{x}}\,\,\,\,\,\,\,\,\,\,\, \nonumber \\
 \quad{}\,\,\,\,\,+\bar{A}(t)e^{-i(1+q+k,l)\cdot\mathbf{x}}+ \bar{B}(t)e^{i(1+q-k,-l)\cdot\mathbf{x}},\\
 \label{eq:perturbomega}
 \omega^{'}=C(t)e^{i(k,l)\cdot\mathbf{x}}+ \bar{C}(t)e^{i(-k,-l)\cdot\mathbf{x}}.
  \end{eqnarray} We can calculate the stream function $\zeta$ from $\omega(x,y,t)=-\nabla^2\zeta$ by inverting the Laplacian, and hence obtain the mean flow:
 \[\mathbf{U}=\frac{i}{k^{2}+l^{2}}\left(C(t)e^{i(k,l)\cdot\mathbf{x}}-\bar{C}(t)e^{-i(k,l)\cdot\mathbf{x}}\right)(l,-k)\;.\]

We substitute the expression above into the two models
and linearize (assuming that $A$, $B$ and $C$ are small). Examining the
coefficients of $e^{i(1+q+k,l)\cdot\mathbf{x}}$  and
$e^{i(-1-q+k,l)\cdot\mathbf{x}}$ results in linear ODEs for $A$ and~$B$:
 \begin{eqnarray}
 \label{eq:perturbA}
 \dot{A}=\left(\mu-\left[1-\left(\left(1+q+k\right)^{2}+l^{2}\right)\right]^{2}-6\beta\right)A\nonumber\\
          -3\beta B+\frac{l(1+q)\sqrt{\beta}}{k^{2}+l^{2}}C,\\
 \dot{B}=\left(\mu-\left[1-\left(\left(-1-q+k\right)^{2}+l^{2}\right)\right]^{2}-6\beta\right)B\nonumber\\
          -3\beta A-\frac{l(1+q)\sqrt{\beta}}{k^{2}+l^{2}}C\;.
 \label{eq:perturbB}
 \end{eqnarray}
The term $k^{2}+l^{2}$, which appears in the denominator of the governing
equations for $\dot{A}$ and $\dot{B}$, arises from inverting the Laplacian when
calculating~$\zeta$ and hence~$\mathbf{U}$ in the linearisation of the 
$\left(\mathbf{U}\cdot\mathbf{\nabla}\right)\psi$ term.
The equation for~$C$ differs between the two models. In model~1, the
coefficient of $e^{i(k,l)\cdot\mathbf{x}}$ yields
\begin{widetext} 
 \begin{equation}
 \dot{C}=-Pr\left(k^{2}+l^{2}+c^2\right)C
         +\mathcal{F}_{\gamma}g_ml(1+q)\sqrt{\beta}\left(\left[-\left(1+q-k\right)^{2}-l^{2}+(1+q)^2\right]B+\left[\left(1+q+k\right)^{2}+l^{2}-(1+q)^2\right]A\right).
 \label{eq:perturbCmodel1}
 \end{equation} \normalsize In model~2, with no intrinsic dynamics for~$\omega$, 
there is an algebraic relation between $C$, $A$ and~$B$:
 \begin{equation}
 C=  \mathcal{F}_{\gamma}gl(1+q)\sqrt{\beta}
                 \left(       ( (1+q+k)^{2}+l^{2}-(1+q)^2)A
                            + (-(1+q-k)^{2}-l^{2}+(1+q)^2)B\right).
 \label{eq:perturbCmodel2}
 \end{equation}
 \end{widetext} In Fourier space, the effect of the filtering, $\mathcal{F}_{\gamma}$, is to
reduce the amplitude of a Fourier component; for the function
$e^{i(k,l)\cdot\mathbf{x}}$, $\mathcal{F}_{\gamma}$~reduces the amplitude by
the multiplicative factor~$e^{-\gamma^2(k^{2}+l^{2})}$.

Equations (\ref{eq:perturbA}--\ref{eq:perturbCmodel1}) for model 1 can be succinctly
expressed as:
 \[
 \left( \begin{array}{c} \dot{A} \\ \dot{B} \\ \dot{C} \end{array}\right) 
 = \begin{bmatrix} M_{1} & M_{2} &  M_{3} \\
                   M_{2} & M_{4} & -M_{3} \\
                   g_m M_{5} & g_m M_{6} & M_{7}
   \end{bmatrix}
 \left( \begin{array}{c} A \\ B \\ C \end{array} \right)
 = J_1 \left( \begin{array}{c} A \\ B \\ C \end{array} \right).
 \]
Equations (\ref{eq:perturbA}--\ref{eq:perturbB}) and (\ref{eq:perturbCmodel2}) for model 2
yield:
 \begin{align*}
 \left( \begin{array}{c} \dot{A} \\ \dot{B}  \end{array} \right)
 = \begin{bmatrix} M_{1}+gM_{3}M_{5}  & M_{2}+gM_{3}M_{6} \\ 
                   M_{2}-gM_{3}M_{5}  & M_{4}-gM_{3}M_{6} 
   \end{bmatrix} 
 \left( \begin{array}{c} A \\ B  \end{array} \right)\\
 = J_2 \left( \begin{array}{c} A \\ B  \end{array} \right).
\end{align*} Here, we use the abbreviations:
 \begin{align*}
 M_{1}&=\mu-[1-((1+q+k)^{2}+l^{2})]^{2}-6\beta, \\ 
 M_{2}&=-3\beta, \\ 
 M_{3}&=-l(1+q)\sqrt{\beta}/(k^{2}+l^{2}), 
 \end{align*}
 \begin{align*}
 M_{4}&=\mu-[1-((-1-q+k)^{2}+l^{2})]^{2}-6\beta,\\ 
 M_{5}&=e^{-\gamma^{2}(k^{2}+l^{2})}l(1+q)\sqrt{\beta}[2k(1+q)+k^{2}+l^{2}], \\ 
 M_{6}&=e^{-\gamma^{2}(k^{2}+l^{2})}l(1+q)\sqrt{\beta}[2k(1+q)-k^{2}-l^{2}], \\
 M_{7}&=-Pr(k^{2}+l^{2}+c^2).
 \end{align*}We note that in the limit $(k,l)\rightarrow(0,0)$, we have $M_{7}\approx-Pr\,c^{2}$, so it is not surprising (looking at the bottom line of the 3$\times3$ matrix for model 1) that long-wavelength instabilities in model 1 will depend only on the combination $g_{m}/(Pr\,c^{2})$.

\subsection{Determinants in the limit of small $k$ and $l$}

The characteristic polynomials (and hence the eigenvalues, traces and 
determinants) of each of
these Jacobian matrices are even in $k$ and~$l$. Bifurcations occur when an 
eigenvalue crosses through zero. The stripe solution is stable only if all eigenvalues are negative for all $(k,l)$, so we are interested in extreme values of the eigenvalues as functions of $k$ and $l$. It can be readily checked that a zero extreme value of the eigenvalue corresponds to a zero extreme value of the determinant. 
Consequently, we use the determinants of $J_1$ and~$J_2$  to assist our analysis of instabilities. The determinant of~$J_1$, $Det(J_1)$, is: 
\begin{widetext} \small \[
 \frac{
   \left(P_{1}^{(1)}k^{4}+P_{2}^{(1)}k^{2}l^{2}+P_{3}^{(1)}l^{4}\right)
 + \left(Q_{1}^{(1)}k^{6}+ \dots
       + Q_{4}^{(1)}l^{6}\right)
 + \left(R_{1}^{(1)}k^{8}+ \dots
       + R_{5}^{(1)}l^{8}\right)
 + \left(S_{1}^{(1)}k^{10} + \dots
       + S_{6}^{(1)}l^{10}\right)
 - Pr(k^{2}+l^{2})^{6}}
 {k^{2}+l^{2}},
 \] \normalsize where all coefficients $P_{i}^{(1)}$ \hbox{etc.} are functions of $\mu$, $q$,
$Pr$, $g_m$, $c$ and the filtering~$e^{-\gamma^2(k^{2}+l^{2})}$.
The determinant of~$J_2$, $Det(J_2)$, is:
 \[
 \frac{
   \left(P_{1}^{(2)}k^{4}+P_{2}^{(2)}k^{2}l^{2}+P_{3}^{(2)}l^{4}\right)
 + \left(Q_{1}^{(2)}k^{6}+Q_{2}^{(2)}k^{4}l^{2}+Q_{3}^{(2)}k^{2}l^{4}
        +Q_{4}^{(2)}l^{6}\right)
 + \left(R_{1}^{(2)}k^{8}+\dots
        +R_{5}^{(2)}l^{8}\right)
 +(k^{2}+l^{2})^{5}}
 {k^{2}+l^{2}},
 \]\end{widetext} 
\normalsize where all coefficients $P_{i}^{(2)}$ \hbox{etc.} are functions of $\mu$, $q$, $g$ and the filtering~$e^{-\gamma^2(k^{2}+l^{2})}$.
The traces of the two Jacobians can be written as
 \begin{eqnarray*}
 Tr(J_{1})=
 -6\beta-Pr\,c^{2}
 +\left(4-Pr-12(1+q)^2\right)k^{2}\\
 +\left(4-Pr-4(1+q)^2\right)l^{2}
 -2\left(k^{2}+l^{2}\right)^{2}
 \end{eqnarray*}
and
 \begin{eqnarray*}
 Tr(J_{2}) = 
 -6\beta
 +\left(4-12(1+q)^{2}\right)k^{2}
 +\left(4-4(1+q)^2\right.\\\left.-2\beta{g}(1+q)^2e^{-\gamma^2(k^{2}+l^{2})}\right)l^{2}
 -2\left(k^{2}+l^{2}\right)^{2},
 \end{eqnarray*} where we recall the relationship in (\ref{eq:nonlinearstripes}) between $\beta$ 
and~$\mu$.

Note that at this point, no approximations or truncations have been made in the
linear stability problem of the stripe solution, by virtue of having an exact
solution. Our task is now to work out the most unstable eigenvalues in the
limit of small $k$ and~$l$; this is made more challenging by the presence of
$k^2+l^2$ in the denominators of the determinants above.

We note that explicit expressions for eigenvalues of $J_{1}$ are not, in
general, analytically attainable (though the eigenvalues can be calculated
numerically). For the matrix $J_{2}$, in the limit
$(k,l)\rightarrow(0,0)$, the trace is $Tr(J_2)=-6\beta$ and the determinant is
zero, so, for small~$(k,l)$, one eigenvalue will be
$-6\beta+\mathcal{O}(k^2+l^2)$, which is bounded away from zero for a
finite-amplitude stripe. The other eigenvalue will be close to zero, approximately $Det(J_2)/Tr(J_2)$. Similarly, in the limit $(k,l)\rightarrow(0,0)$, $Det(J_1)=0$, so $J_1$ will have an eigenvalue close to zero for small $(k,l)$. Since bifurcations occur when an eigenvalue is equal to zero, this can be detected in both cases by considering only the determinants of the two matrices. Hopf bifurcations (see section~(\ref{osv})) require additional consideration.
We will expand $Det(J_1)$ and $Det(J_2)$ in powers of $k$ and $l$, including the  filtering $e^{-\gamma^2(k^2+l^2)}$ in the expansion. This yields expressions of the form,
 \begin{widetext}
 \begin{equation}
 Det(J_{1,2})=
 \frac{\left(A_{1,2}k^{4}+B_{1,2}k^{2}l^{2}+C_{1,2}l^{4}\right) 
      +\left(D_{1,2}k^{6}+E_{1,2}k^{4}l^{2}+F_{1,2}k^{2}l^{4}+G_{1,2}l^{6}\right)
      +\mathcal{O}((k^2+l^2)^4)}
      {k^{2}+l^{2}},
 \label{eq:sig2}
 \end{equation}
 
 \end{widetext} where in model 2, the coefficients are:
 \begin{align*}
  A_{2}&=12\beta(3(1+q)^2-1) -16q^2(1+q)^2(2+q)^2, 
   \end{align*}
  \begin{align*} 
  B_{2}&=-24\beta(1-2(1+q)^2) -16q^2(1+q)^2(2+q)^2 \\
  &\quad{}+4(1+q)^2\beta\left(3\beta-4q(1+q)^2(2+q)\right)g,\\ 
   C_{2}&=12\beta\left(\beta g(1+q)^2 + q(2+q)\right), 
  \end{align*}
  \begin{align*}  
  D_{2}&=6\beta +4+4(1+q)^2\left((1+q)^2+2\right), \\
  E_{2}&=g\left[4\gamma^2(1+q)^2 \beta\left( 4(1+q)^4-3\beta-4(1+q)^2\right)\right.\\
  &\left.\quad{}-4\beta(1+q)^2 \left((1+q)^2 +1 \right)\right]+12+8(1+q)^2\\
  &\quad{}+18\beta-4(1+q)^4,   
  \end{align*} and in model 1, the coefficients are:
  \begin{align*}
 A_{1}&=-Pr\,c^2A_{2}, \\ 
  B_{1}&=-Pr\,c^2B_{2}, \\ 
  C_{1}&=-Pr\,c^2C_{2}, \\    
  D_{1}&=-Pr\,c^2D_{2}-PrA_{2}, \\ 
  E_{1}&=-Pr\,c^2E_{2}+\beta(-48-84q^2-168q)+32q^2\\
  &\quad{}(q+1)^2(q+2)^2.
 \end{align*} We will hence  refer to the determinants of two matrices in the general form with no subscripts:
\begin{align}
 Det(J)=
 \left[\left(Ak^{4}+Bk^{2}l^{2}+Cl^{4}\right) 
      +\left(Dk^{6}+Ek^{4}l^{2}+Fk^{2}l^{4}\right.\right.\nonumber\\\left.\left.+Gl^{6}\right)
      +\mathcal{O}((k^2+l^2)^4)\right]/
      (k^{2}+l^{2}).
 \label{eq:sig2}
 \end{align} The values of $F$ and $G$ are not needed subsequently. 
Stripes are stable if all eigenvalues are less than zero, corresponding to  $Det(J_{1})<0$ and $Det(J_{2})>0$. To simplify the presentation, we  focus on $Det(J_{1})$ since the sign of $Det(J_{1})$ coincides with the sign of the most unstable eigenvalue.


\section{LONG-WAVELENGTH INSTABILITIES: Zigzag and Eckhaus }
\begin{figure}[t]
\begin{center}
\scalebox{0.38}{\includegraphics*[viewport=0.48in 0.6in 10in 6.34in]{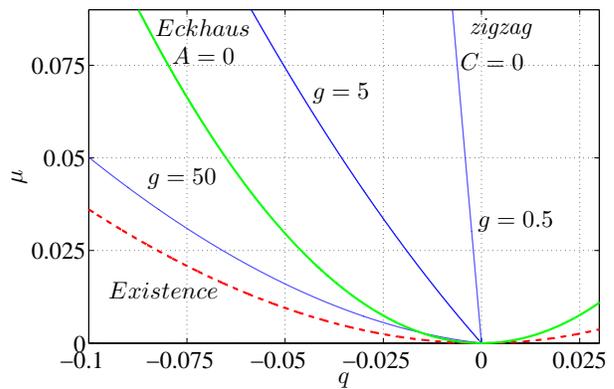}}
\end{center}
\vspace{-0.15in}
\caption{(Color online) Location of the Eckhaus and zigzag stability 
boundaries in 
the $(q,\mu)$ plane ($q<0$),
for $g=0.5, 5$ and $ 50$. The Eckhaus boundary derived from $A=0$, is denoted in green thick curve. The dashed red curve is the existence boundary. The 
zigzag boundary (blue thin curves) crosses the Eckhaus boundary at $q\approx-0.015$  when  $g=50$. The
behavior for small $q$ is approximately
$\mu=-6q/g$. Stripes are stable to the right of the Eckhaus and zigzag 
boundaries.}
 \label{fig:2}
\end{figure}
Bifurcation points correspond to parameter values for which an eigenvalue has zero
real part. We first investigate how the determinants depend on $k$ and~$l$, and then use this information to
explore how the bifurcation lines depend on the other parameters $\mu$, $q$ and either $g$ or $g_{m}$, $Pr$ and $c$. In this section we examine the well known Eckhaus and zigzag
instabilities, which correspond to perturbations with $l=0$ and $k=0$
respectively.

In the Eckhaus case, with $l=0$, (\ref{eq:sig2}) yields
$Det(J)=Ak^{2}+Dk^{4}+\ldots$. For small~$k$, this is positive  when $A>0$. Thus for model 1, instability corresponds to  $A_1>0$ and for model 2, instability corresponds to $A_2<0$.

Accordingly, in both cases, the Eckhaus instability, which is independent of the mean flow,
occurs  when $A=0$:
\[\mu_{Eck}=\frac{q^{2}(7q^{4}+42q^{3}+90q^{2}+80q+24)}{(3q^{2}+6q+2)}.\]
 \begin{flushleft}
Note that in the limit of  $q\rightarrow0$, $\mu_{Eck}\rightarrow 12q^{2}$ and hence
$\mu_{Eck}\rightarrow3\mu_{Existence}$. 
\end{flushleft}

Similarly, in the zigzag case, with $k=0$, (\ref{eq:sig2}) yields
$Det(J)=Cl^{2}+Gl^{4}+\dots$, which is positive  for small~$l$ when $C>0$. Thus for model 1, instability corresponds to  $C_1>0$ and for model 2, instability corresponds to $C_2<0$.
Accordingly, the zigzag instability in both models
occurs when $C=0$:
 \begin{align*}
 \mu_{zigzag}&=\mu_{Existence}-\frac{3q(2+q)}{g(1+q)^{2}},
  \end{align*} where for model 1 we have identified $g=g_m/Pr\,c^2$.
 
Unlike the Eckhaus instability, the zigzag instability is affected by the mean flow. Vorticity and mean flows act as a stabilizing influence on the zigzag
instability, which is suppressed for larger values of~$g$, resulting in a
larger region of stable stripes for $q<0$ in the $(\mu,q)$ stability diagram.
Figure \ref{fig:2} shows how the zigzag instability boundary behaves for
different values of~$g$. Note that for large enough $g$, it no longer forms the lower stability
boundary except for very small $\mu$.

The  zigzag and Eckhaus instabilities cross in parameter space when $g$ and $q$ are related by 
 \[
 g=\frac{3(3(1+q)^2-1)}{4q(1+q)^4(2+q)},
 \]
where $q<0$. Furthermore, when $g\rightarrow0$, we recover the standard result
for the SHE without mean flow, and when
$g\rightarrow\infty$, we have $\mu_{zigzag}\rightarrow\mu_{Existence}$. However,
even in this limit, for small $\mu$ and $q$, $\mu_{zigzag}$ is
approximately~$-6q/g$, so the zigzag stability boundary emerges from
$(q,\mu)=(0,0)$ in a straight line of slope $-6/g$, as can be seen in
figure~\ref{fig:2}.

\section{LONG-WAVELENGTH INSTABILITIES: Skew-varicose}

The skew-varicose instability, which is driven by  the inclusion of mean flow, the strength of which is determined by $g$, is associated with modes for which the maximum positive growth rate occurs when $k\neq0$ and $l\neq0$. Two conditions are required to characterize the SVI: the determinant should be zero and should have maximum or minimum value (for model 1 and model 2 respectively) for $k\neq0$ and $l\neq0$.
We first express these conditions for the SVI in terms of the coefficients $A-G$ of the expression (\ref{eq:sig2}), the power series expansion of the determinants of $J_{1}$ and $J_{2}$, which we denote simply by $Det$. We then  express these conditions in terms of the parameters $\mu$, $q$ and either $g$ or $g_{m}$, $c$ and $Pr$, in order to locate the  SVI boundary in the $(\mu,q)$ plane.

\begin{figure}[t]
\mbox{
\subfigure[]{
\scalebox{0.37}{\includegraphics*[viewport=0in 0.6in 9.25in 6.45in]{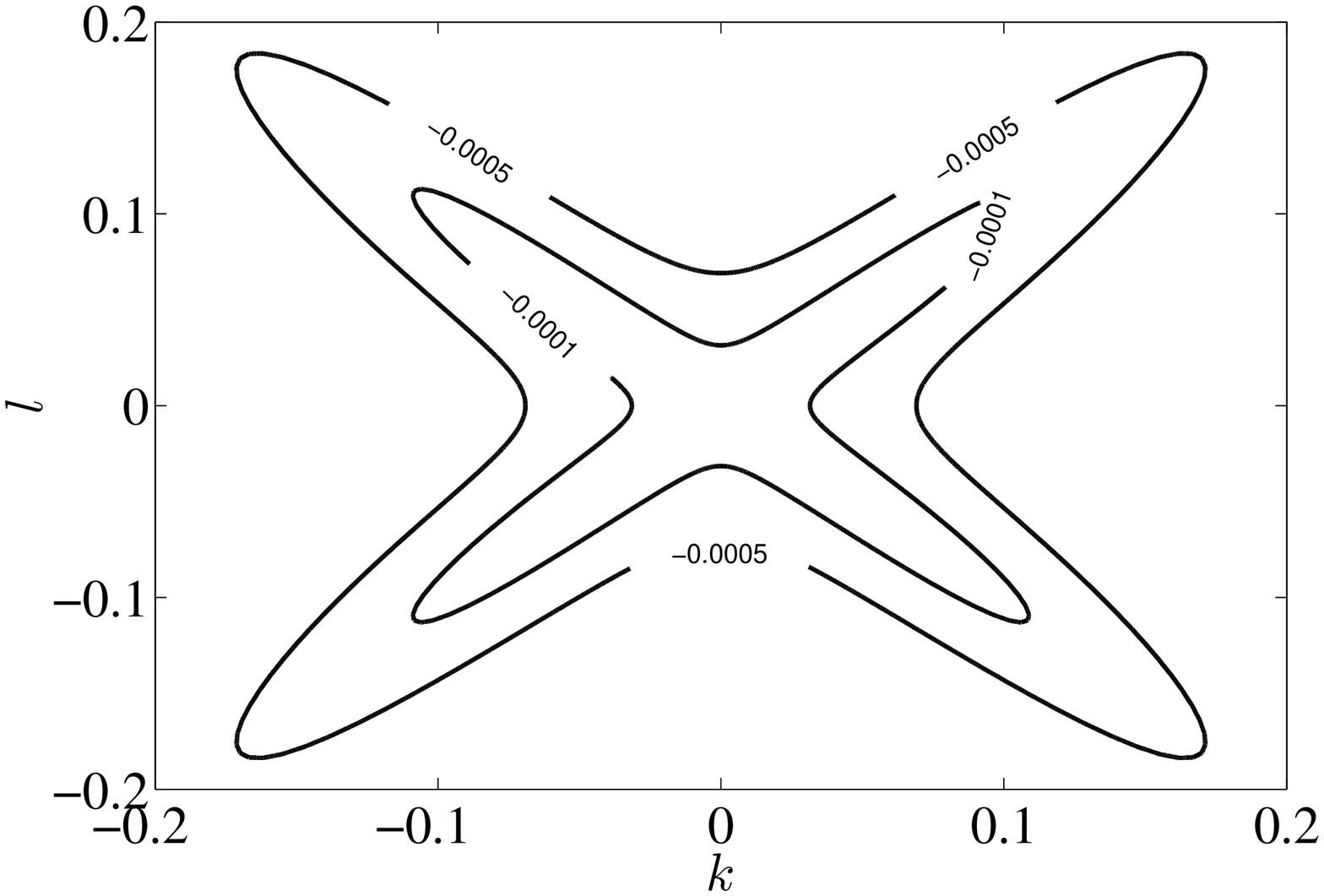}}
\label{fig:3a}}
}
\mbox{
\subfigure[]{
\scalebox{0.37}{\includegraphics*[viewport=0in 0.6in 9.25in 6.45in]{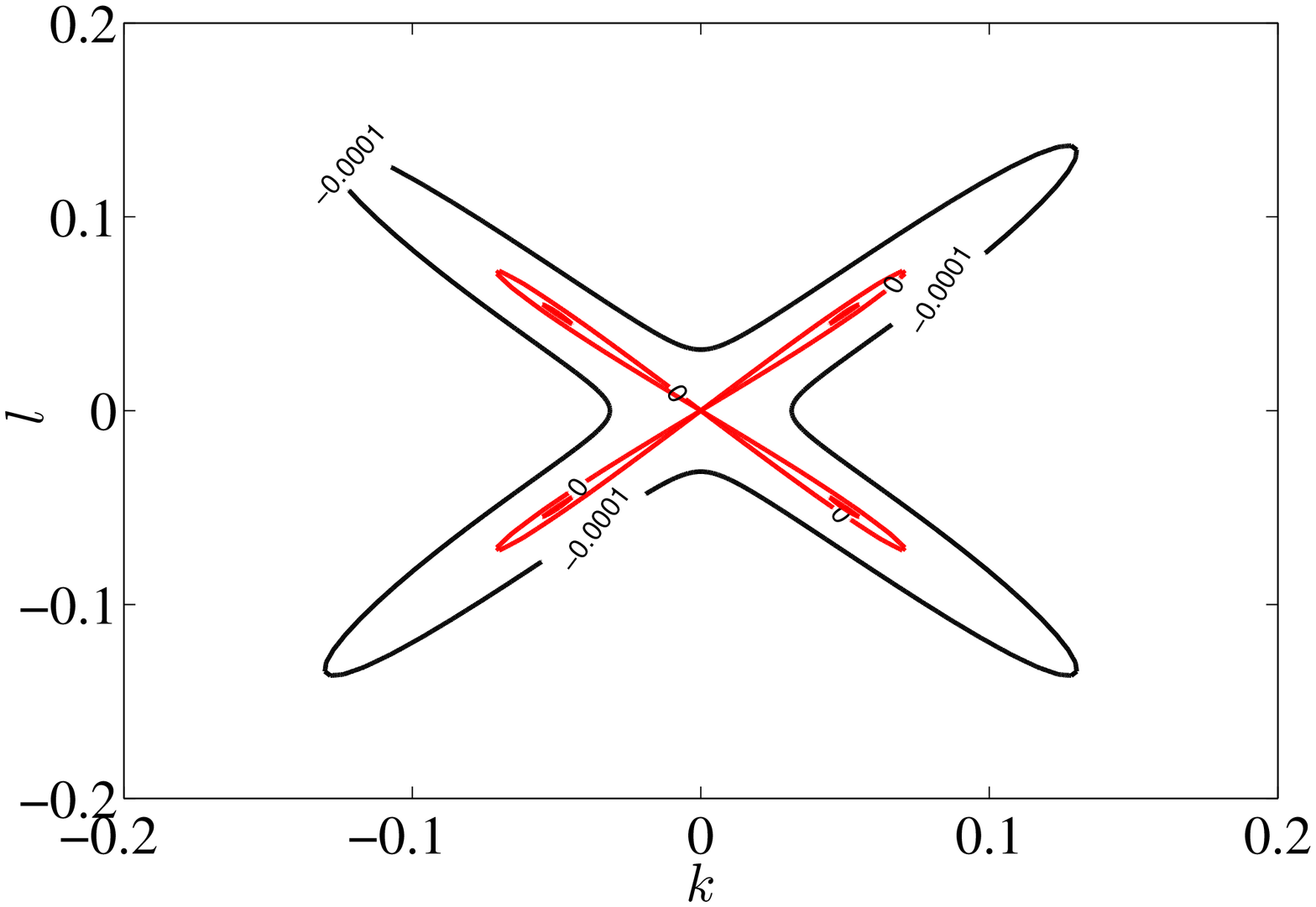}}
\label{fig:32}}
}
\caption{(Color online) Case I of the SVI: behavior of the $Det$ in the $(k,l)$ plane for $A=-0.1$ and $C=-0.1$ ($A<0$ and $C<0$ makes stripes Eckhaus and zigzag stable). The other coefficients are $D=-1$, $E=2$, $F=-1$ and $G=-1$. (a) $B=0.195$, giving stable stripes $(B^{2}<4AC)$. (b) $B=0.205$, giving stripes that are unstable to the SVI $(B^{2}>4AC)$. The positive maximum of the determinant emerges from  $(k,l)=(0,0)$ but occurs with  $k\neq0$ and $l\neq0$. A negative value of the determinant is indicated by black contours while zero and positive values of the determinant are in red(gray). }
 \label{fig:3}
\end{figure}

There are two different manifestations of the skew-varicose
instability.  In the first (case I), the instability emerges from
$k=l=0$, as illustrated in figure
\ref{fig:3}. In this case, the determinant is negative for $(k,l)$ close to $(0,0)$, corresponding to $A<0$ and $C<0$. If we suppose in the first instance that we can write the leading order terms in the determinant as:
$Det=\frac{Ak^{4}+Bk^{2}l^{2}+Cl^{4}+\mathcal{O}\left((k^{2}+l^{2})^{3}\right)}{k^{2}+l^{2}}$,
then imposing $Det=0$, along with $\partial Det/\partial k=0$ and $\partial Det/\partial l=0$, results in 
$(B^{2}-4AC)(A-B+C)=0$ and $\frac{k^{2}}{l^{2}}=\frac{B-2C}{B-2A}$. In order for ${k^{2}}/{l^{2}}$ to be positive, we need $B>max(2A,2C)$, which excludes  $B=A+C$ and $B=-\sqrt{4AC}$. We conclude that
\begin{equation}B^{2}-4AC>0,\ A<0,\  C<0\ \text{and}\ B>0\,\text{with}\,\,  \frac{k^{2}}{l^{2}}=\sqrt{\frac{C}{A}} \label{Eq:1}\end{equation}is the condition for the \hbox{SVI}. However, the truncation of $Det$ above is degenerate: the conditions are satisfied along a line in the $(k,l)$ plane, rather than at a point. This degeneracy is resolved by restoring the higher order terms, as illustrated in figure \ref{fig:3}.

\begin{figure}[t]
\mbox{
\hspace{-0.8in}
\subfigure[]{
\scalebox{0.563}{\includegraphics*[viewport=0.4in 0in 6.2in 3.8in]{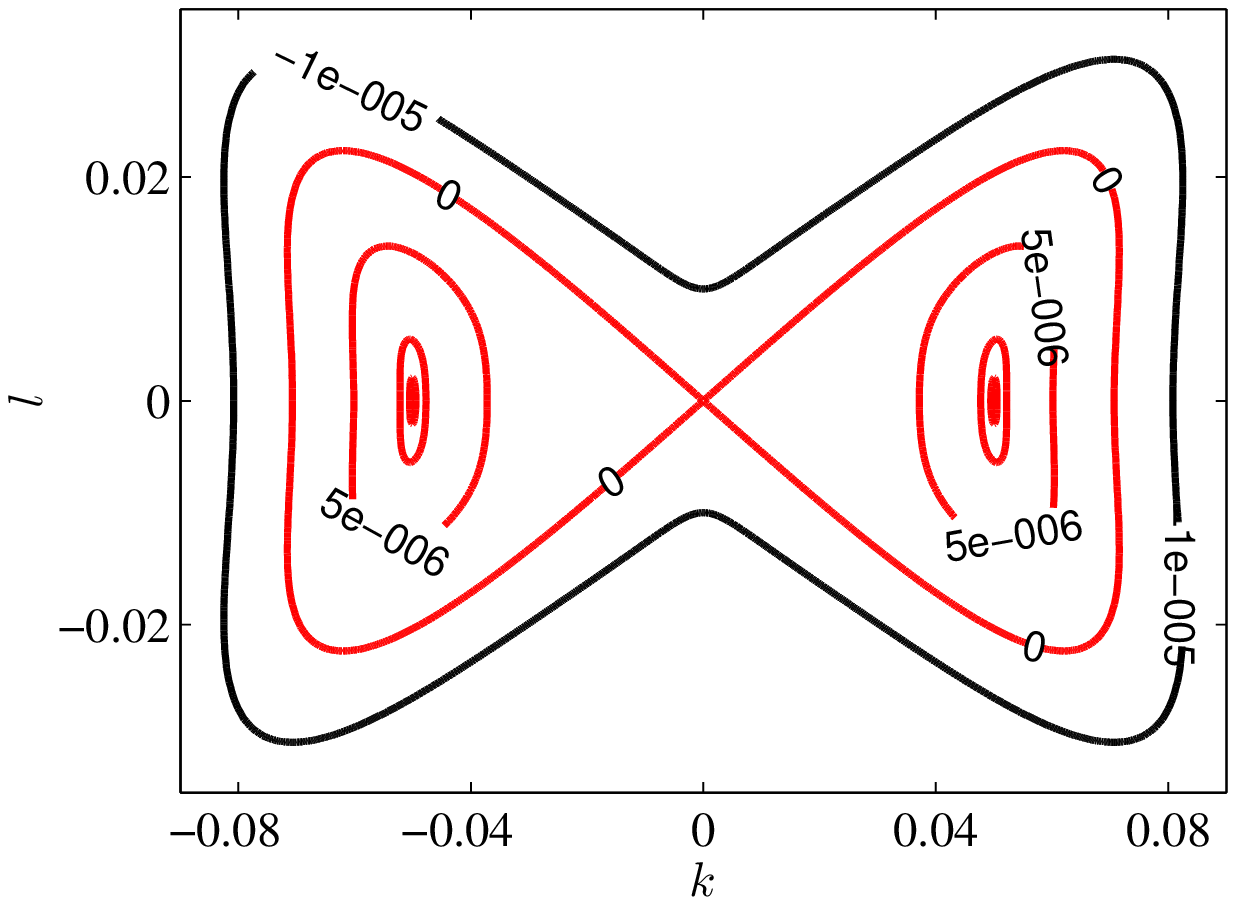}}

\label{fig:41}}
}
\mbox{
\hspace{-0.5in}
\subfigure[]{
\scalebox{0.58}{\includegraphics*[viewport=1in 0in 6.2in 3.8in]{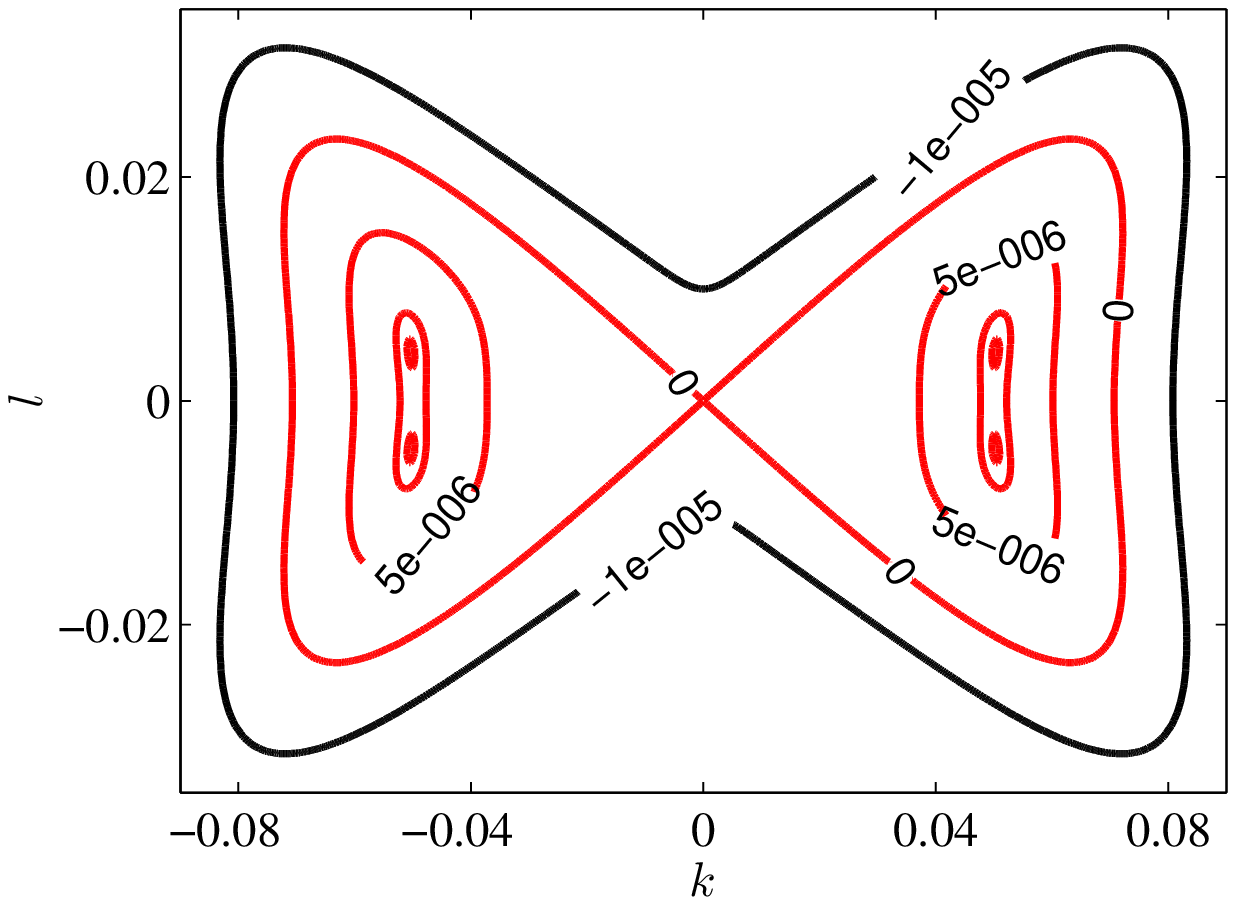}}

\label{fig:42}}
}
\caption{(Color online) Case II of the SVI: behavior of the $Det$ in the $(k,l)$ plane for $A=0.005$ and $C=-0.1$ ($A>0$ and $C<0$ makes stripes  Eckhaus unstable but stable to zigzags). The other coefficients are $D=-1$, $E=2$, $F=-1$ and $G=-1$. (a) $B=-0.003$, giving SV stable and Eckhaus unstable stripes $\left(B<\frac{D+E}{2D}A\right)$; the maximum occurs with $l=0$. (b) $B=-0.001$, giving stripes that are unstable to both SV and Eckhaus instabilities $\left(B>\frac{D+E}{2D}A\right)$; the maximum moves off axis, and $(k_{max},0)$ is now a saddle. A negative value of the determinant is indicated by black contours while zero and positive values of the Determinant are in red (gray).}
 \label{fig:4}
\end{figure}

\begin{figure}[t]
\hspace{-0.3in}
\scalebox{0.37}{\includegraphics*[viewport=0in 0.6in 9.25in 6.45in]{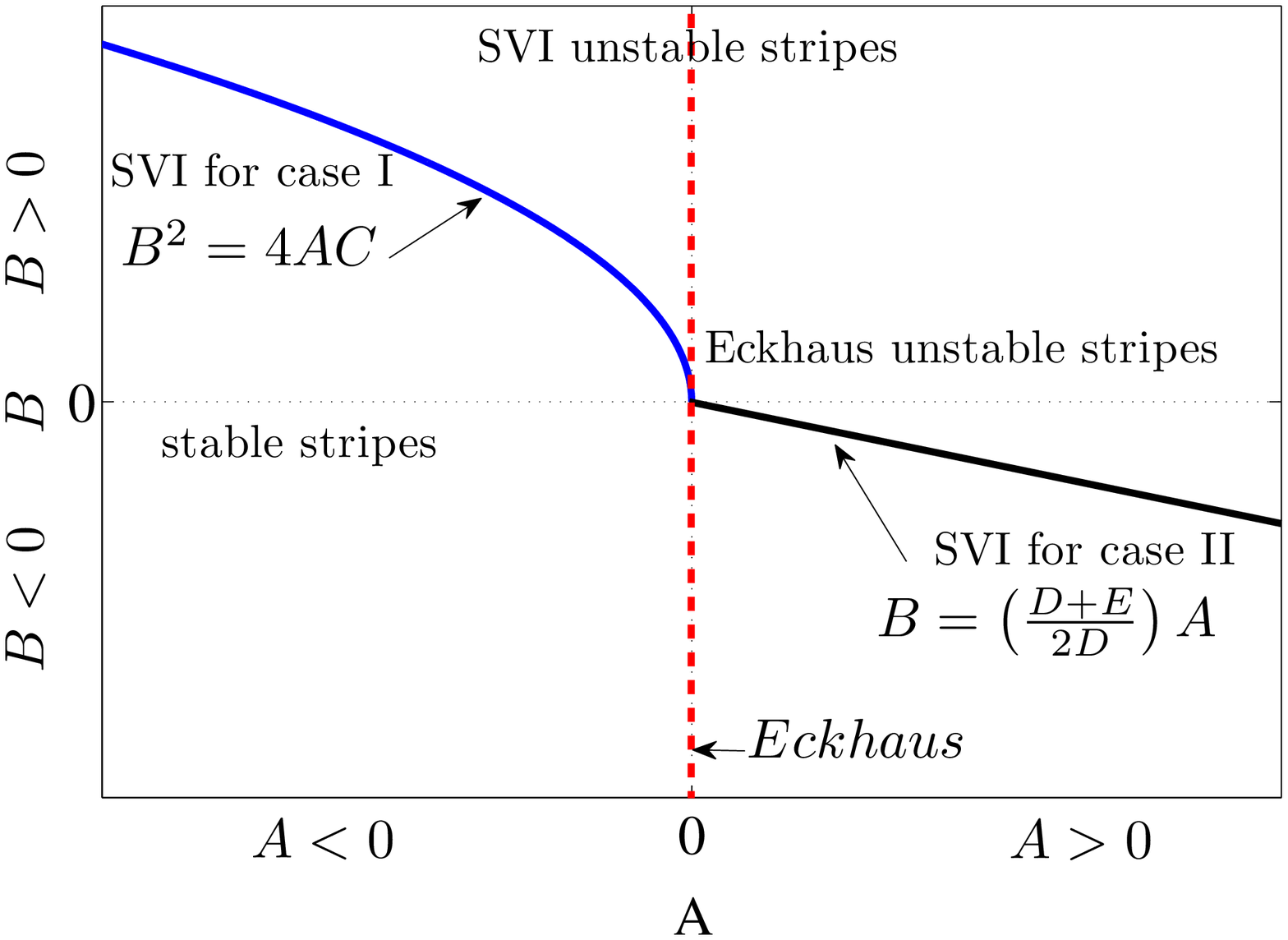}}
\caption{(Color online) Schematic diagram of the conditions for the Eckhaus and SV instabilities in cases I and II in $(A,B)$ plane. The Eckhaus instability occurs when $A=0$ (dashed red line). The blue curve shows the SV stability boundary: $B^{2}=4AC$ in case I ($A<0$), and $B=\left(\frac{D+E}{2D}A\right)$ in case II ($A>0$). At the intersection of the two stability boundaries, $(A,B)=(0,0)$.}
 \label{fig:5}
\end{figure}

The other possibility for the skew-varicose instability is that it can
accompany the Eckhaus instability. This occurs
when \emph{$A>0$} and \emph{$C<0$}: for $l=0$ and $D<0$, $Det$ is positive for a range
of $k$ and attains its maximum on the $l=0$ axis at a finite 
$k=k_{max}$. In the $(k,l)$ plane, $(k_{max},0)$ can either be a maximum or a saddle, as illustrated in figure \ref{fig:4}. We define the SVI (case II) to be the point at which $(k_{max},0)$ changes from a maximum to a saddle; at this point the maximum eigenvalue moves off the $k$ axis.

Unlike in the previous case, the growth rate at the SVI is positive, since stripes are already Eckhaus unstable. Therefore the 
SVI occurs when there is a degenerate maximum at $l=0$ and
$k=k_{max}$, about to become saddle. The conditions
$\frac{\partial Det}{\partial k^{2}}=0$  and
$\frac{\partial Det}{\partial l^{2}}=0$ at this point yield the parameter values at which this variant of the
skew-varicose instability occurs. 

At $l=0$, the first condition implies
${k^{2}_{max}}=-\frac{A}{2D}+\mathcal{O}(A^2) $, and the second condition implies
$(B-A)+(E-D)\left(-\frac{A}{2D}\right)=0$ at $k=k_{max}$. Hence in the limit of small $A$, the condition for case II of the SVI is 

\begin{eqnarray}B=\left(\frac{D+E}{2D}\right)A+\mathcal{O}(A^{2}),\, A>0 \, \text{and}\,  D<0,\,\, \nonumber\\\text{with}\,\, k^{2}=\frac{-A}{2D}+\mathcal{O}(A^{2})\, \text{and}\,\, l=0\, .\label{Eq:2}\end{eqnarray} Contours of $Det$ on either side of this boundary are illustrated in figure \ref{fig:4}. Moreover, the conditions (\ref{Eq:1}) and (\ref{Eq:2}) are summarized in  the schematic diagram in figure~\ref{fig:5}, which also indicates the regions affected by the SVI and the point $(A,B)=(0,0)$, where the two cases coincide.

 \begin{figure}[t]
\scalebox{0.37}{\includegraphics*[viewport=0.5in 0.6in 9.25in 6.45in]{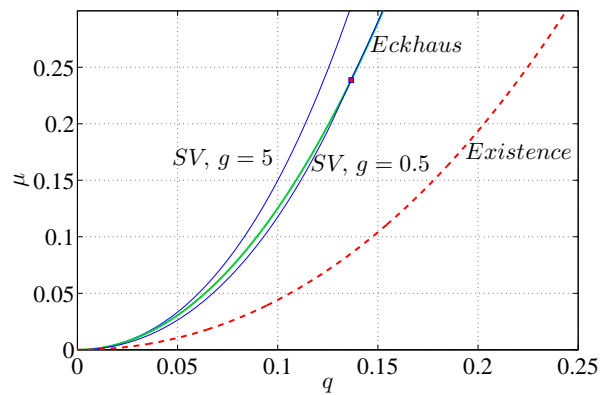}}
\caption{(Color online) Numerical computation (in model 2) of the SV stability boundary in the $(\mu,q)$ plane. For $g=5$, the SVI pre-empts the Eckhaus instability for all $\mu$ and hence the region of stable stripes is bounded by the 
skew-varicose instability curve. However, for $g=0.5$,
the region of stripe stability is bounded by the skew-varicose instability curve only when $\mu>0.2387$.
The crossing point, $(0.1367,0.2387)$, of the two boundaries is denoted as a  red square.
For $0<\mu<0.2387$, the Eckhaus precedes the SV curve, which reaches  the origin as a parabola $\mu=9.33q^{2}$. 
The green thick curve denoted by E is for the Eckhaus boundary whereas the dashed red curve is the boundary of existence of
stripes.}
 \label{fig:6}
\end{figure}

\begin{figure}[t]
\mbox{
\subfigure[]{
\scalebox{0.37}{\includegraphics*[viewport=0.3in 0.6in 9.45in 6.45in]{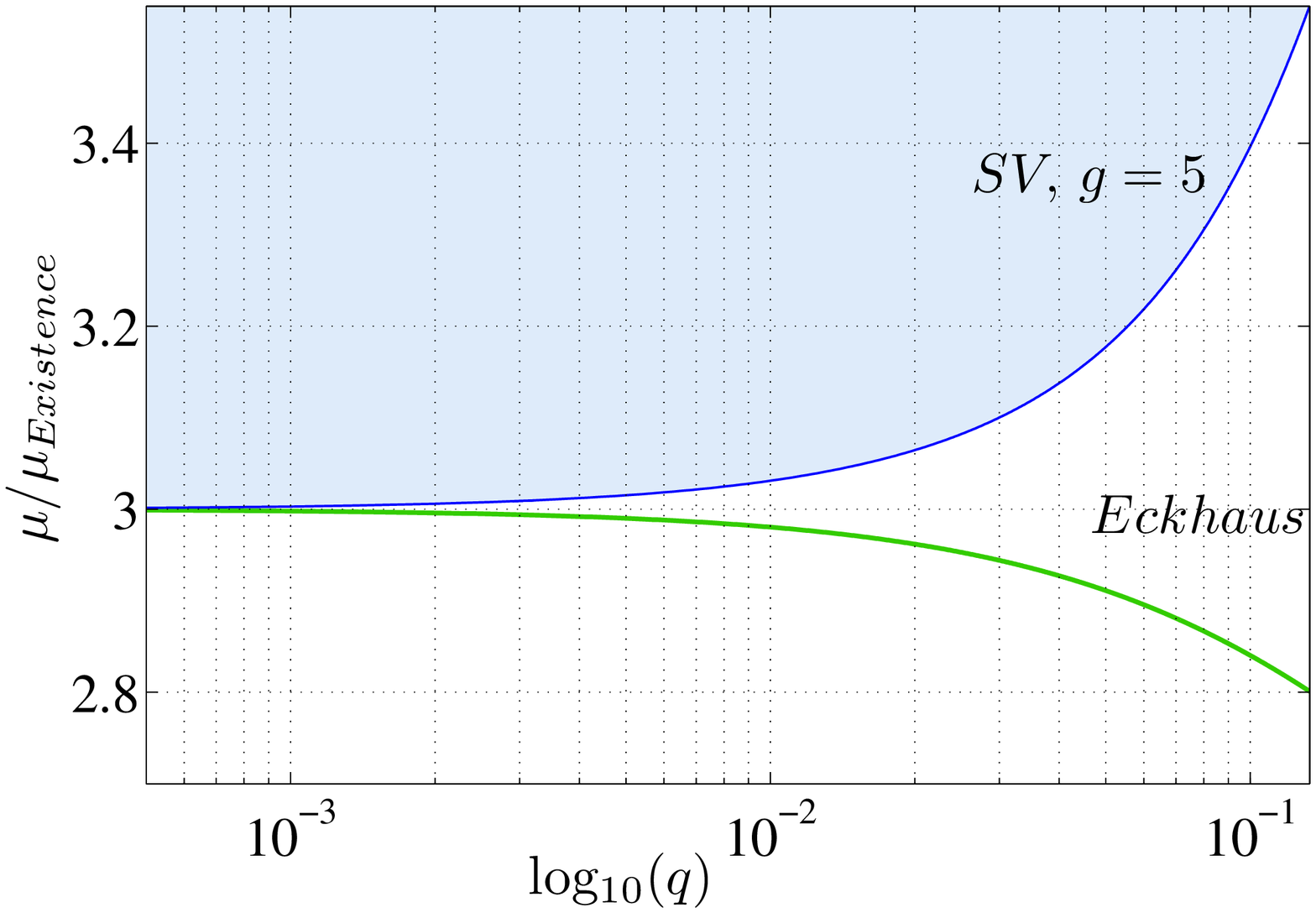}}
\label{fig:7a}}}
\mbox{
\subfigure[]{
\scalebox{0.37}{\includegraphics*[viewport=0.3in 0.6in 9.45in 6.45in]{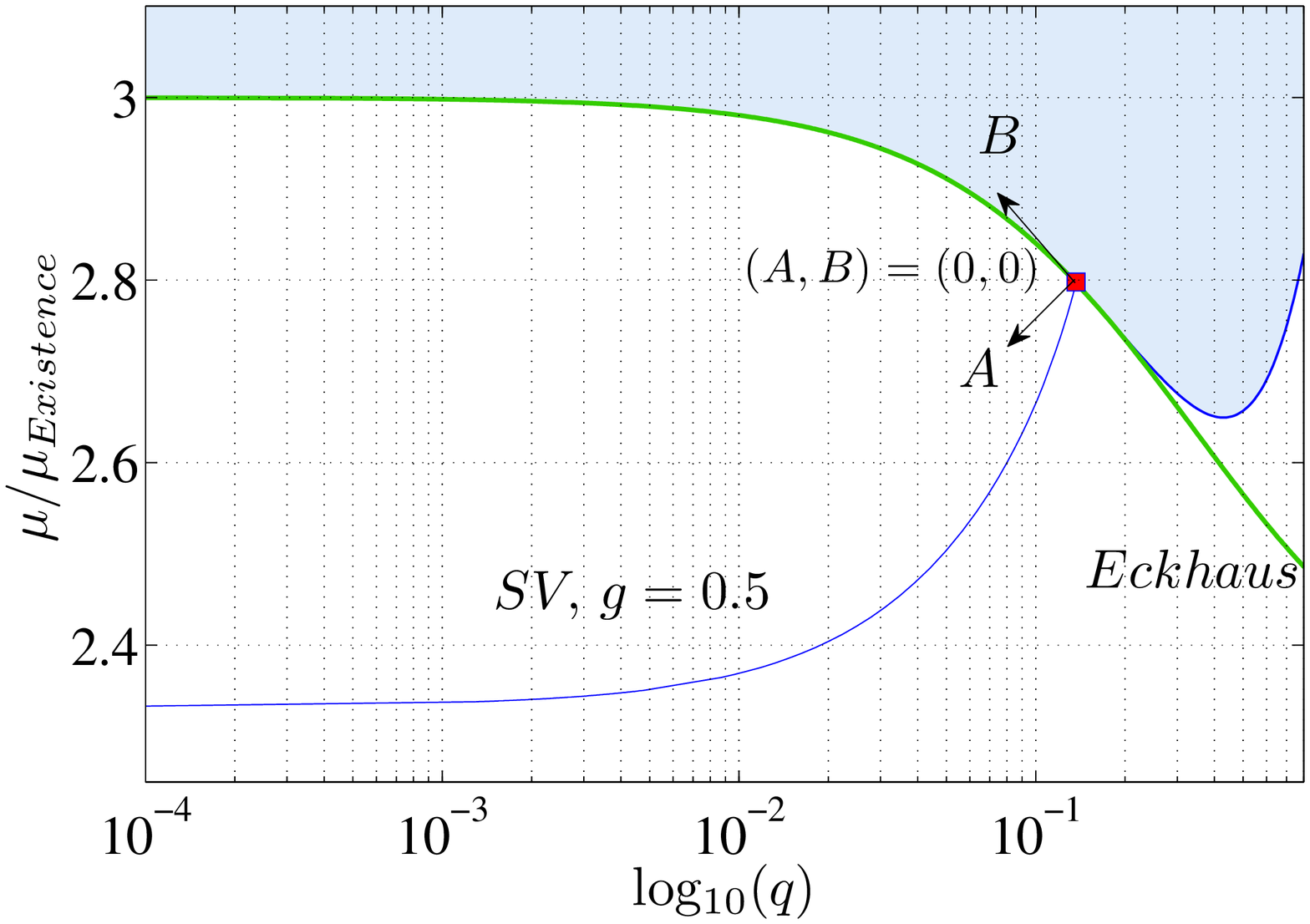}}
\label{fig:7b}
}}
\caption{(Color online) Numerical computation (in model 2) of the SVI boundary as a function of
$\mu/\mu_{Existence}$ and $\log_{10}(q)$. $\mu_{Existence}=(1-(1+q)^{2})^{2}$.
The Eckhaus boundary, $A=0$, is denoted by a green thick curve and the shaded region corresponds to stable stripes. (a) $g=5>g_{critical}$. Here,  $\mu_{SV}/\mu_{Existence}\rightarrow3 $  as $q\rightarrow0$. The SVI precedes the Eckhaus for all $q$ values. (b) $g=0.5<g_{critical}$. Here, $\mu_{SV}/\mu_{Existence}\rightarrow2.3333 $  as $q\rightarrow0$, which in turn becomes $\mu_{SV}\rightarrow9.3333q^{2} $  as $q\rightarrow0$. The point of intersection of the SVI boundary with the Eckhaus boundary is denoted by $(A,B)=(0,0)$. A schematic illustration of the $A$ and $B$ axes at the crossing point is also shown (see figure~\ref{fig:5}).}
\label{fig:7}
 \end{figure}

We then translate the conditions (\ref{Eq:1}) and (\ref{Eq:2}) into bifurcation lines in the $(\mu,q)$ plane.  In order to derive bifurcation lines we have used a  branch-following package, MATCONT~\cite{Dhooge:2004}. An example of the calculation for model 2 is given in figure \ref{fig:6}, which was computed working directly with numerically determined eigenvalues and computing $\frac{\partial }{\partial k^{2}}$ and $\frac{\partial }{\partial l^{2}}$ numerically. In  case I, the SVI boundaries derived using  condition (\ref{Eq:1}) coincide with the numerical computation. However, in case II, condition (\ref{Eq:2}) agrees with the numerical computation only when $A$ is close to zero, as would be expected.  The transition from case I to case II occurs at $(\mu,q)=(0.1367,0.2387)$ for $g=0.5$: for smaller $\mu$, the stability region of stripes is bounded on the right by the Eckhaus instability, while for larger $\mu$, the Eckhaus instability is preempted by the \hbox{SVI}. For $g=5$, the SVI pre-empts the Eckhaus instability for all $\mu$.

A feature of  the SVI is that it  makes stripes unstable only for $q>0$ and, for any given value of $\mu$, $g$ must be large enough for the SVI to preempt the Eckhaus instability. 
The crossing point between the Eckhaus and skew-varicose instabilities for a given value of $q$ occurs for some  $g$, say $g_{Eck}$, which  is given by
\begin{equation*}g_{{Eck}}=\frac{3}{8}\left[\frac{3(1+q)^{2}-1}{(1+q)^{6}}\right].\end{equation*} The condition $A=0$ could be used to express $g_{Eck}$ as a function of $\mu$ if desired. In the limit $\mu\rightarrow0$ and  $q\rightarrow0$,  $g_{Eck}$ goes to $0.75$, which we call $g_{critical}$. The expression for  $g_{Eck}$  and the value $g_{critical}=0.75$ is the same in models 1 and 2 provided we use the relation
$g=g_m/(Pr\,c^{2})$. The Eckhaus instability precedes the SVI for  some range of $\mu$ only if $g<g_{critical}$. We have found that the SVI boundary for $g>g_{critical}$ approaches the origin as
$\mu=12q^{2}$, as does the Eckhaus  curve, whilst  it approaches  as  $\mu=nq^{2}$, with $4<n<12$ when $g<g_{critical}$. A detailed presentation of this asymptotic result will be discussed in section (\ref{sec:asy}). The distinction between   $g>g_{critical}$ and $g<g_{critical}$ is illustrated in figure \ref{fig:6} where we present SVI boundaries for $g=5>g_{critical}$ and  $g=0.5<g_{critical}$.

 Interestingly, for a fixed $q$, $g\rightarrow0$ implies $\mu_{SV} \rightarrow \mu_{Existence}$. Therefore when $g\rightarrow0$,
the SVI boundary coincides with the existence curve of stable stripes for
small $\mu$. 
 
The Eckhaus and SVI boundaries are often very close, so we present our result in an alternative way in figure~\ref{fig:7}.
Figures \ref{fig:7a} and  \ref{fig:7b} present the SVI
boundaries for $g=5$ and $g=0.5$ in the ($\mu/\mu_{Existence}$ , $\log_{10}(q)$) plane. 
This is a better way of illustrating the  regions of stable stripes (shaded regions) and the behavior of  the Eckhaus and SVI boundaries as $q\rightarrow0$.
Moreover, it  shows how the coordinates axes $A$
and $B$ from (14) can be defined near the SVI--Eckhaus crossing point. 
 
\subsection{Asymptotic analysis of the SVI boundary}
\label{sec:asy}
In this section we focus on the asymptotic behavior of the SVI boundary in model 2, in the two cases discussed above, $g>g_{critical}$ and $g<g_{critical}$. Let us  consider the case when $g>g_{critical}$, where the SVI boundary  pre-empts the Eckhaus instability boundary. We use the condition $B^{2}-4AC=0$, which can be written as 
\begin{equation}F_{1}g^{2}+F_{2}g+F_{3}=0,\label{eq:before}\end{equation}
where the $F_{i}$'s are functions of $\mu$ and $q$. In the limit of  very small $\mu$ and $q$, 
\begin{subequations}
\begin{gather}
F_{1}\approx\frac{1}{9}({16384}q^{6}-{8192}\mu q^{4}+{1024}\mu^{2}q^{2}-256\mu^{3}q+16\mu^{4}),\\
F_{2}\approx\frac{1}{3}({-24576}q^{5}+{8192}\mu q^{3}-{512}\mu^{2}q-64\mu^{3}),\\
F_{3}\approx({9216}q^{4}-1536\mu q^{2}+{64}\mu^{2}).
\end{gather}
\label{eq:f}
\end{subequations} We note at this point that equation (\ref{eq:before}) is valid for model 1 (with the same values of $F_1$, $F_2$ and $F_3$) when we identify $g$ with $g_{m}/{Pr\,c^2}$.
For smallish $g\geq g_{critical}$, as $\mu\rightarrow0$, $\mu\sim q^{2}$  as shown in figure \ref{fig:7a} and hence functions in equation (\ref{eq:f}) 
can be taken as $F_{1}\sim q^{6}$, $F_{2}\sim q^{5}$ and $F_{3}\sim q^{4}$. Therefore the equation (\ref{eq:before}) can be approximated as  $F_{3}=0$, in which case $\mu=12q^{2}$. This can be improved by including $F_{2}$ and hence $F_{2}g+F_{3}=0$ in which case we obtain
 \begin{equation}\mu_{SV}=12q^{2}\left(\frac{\displaystyle9-8qg}{\displaystyle9-24qg}\right),\,(g>g_{critical}\,,q\ll1/g)\label{small}\end{equation}
and hence $\mu\rightarrow12q^{2}$ when $q\rightarrow0$. Therefore when $g>g_{critical}$, the SVI boundary for small~$\mu$ has the same curvature as the Eckhaus boundary. Note that figure~\ref{fig:7a} shows how the SVI boundary for $g=5$ coincides with the Eckhaus boundary as $q\rightarrow0$. 
\begin{figure}[t]
\begin{center}
\scalebox{0.37}{\includegraphics*[viewport=0.3in 0.4in 9.49in 6.45in]{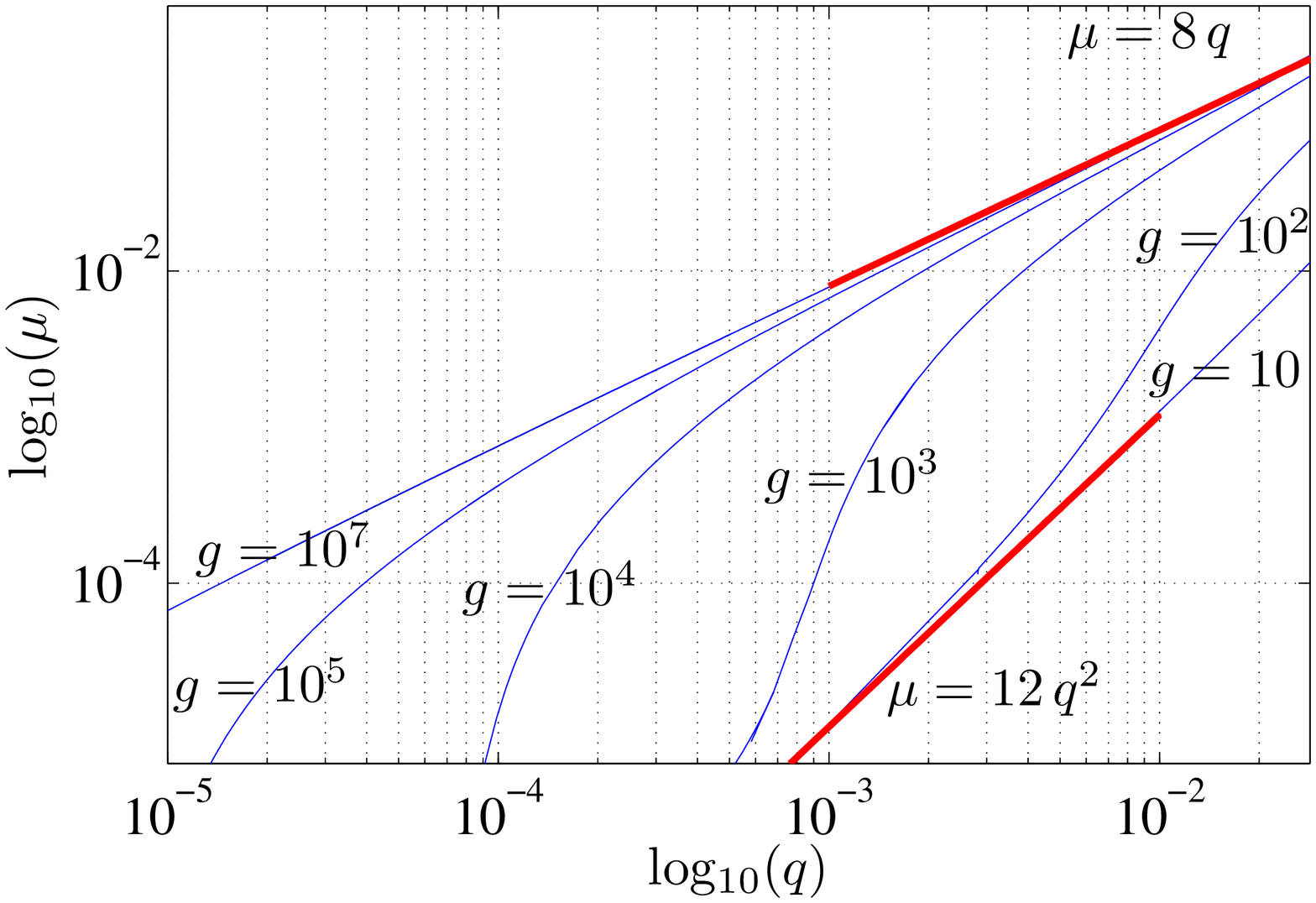}}
\end{center}

\caption{(Color online) Numerical computation of the SVI boundary on a logarithmic scale for $g=10^{i}$  for $i=1,2,3,4,5$ and $7$. The transition from $\mu\sim q$ to $\mu\sim q^2$ occurs when $g\propto 1/q$. In the limit of small $q$, the SVI curve is
tangent to $\mu=12q^{2}$,  whereas in the limit of large $g$, the SVI curve goes as
$\mu\sim 8q $ as $\mu$ increases. Both asymptotes, $\mu\sim 8q$ and $\mu\sim 12q^2$ are denoted by
red thick lines.}
 \label{fig:8}
\end{figure}
The limit of very large $g$ is also of interest, as it corresponds to stress-free boundary conditions (see section (\ref{osv}) below). In this limit, we expect the SVI boundary (as shown in figures~\ref{fig:14} and~\ref{fig:16} below) should have $\mu\sim q$ \cite{Bernoff:1994}. From equation (\ref{eq:f}) with $\mu\sim q$, we have $F_{1}\sim q^{4}$, $F_{2}\sim q^{3}$ and $F_{3}\sim q^{2}$ and $F_{1}g^{2}+F_{2}g+F_{3}\sim q^{2}\left(q^{2}g^{2}+qg+1\right)$. If $qg\ll1$, we set $F_{3}=0$ and recover $\mu=12q^{2}$. If $qg\gg1$, we set $F_{1}=0$ and obtain $\mu\sim8q$. Again, this can be improved by setting $F_{1}g+F_{2}=0$, and we obtain
 \begin{equation} \mu_{SV}=8q\left(\frac{2qg+3}{2qg-3}\right),\,(g>g_{critical}\,,q\gg1/g).\label{large} \end{equation}
Finally, we note that the transition between $\mu\sim12q^{2}$ and $\mu\sim8q$ will occur when $qg$ is of order unity, so $q_{transition}\sim\frac{1}{g}$. These three regimes are illustrated in figure \ref{fig:8}.

In the skew-varicose mechanism for $g>g_{critical}$, the maximum of $Det$ is attained for perturbations of a mode, say $(k_{max},l_{max})$, as $k_{max}\rightarrow0$ and $l_{max}\rightarrow0$. As shown by the equation (\ref{Eq:1}), $\frac{{k_{max}}^{2}}{{l_{max}}^{2}}=\sqrt{\frac{C}{A}}$ and so $\frac{{k_{max}}}{{l_{max}}}=\mathcal{O}(1)$.    

Secondly, let us consider the case $g<g_{critical}$, for which the SVI boundary lies in the Eckhaus band for small $\mu$. We do not have an exact criterion for the SVI in this case, though equation (16) is an approximate criterion. However, for small $\mu$, we expect $\mu_{SV}$ to depend approximately linearly on $g$. We know from equation (19) that for $g=g_{critical}$,  $\mu_{SV}\rightarrow 12q^{2}$  as  $q\rightarrow0$ and from conditions in equation (16) in the limit $g\rightarrow0$, we have $\mu_{SV}\rightarrow 4q^{2}$  as  $q\rightarrow0$. A linear interpolation between these yields 
\begin{equation}\mu_{SV}\approx4\left(\frac{2}{g_{critical}}g+1\right)q^{2},\,(g<g_{critical}\,,q\ll1).\label{eq:cur1}\end{equation} This relation is shown numerically using MATCONT~\cite{Dhooge:2004} to be a very good approximation; $\mu_{SV}/\mu_{Existence} $ behaves linearly with $g$ with a gradient $8/3$. In addition, the numerical simulation illustrated in figure~\ref{fig:7b}, shows  at $g=0.06$, $\mu/\mu_{Existence}\rightarrow2.6$ as $q\rightarrow0$ which agrees
well with equation (\ref{eq:cur1}). All these explicit results are for model 2. When $g>g_{critical}$, the expressions for conditions given in (\ref{Eq:1}) are the same for model 1 with the relation $g=g_m/Pr\,c^{2}$. Therefore, equations (\ref{small}) and (\ref{large}) are the same for model 1.
When $g<g_{critical}$, for model 1, the condition for the SVI given by  equation (\ref{Eq:2}) with $g=g_m/Pr\,c^{2}$ is not the same as in model 2. However, for fixed $Pr$ and $c\neq0$, we find that $\mu_{SV}/\mu_{Existence} $ behaves linearly with $g_{m}$  and hence equation (\ref{eq:cur1}) holds for model 1. The case where $c=0$ is considered in the next section.

\section{LONG-WAVELENGTH INSTABILITIES: Skew-varicose and Oscillatory skew-varicose with stress-free boundary conditions. }
\label{osv}
We now consider the case that models convection with stress-free boundary condition. In model 1, $c$ is the parameter that accounts for the boundary conditions at the top and bottom,  and 
stress-free boundary conditions correspond to $c=0$. Using the relation  $g=g_m/Pr\,c^{2}$ to connect the two models, taking the limit $g\rightarrow\infty$ corresponds to stress-free boundary condition in model 2. 

When $c=0$, for any coupling constant $g_{m}$, the SVI boundary always pre-empts the Eckhaus boundary in the $(\mu,q)$ plane. To show this we can focus on model 1, for which the SVI condition is  $F_{1}{g_m}^{2}+F_{2}g_mPr\,c^{2}+F_{3}(Pr\,c^{2})^{2}=0$  with $F_{1}$, $F_{2}$ and $F_{3}$  given in equation (\ref{eq:f}); this confirms the relation $g=g_m/Pr\,c^{2}$. Hence $c=0$ implies $F_{1}=0$, and this gives the criterion for instability  as $\mu=8q$ for any $g_{m}$, as above. A remarkable property of the SVI in model 1 with stress-free boundary conditions is that the stability boundary is independent of $Pr$ and $g_{m}$.  

Another instability  of interest in the stress-free case is the oscillatory skew-varicose (OSV) instability, which consists of a long-wavelength transverse oscillations of  the stripes that propagate along their axis \cite{Bolton:1985}. For  OSV modes, the associated eigenvalues are complex. This instability does not occur in model 2  because $Tr(J_{2})<0$ for small $k$ and $l$ and so complex eigenvalues are not possible. In model $1$, the OSV instability does not appear for $c^2=2$.  With $c=0$, the regions in the $(k,l)$ plane with  positive growth rate and nonzero frequency emerge from  $(k,l)=(0,0)$ in a way that resembles the contours for the SVI, shown in figure \ref{fig:3}. 
\begin{figure}[t]
\scalebox{0.37}{\includegraphics*[viewport=0.5in 0.7in 9.6in 6.41in]{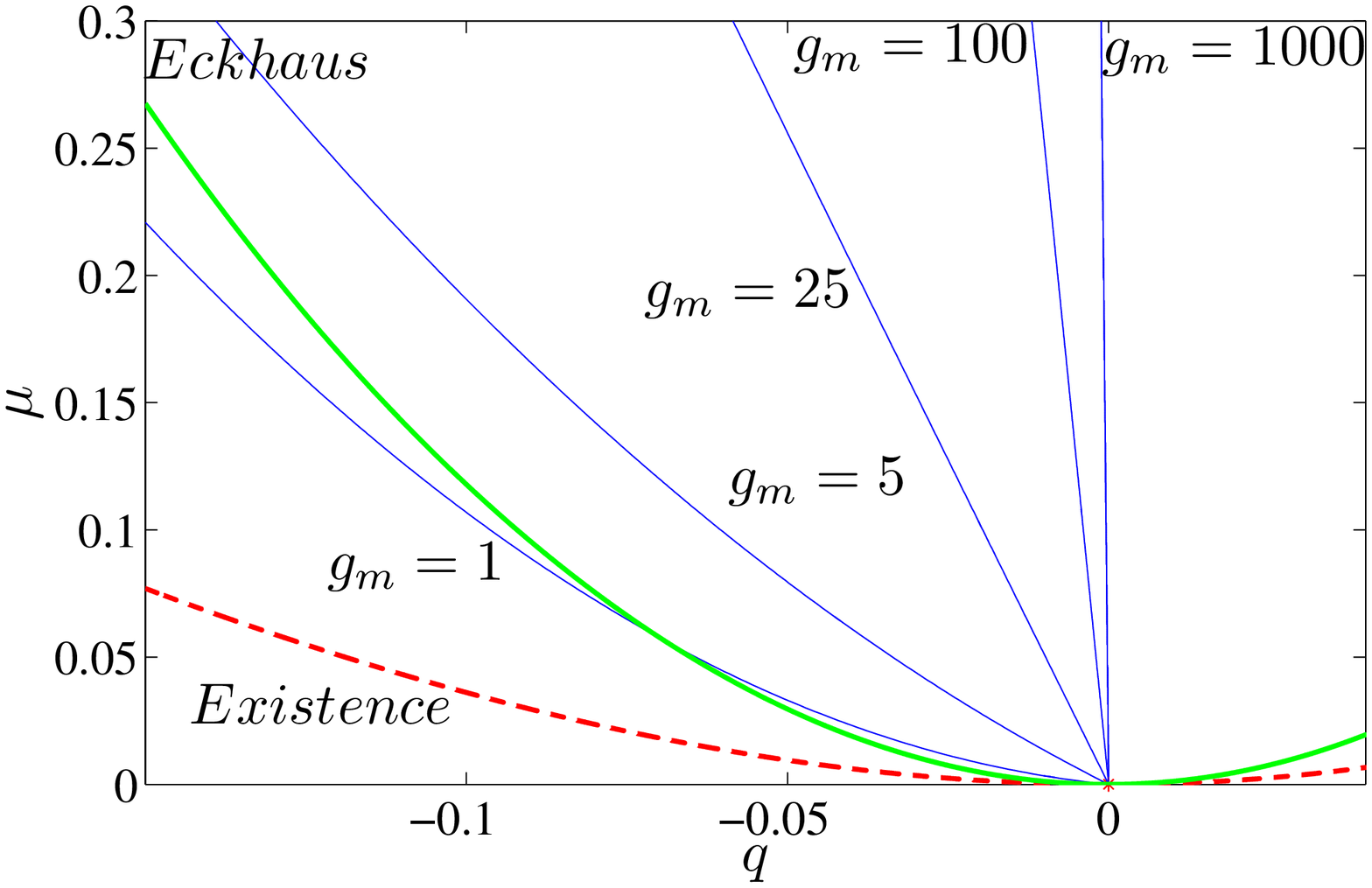}}

\caption{(Color online) The location of the OSV instability boundary for model 1 for $c=0$ (stress-free boundary conditions), $Pr=1$ and  $g_{m}=1000,100,25,5$ and $1$. Stripes are OSV unstable to the left of the OSV boundary. For small $\mu$, the boundary is asymptotic to $\mu=\left(\frac{-3+\sqrt{5}}{3}\right)qg_{m}$.  The Eckhaus boundary is denoted in green thick and the existence curve is in dashed red.}
 \label{fig:osv}
\end{figure}
 
We derive the asymptotic behavior of the OSV instability boundary. At the point of instability, the eigenvalues are purely imaginary; this gives two conditions, $C-AB=0$ and $B>0$, where the characteristic equation of $J_{1}$ is $\lambda^{3}+A\lambda^{2}+B\lambda+C=0$. The coefficients  $A$, $B$ and $C$ are functions of $g_{m}$, $q$, $\mu$, $k$ and $l$. We have set $Pr=1$ for illustrating this calculation. For small $k$ and $l$, the condition $C-AB=0$  gives, 
\begin{equation}G_{1}\left({\frac{g_{m}q}{\mu}}\right)^{2}+G_{2}{\left({\frac{g_{m}q}{\mu}}\right)}+G_{3}=0,\label{eq:osvbefore}\end{equation}
after maximizing $C-AB$ over $(k,l)$. The $G_{i}$'s are functions of $\mu$ and $q$. In the limit of  very small $\mu$ and $q$, we find,  
\begin{subequations}
\begin{gather}
G_{1}\approx\frac{1}{9}\left({1024}-{8192} \left(\frac{q^{2}}{\mu}\right)+{16384}{\left(\frac{q^{2}}{\mu}\right)}^{2}\right)\;,\\
G_{2}\approx{512}-\frac{26624}{3} \left(\frac{q^{2}}{\mu}\right)+{46421}{\left(\frac{q^{2}}{\mu}\right)}^{2}-{76459}{\left(\frac{q^{2}}{\mu}\right)}^{3}\;,\\
G_{3}\approx{256}-{8192} \left(\frac{q^{2}}{\mu}\right)+{90112}{\left(\frac{q^{2}}{\mu}\right)}^{2}-{393216}{\left(\frac{q^{2}}{\mu}\right)}^{3}\nonumber\\+{589824}{\left(\frac{q^{2}}{\mu}\right)}^{3},
\end{gather}
\label{eq:g}
\end{subequations} 
\vspace{-0.1in}
\begin{flushleft}
where we have dropped terms that can be shown to be smaller than those retained. 
\end{flushleft}
\vspace{-0.05in}
We note at this point that equations (\ref{eq:osvbefore}) and (\ref{eq:g})  with $(q^{2}/\mu)\ll1$ give $\frac{1024}{9}\left({\frac{g_{m}q}{\mu}}\right)^{2}+512{\left({\frac{g_{m}q}{\mu}}\right)}+256=0$, which implies  $\mu=\left(\frac{-3+\sqrt{5}}{3}\right)qg_{m}$.
 Hence, the OSV instability  boundary  has a linear relationship between $q$ and $\mu$ for small $\mu$, and  it bounds the region of stable wavenumbers for negative $q$. Stripes are stable on the right of the boundary. In this asymptotic limit, the point of maximum growth rate in the $(k,l)$ plane can be found at the point of maximum of $C-AB$; this point satisfies $l/k=\sqrt5$. Figure \ref{fig:osv} shows the behavior of the OSV instability boundary in model 1 with $c=0$, for different values of the coupling constant $g_{m}$.

When $c$ is increased from zero, the OSV instability turns into the so-called oscillatory instability \cite{Busse:1972}. This has the nature of an oscillatory cross-roll instability, setting in with non zero $k$ and $l$. The boundary of this oscillatory instability  emerges from  $\beta=0$, the existence curve. However, this instability is prominent only for  $Pr\,c^2\ll1$; for higher values of  $Pr\,c^2$, the instability moves to larger negative $q$. 
\begin{figure}[t]
\begin{center}
\scalebox{0.37}{\includegraphics*[viewport=0.5in 0.8in 9.6in 6.3in]{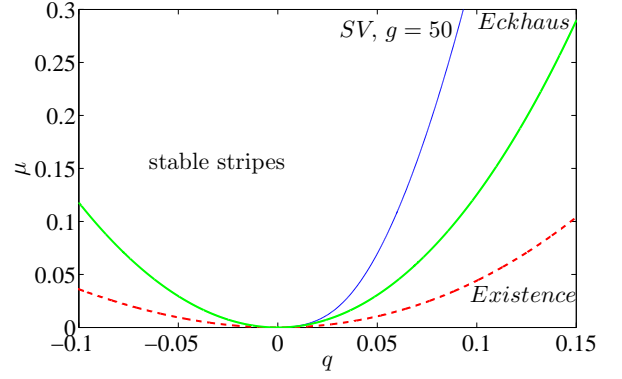}}
\end{center}

\caption{(Color online) Stability diagram in the neighborhood of $\mu=0$ for model 1 with $c^{2}=2$ (no-slip boundary conditions), $Pr=1$, $g_{m}=50$ and $\gamma=2.5$. Stable stripes are in the region indicated, bounded by Eckhaus instability (green thick curve) from below and by the SVI boundary (blue thin curve) from above.  }
 \label{fig:sv50}
\end{figure}
\begin{figure}[h]
\vspace{-0.1in}
\scalebox{0.36}{\includegraphics*[viewport=0.5in 0.7in 9.6in 6.6in]{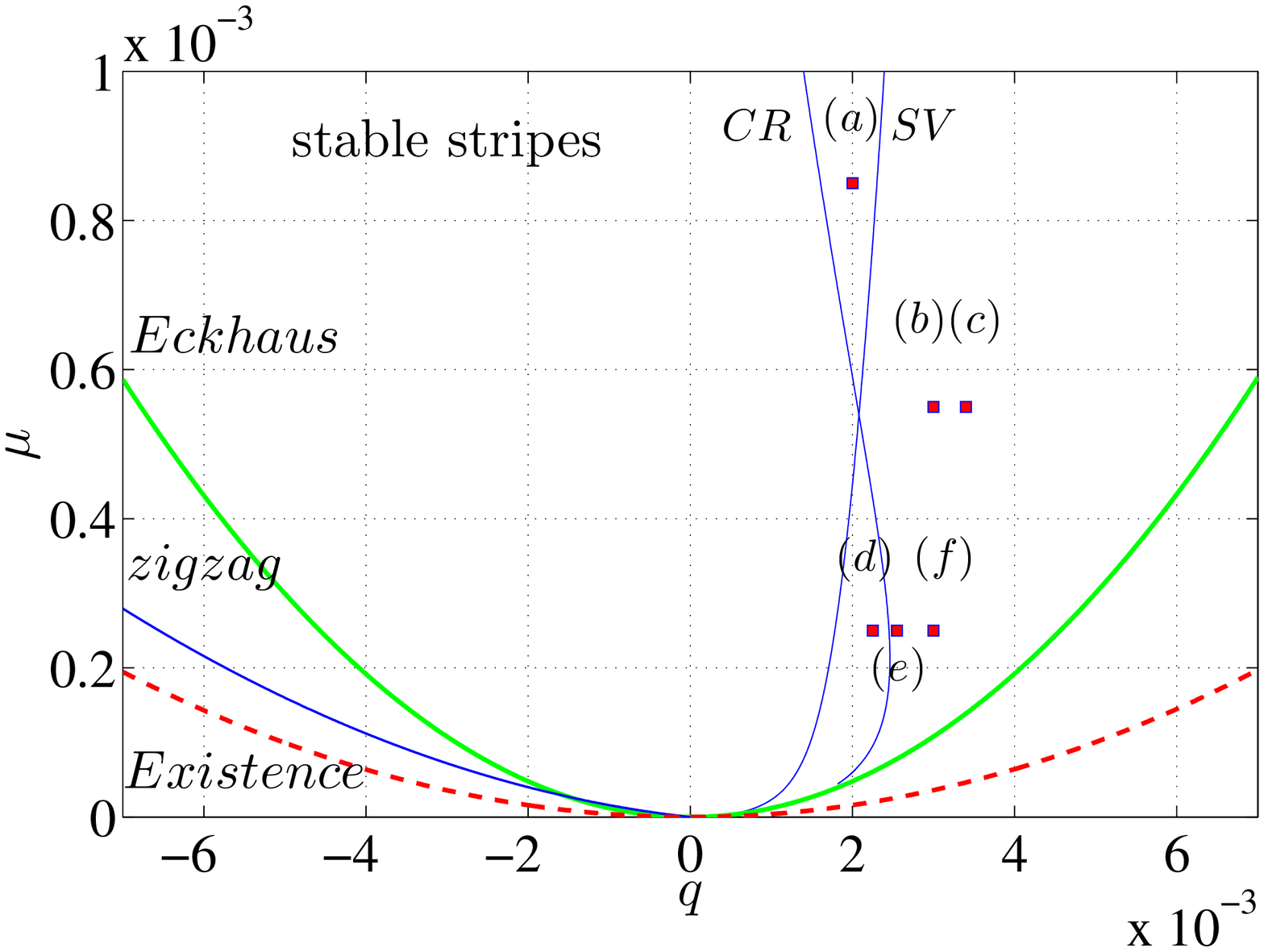}}

\caption{(Color online) Stability diagram in the neighborhood of $\mu=0$ for model 1 with $c^{2}=2$ (no-slip boundary conditions), $Pr=1$, $g_{m}=1000$ and $\gamma=2.5$. Stable stripes are in the region indicated, bounded by the SV and CR instabilities from above, and by the zigzag and  Eckhaus  instabilities from below. Stable stripes exist at $q=0$ for range of $\mu$ close to zero. Growth rates as a function of $k$ and $l$ at the points indicated by (a)--(f) are given in figure \ref{fig:13}. }
 \label{fig:zoom}
\end{figure}

\begin{figure*}[t]
\vspace{-0.15in}
\mbox{
\vspace{-0.4in}
\subfigure[  ]{
\scalebox{0.3}{\includegraphics*[viewport=0.4in 0.7in 7.92in 6.6in]{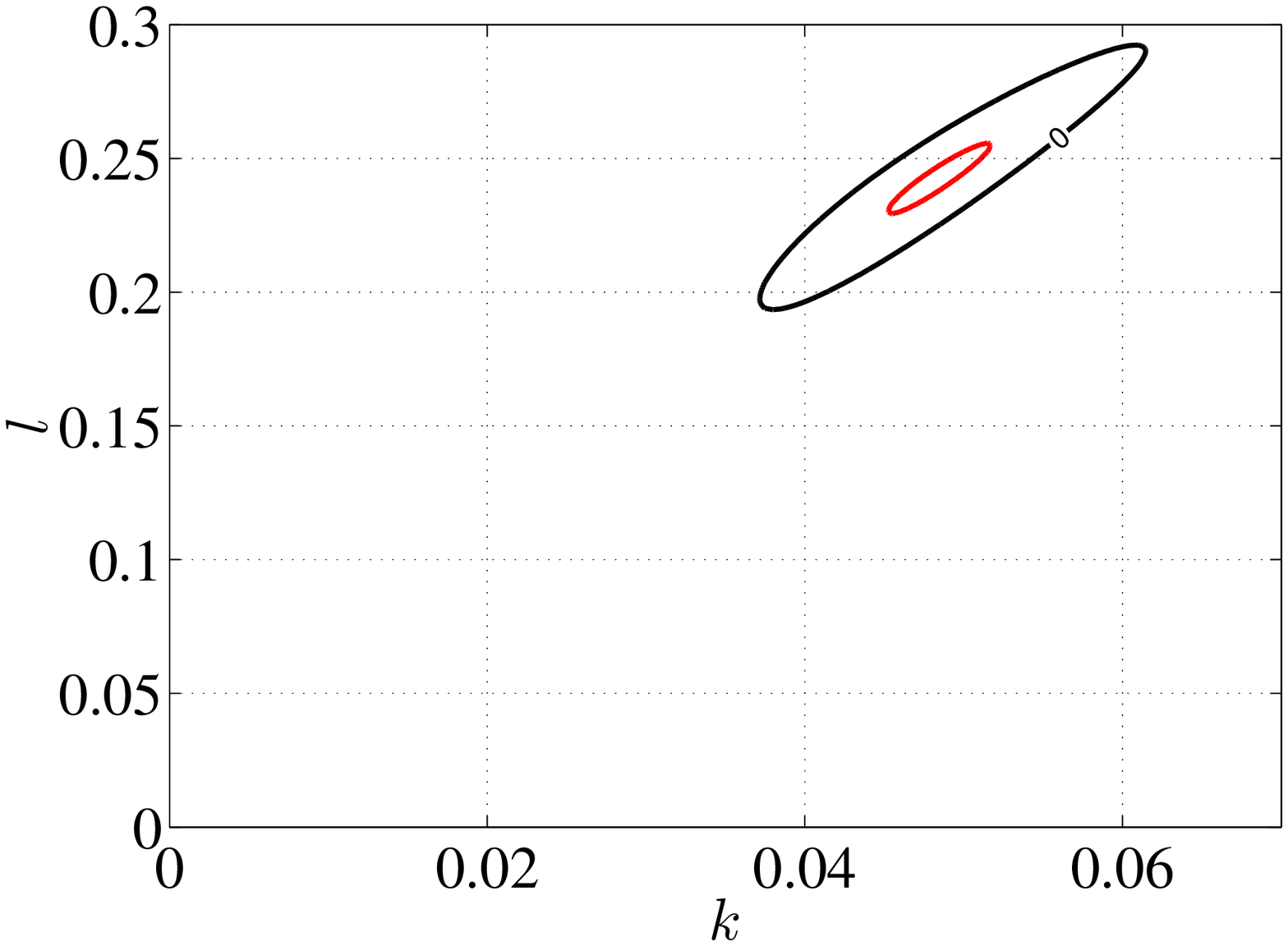}}
\label{fig:gr1}}
\subfigure[  ]{
\scalebox{0.3}{\includegraphics*[viewport=0.4in 0.7in 7.92in 6.6in]{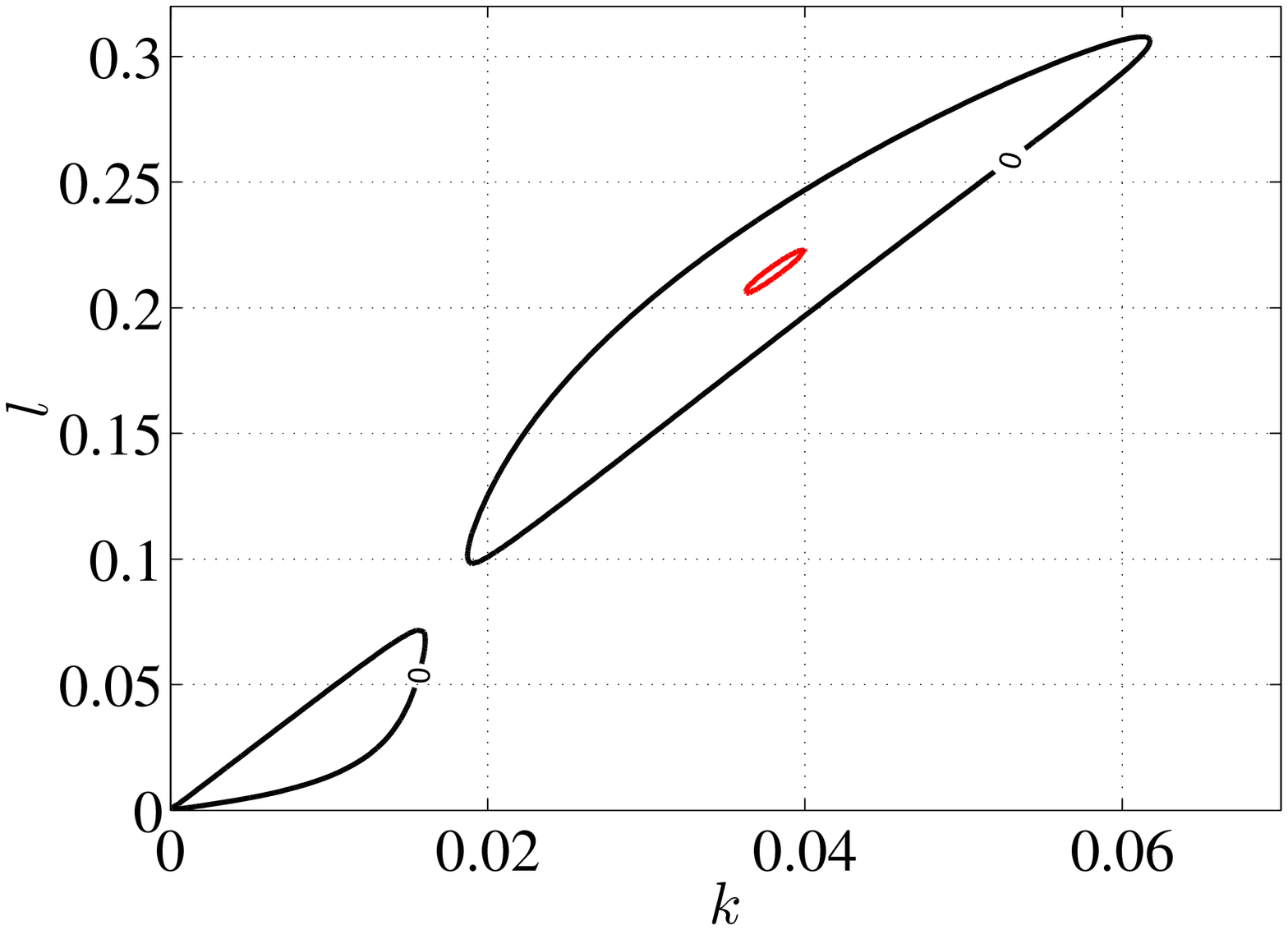}}
\label{fig:gr2}}
\subfigure[  ]{
\scalebox{0.3}{\includegraphics*[viewport=0.35in 0.7in 7.9in 6.6in]{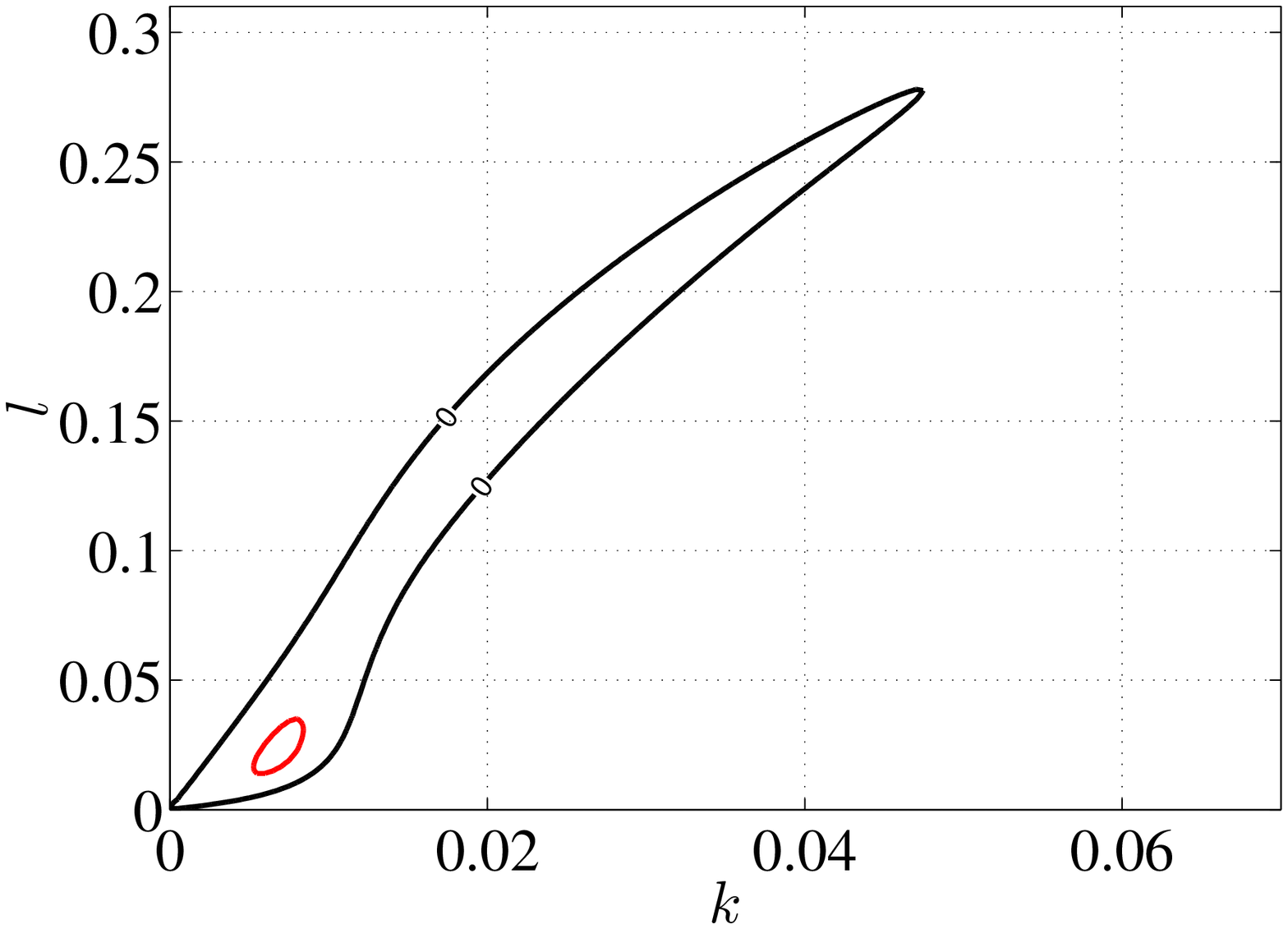}}
\label{fig:gr3}}
}
\mbox{
\vspace{-0.7in}
\subfigure[  ]{
\scalebox{0.3}{\includegraphics*[viewport=0.4in 0.7in 7.9in 6.6in]{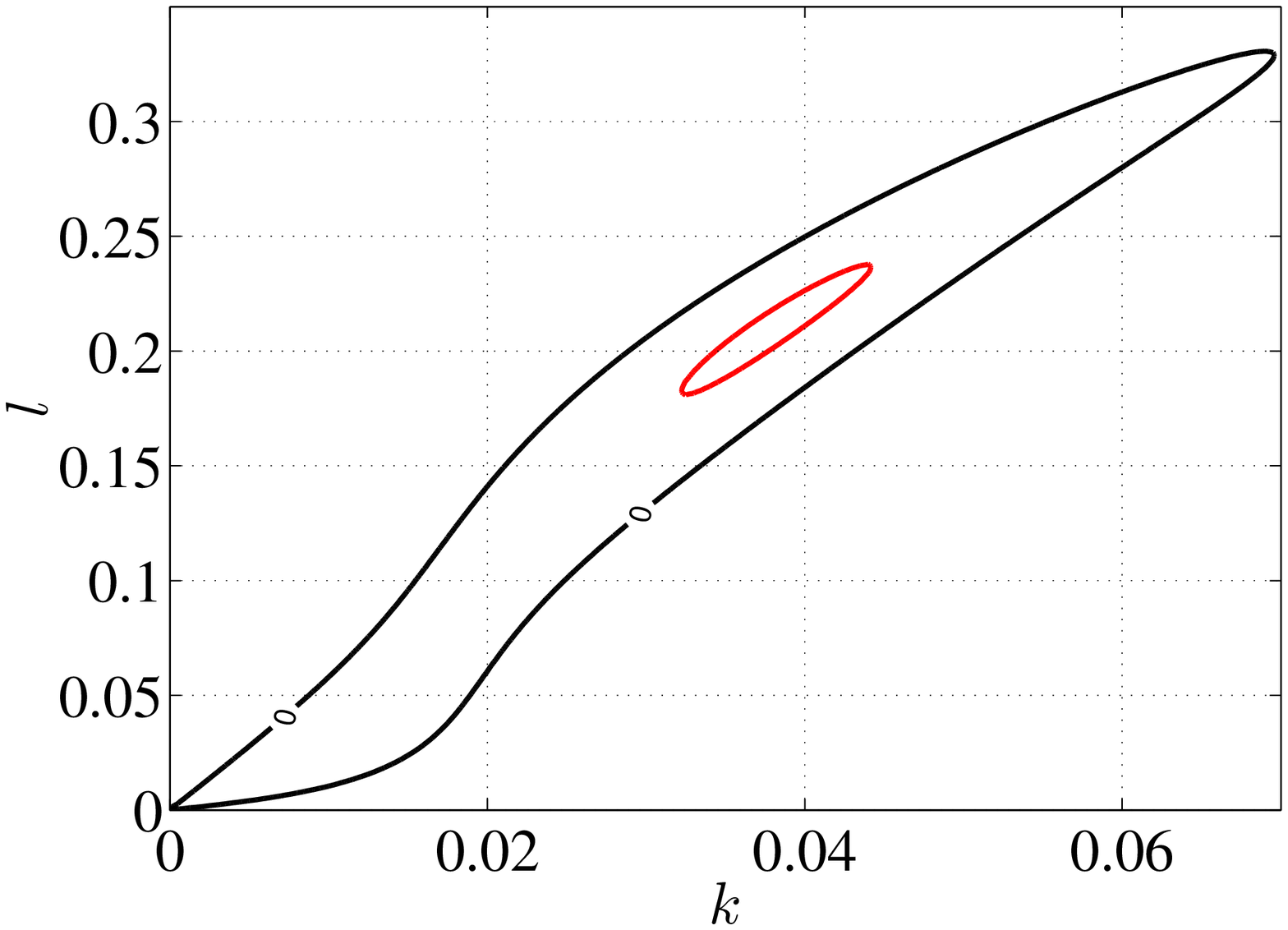}}
\label{fig:gr1}}
\subfigure[  ]{
\scalebox{0.3}{\includegraphics*[viewport=0.4in 0.7in 7.92in 6.6in]{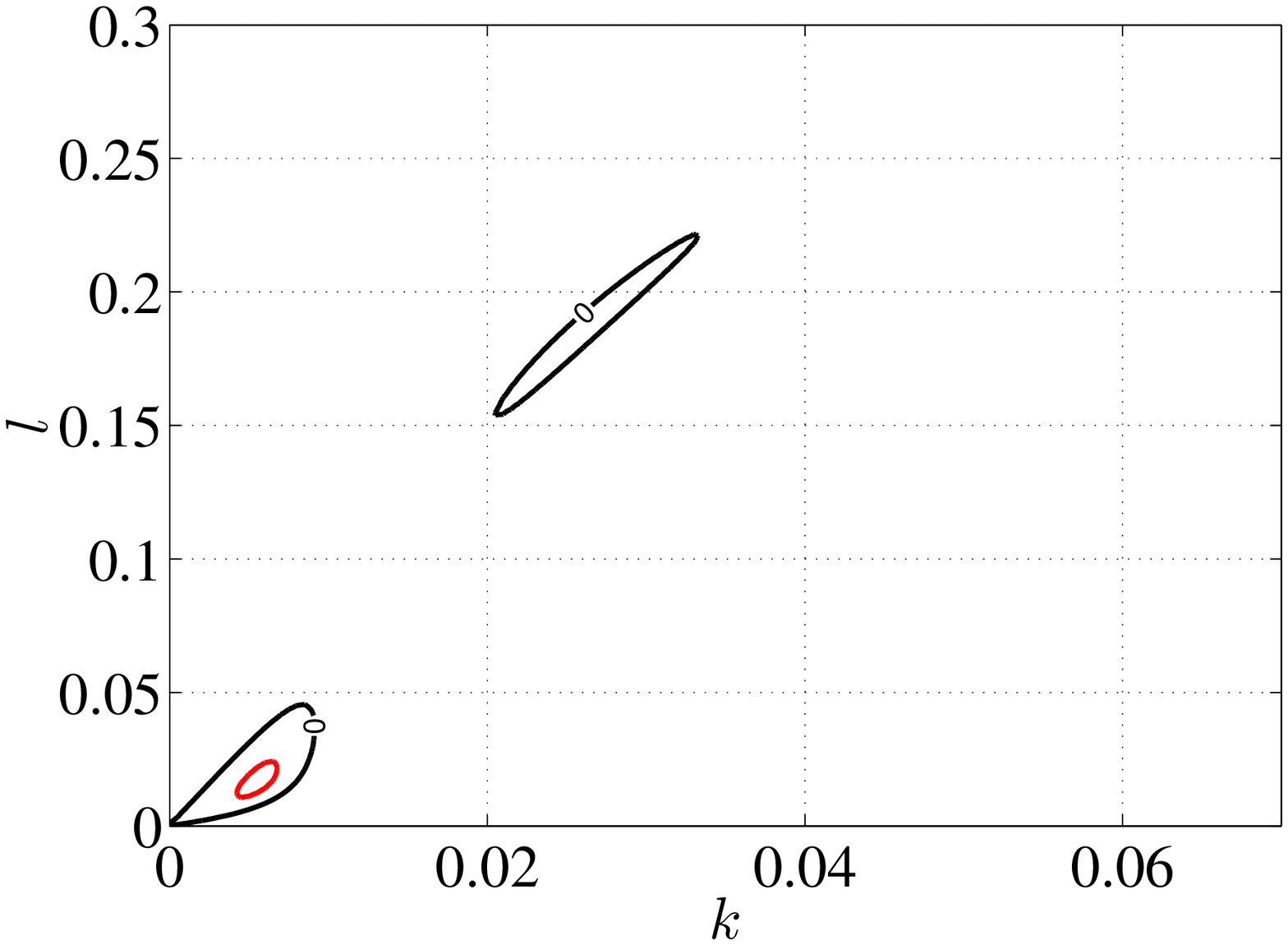}}
\label{fig:gr2}}
\subfigure[  ]{
\scalebox{0.3}{\includegraphics*[viewport=0.35in 0.7in 7.92in 6.6in]{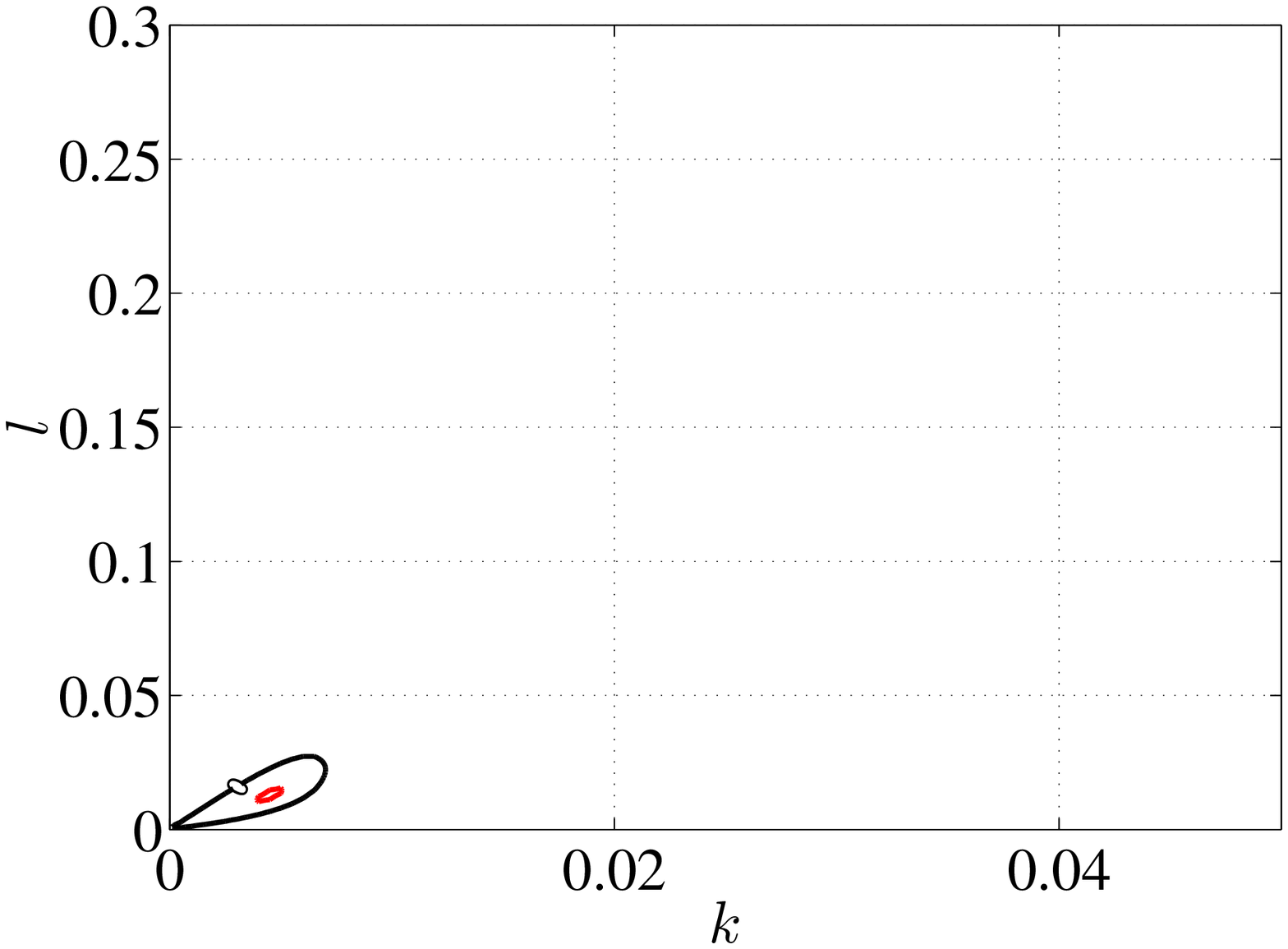}}
\label{fig:gr3}}}

\caption{(Color online) Growth rates of perturbations at selected $(\mu,q)$ indicated by (a)--(f) in figure~\ref{fig:zoom}. (a) $\mu=8.5\times10^{-4}$, $q=0.002$; stripes are CR unstable, but growth rates for $k$ and $l$ close enough to zero are negative. (b) $\mu=5.5\times10^{-4} $, $q=0.003$; stripes are CR and SV unstable. Growth rates for $k$ and $l$ close enough to zero become positive due to the SVI and there are two distinct zero contours. (c) $\mu=5.5\times10^{-4} $, $q= 0.0034 $; stripes are CR and skew-varicose unstable. One large contour encloses both the SV and CR instabilities.  (d) $\mu=2.5\times10^{-4}$, $q=0.00225$; stripes are skew-varicose unstable, and growth rates for $k$ and $l$ close enough to zero are positive for a range of polar angles. (e) $\mu= 2.5\times10^{-4}$, $q=0.0025 $; same as case (b), but the  maximum occurs in the SV region. (f) $\mu= 2.5\times10^{-4}$, $q=0.003 $; same as case (c), but again the maximum occurs in the SV region. The zero contour  is denoted in black (outer) and contours of positive growth rate are denoted in red (inner in gray).}
 \label{fig:13}
 \end{figure*}

\section{SHORT-WAVELENGTH instabilities: The Cross-roll instability.}
\label{cr}
The cross-roll (CR) instability is  so-called because the fastest growing disturbances appear to take the form of stripes perpendicular to the basic steady stripe pattern: these disturbances have  non-zero $k$ and $l$ in the limit $\mu\rightarrow0$~\cite{Busse:1971}. In contrast to the oscillatory cross-roll instability discussed above, the most unstable eigenvalue at the CR instability is zero.

We will show numerically in section (\ref{sec:3}) below that the CR instability only forms a boundary of the region of stability of stripes if $g$ is large. The filtering $\mathcal{F}_{\gamma}$, discussed in section (\ref{sec:2}), was introduced to suppress the CR instability in favour of the SVI, so the location of the CR instability boundary depends on $\gamma$, while the SVI is not influenced by $\gamma$ for small enough $\mu$ and $q$. Since the CR instability sets in with non-zero $k$ and $l$ even with small $\mu$ and $q$, asymptotic analysis of the type carried out above cannot be done. Numerical examples are shown in the next section.
\section{Numerical illustration of combined stability boundary.}
\label{sec:3}
We now present numerical results obtained with MATCONT~\cite{Dhooge:2004}. Details of the specification of the conditions on the eigenvalues for each of the six types of instability are given in the Appendix. We choose two illustrative parameter values: $g_{m}=50$ and  $g_{m}=1000$, and we consider no-slip $(c^{2}=2)$ and stress-free $(c=0)$ boundary conditions.

Figure~\ref{fig:sv50} shows the stability diagram for model 1 with  parameter value $c^2=2$ (corresponding to no-slip boundary conditions), $Pr=1$, $g_m=50$ and $\gamma=2.5$. The region of stable stripes is bounded by the SVI from above and the Eckhaus instability from below. The zigzag instability boundary for this value of $g_{m}$ lies below the left Eckhaus boundary except for very small $\mu$ and $q$ (as in figure 2). The CR instability does not occur in this range of parameters.

\begin{figure}[t]
\vspace{-0.13in}
\scalebox{0.36}{\includegraphics*[viewport=0.5in 0.8in 9.6in 6.3in]{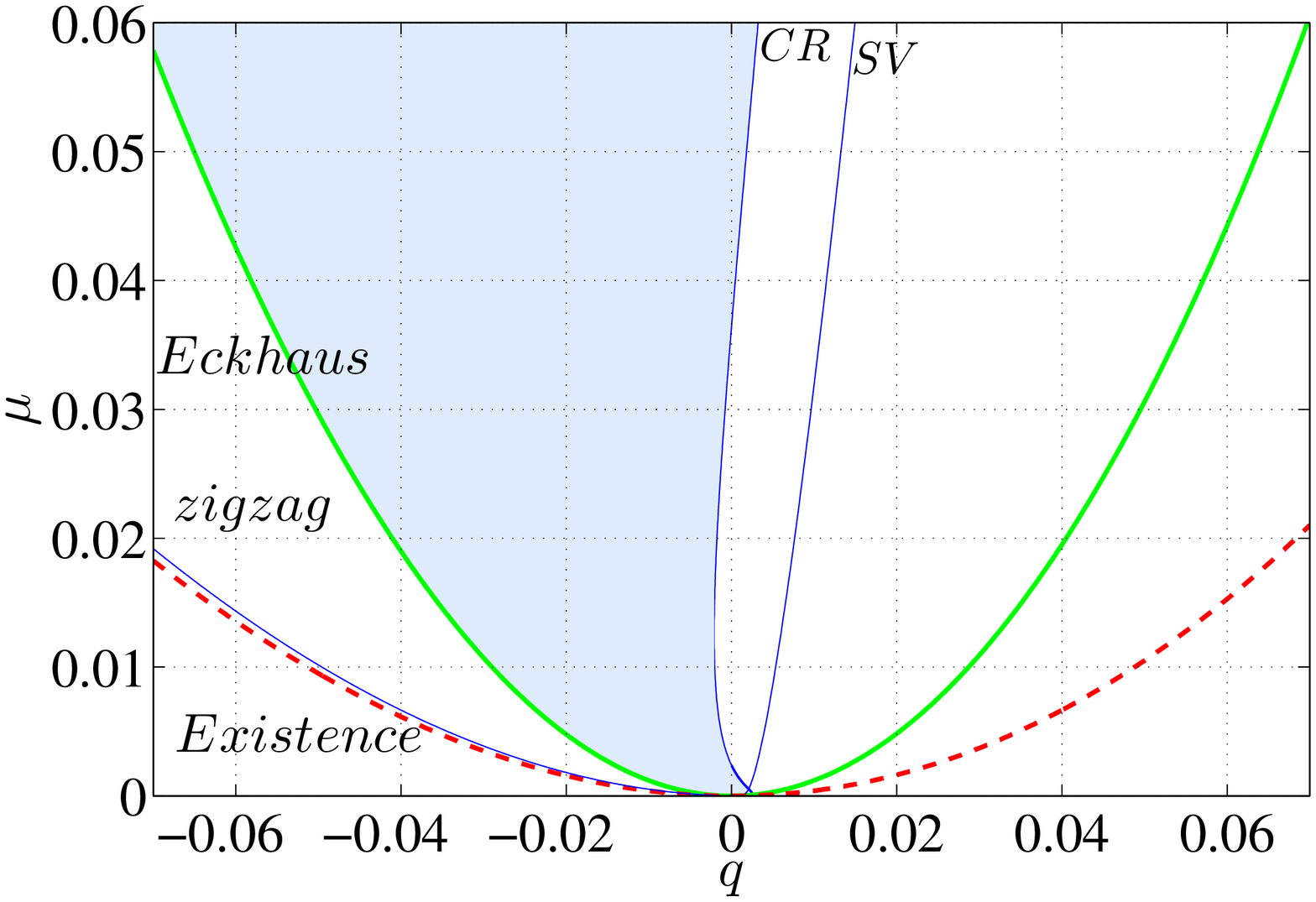}}

\caption{(Color online) Stability diagram  with parameters as in figure \ref{fig:zoom} covering a larger range of $q$ and $\mu$. Stripes are stable in the shaded region, bounded by the Eckhaus and CR instabilities.  }
 \label{fig:rwc}
\end{figure}

Figure~\ref{fig:zoom} and~\ref{fig:rwc} show the stability diagram for  model 1  with parameters $c^2=2$ (no-slip), $Pr=1$, $g_m=1000$ and $\gamma=2.5$.  Close to  $\mu=0$~(figure~\ref{fig:zoom}), a neighborhood of $q=0$ is in the stable regime and the zigzag and the SV instabilities bound the region of stable stripes. For $\mu>5\times10^{-4}$, the region of stable  stripes is bounded  by the CR instability from above and by the Eckhaus instability from below. For this value of  $g_m$, the CR instability boundary crosses the SVI boundary, and the zigzag instability boundary has a linear behavior for small $\mu$.

\begin{figure}[!htbp]
\vspace{-0.13in}
\scalebox{0.36}{\includegraphics*[viewport=0.5in 0.8in 9.6in 6.3in]{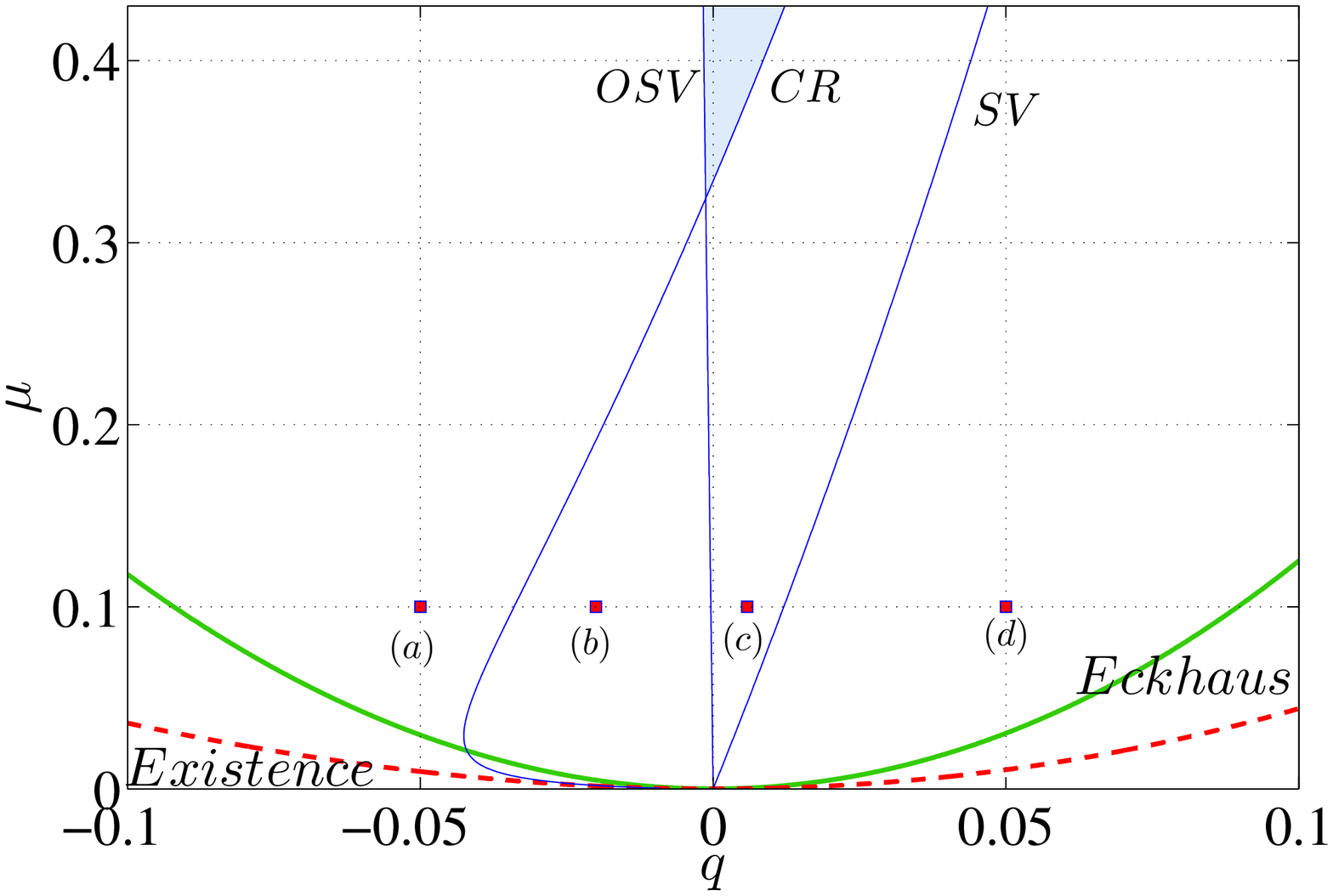}}
\caption{(Color online) Stability diagram for model 1 with $c=0$ (stress-free boundary conditions), $Pr=1$, $g_{m}=1000$ and $\gamma=2.5$. Stripes are skew-varicose unstable to the right of the SVI boundary, OSV unstable to the left of the OSV boundary and CR unstable to the right of the CR boundary.
 Stable stripes exist for  $\mu>0.32$; the region of stable stripes  is shaded, and is bounded by the OSV instability boundary on the  left and the CR boundary on the  right. Growth rates as a function of $k$ and $l$ at the points indicated by (a)--(d) are given in figure 15.}
 \label{fig:14}
\end{figure}

\begin{figure}[t]
\mbox{
\vspace{-0.8in}
\hspace{-0.2in}
\subfigure[  ]{
\scalebox{0.33}{\includegraphics*[viewport=1in 0in 5in 4.1in, width=5in, height=4.5in]{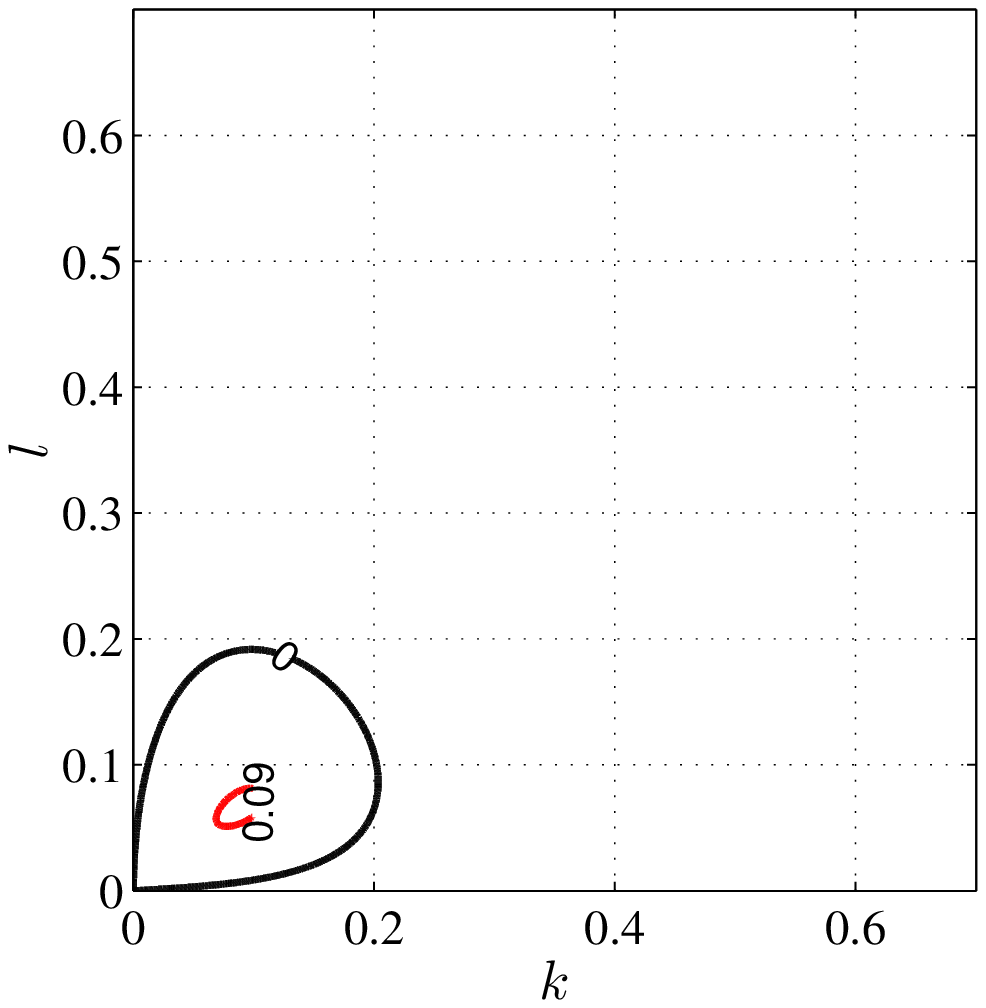}}
\label{fig:st1}}
\subfigure[  ]{
\hspace{-0.2in}
\scalebox{0.33}{\includegraphics*[viewport=1in 0in 5in 4.1in, width=5in, height=4.5in]{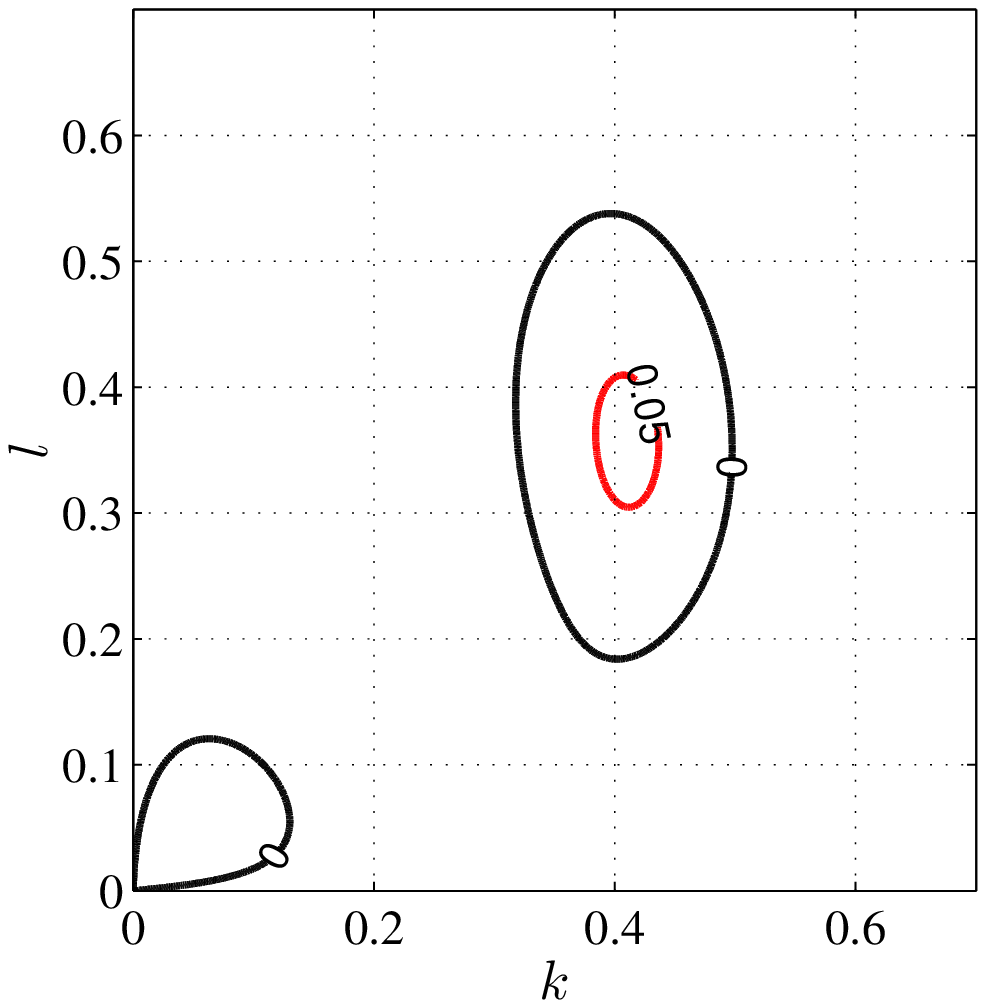}}
\label{fig:st2}}
}
\mbox{
\vspace{-0.7in}
\hspace{-0.2in}
\subfigure[  ]{
\scalebox{0.33}{\includegraphics*[viewport=1in 0in 5in 4.1in, width=5in, height=4.5in]{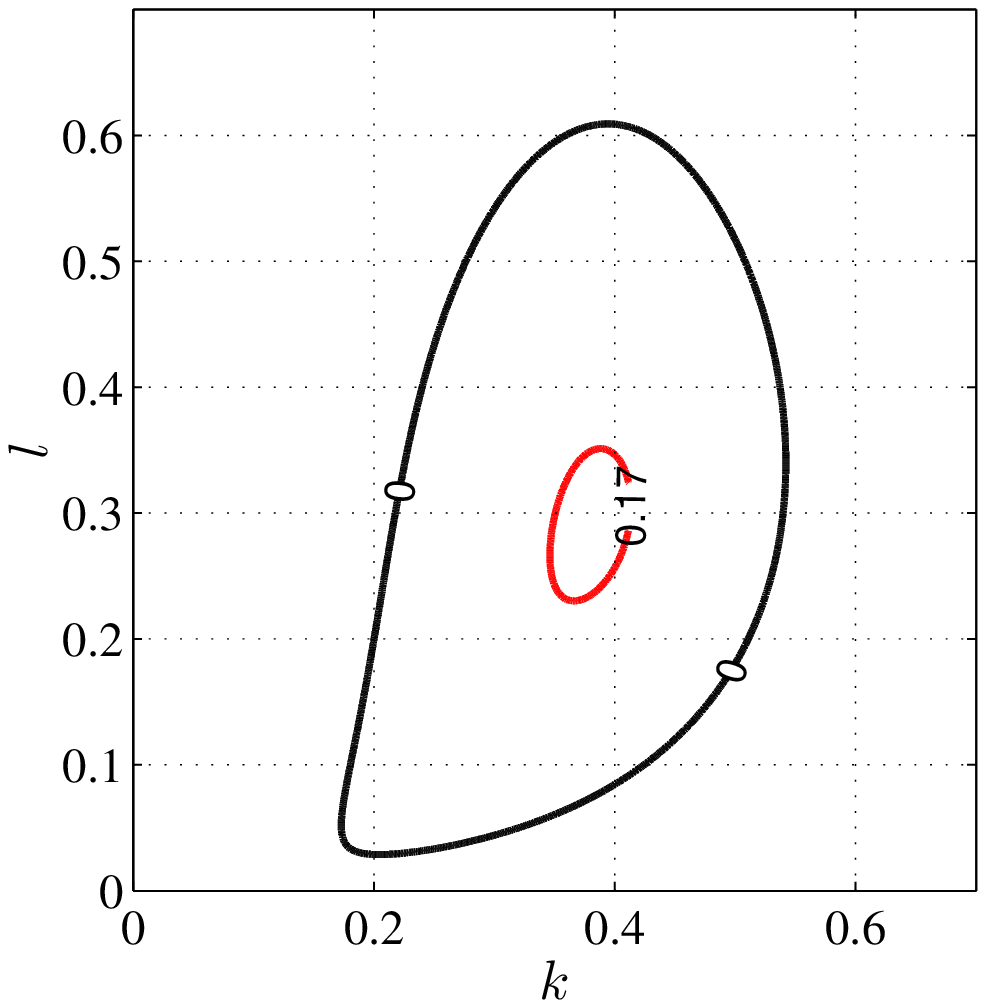}}
\label{fig:st3}}
\subfigure[  ]{
\hspace{-0.2in}
\scalebox{0.33}{\includegraphics*[viewport=1in 0in 5in 4.1in, width=5in, height=4.5in]{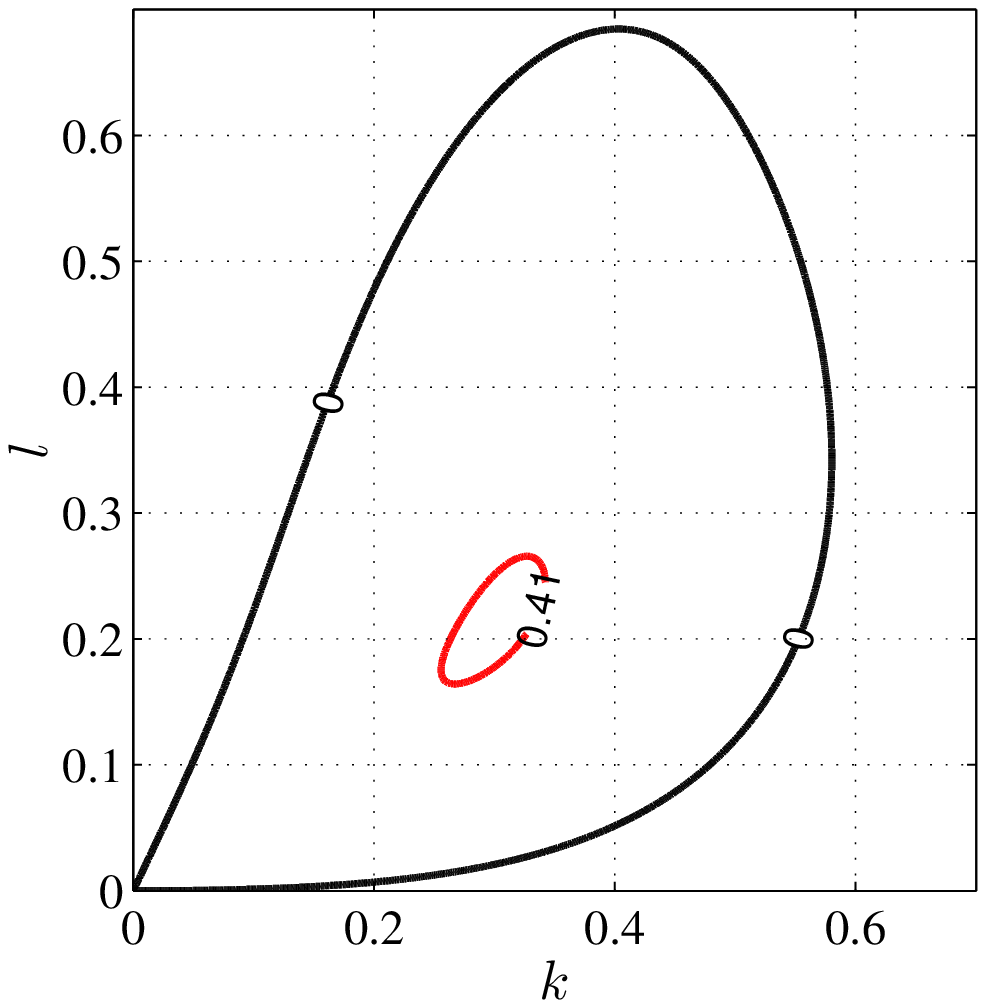}}
\label{fig:st4}}
}

\caption{(Color online) Growth rates of perturbations at selected $(\mu,q)$ indicated in figure \ref{fig:14}, all with $\mu=0.1$. (a) $q=-0.05$, where  stripes are OSV unstable. (b) $q=-0.02$, where stripes are CR and OSV unstable. (c) $q=0.0058$, where stripes are CR unstable but OSV stable. (d) $q=0.05$,  where stripes are CR and skew-varicose unstable, though the distinction between these two instabilities has become blurred. The zero contour of the real part of the eigenvalue is denoted in black (outer) and contours of  positive growth rate are in red (inner in gray).  The eigenvalues in the OSV case are complex.}
 \label{fig:sc}
 \end{figure}
 
\begin{figure}[t]
\scalebox{0.36}{\includegraphics*[viewport=0.5in 0.8in 9.6in 6.6in]{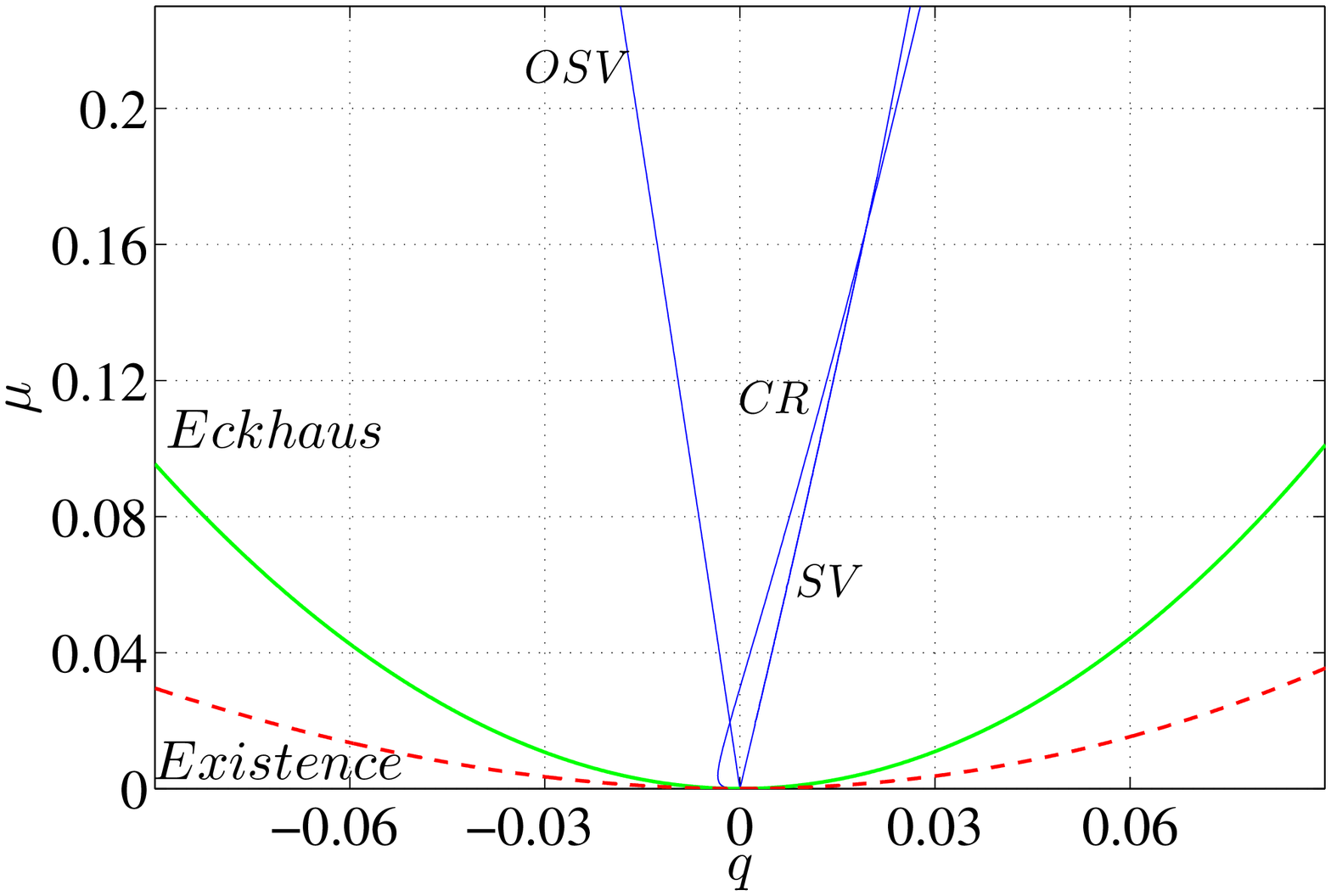}}
\caption{(Color online) Stability diagram for model 1 with $c=0$ (stress-free boundary conditions), $Pr=1$, $g_{m}=50$ and $\gamma=2.5$. Stable stripes exist for  $\mu>0.02$; the region of stable stripes is bounded by the OSV boundary on the left and the CR boundary and then the SV on the right. }
 \label{fig:16}
\end{figure}

\begin{figure}[!h]
\mbox{
\vspace{-0.2in}
\hspace{-0.2in}
\subfigure[  ]{
\scalebox{0.27}{\includegraphics*[viewport=0.4in 0.7in 7in 6.6in,width=6in, height=7in]{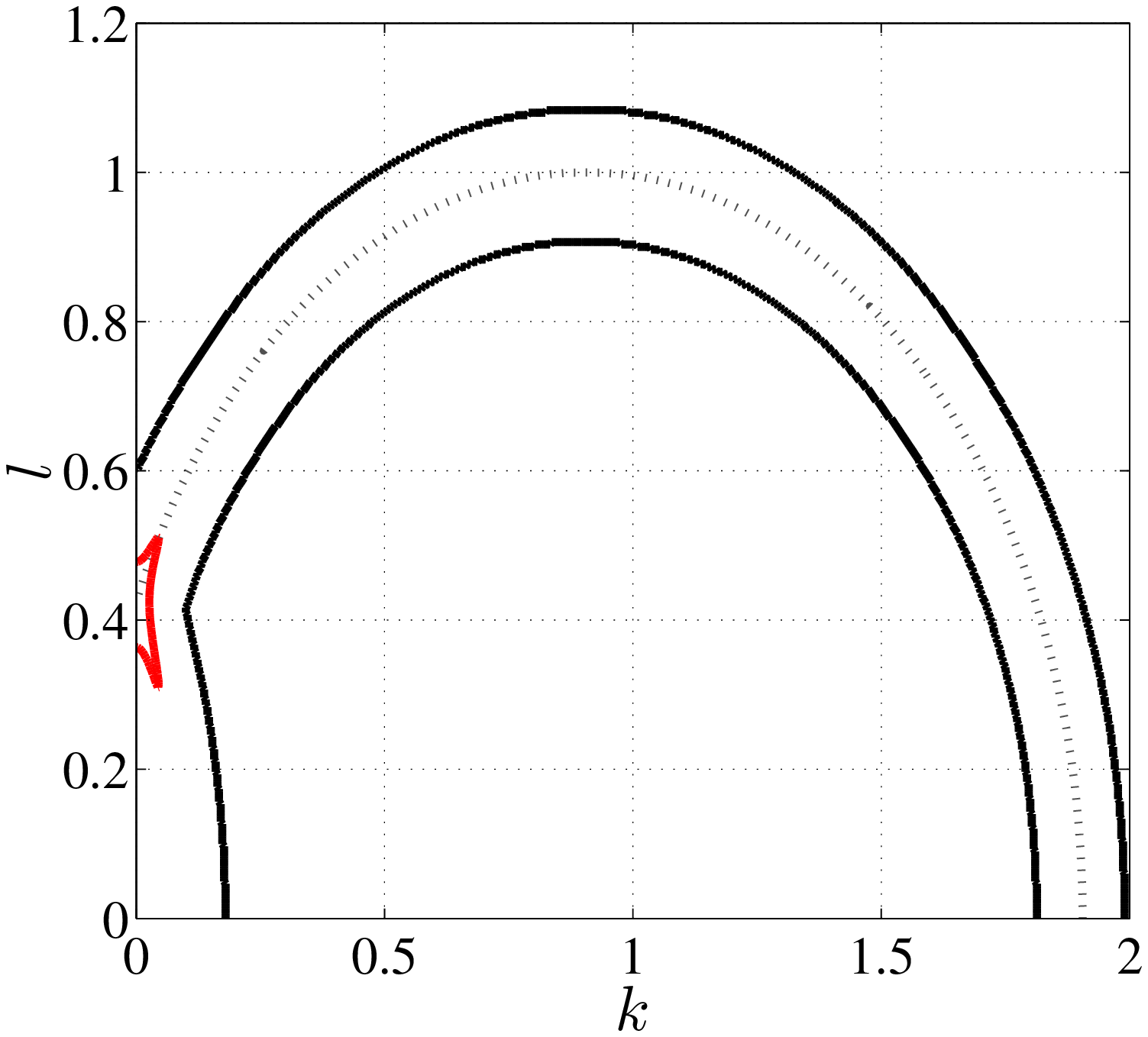}}
\label{fig:17st1}}
\subfigure[  ]{
\hspace{-0.15in}
\scalebox{0.27}{\includegraphics*[viewport=0.4in 0.7in 7in 6.6in,width=6in,  height=7in]{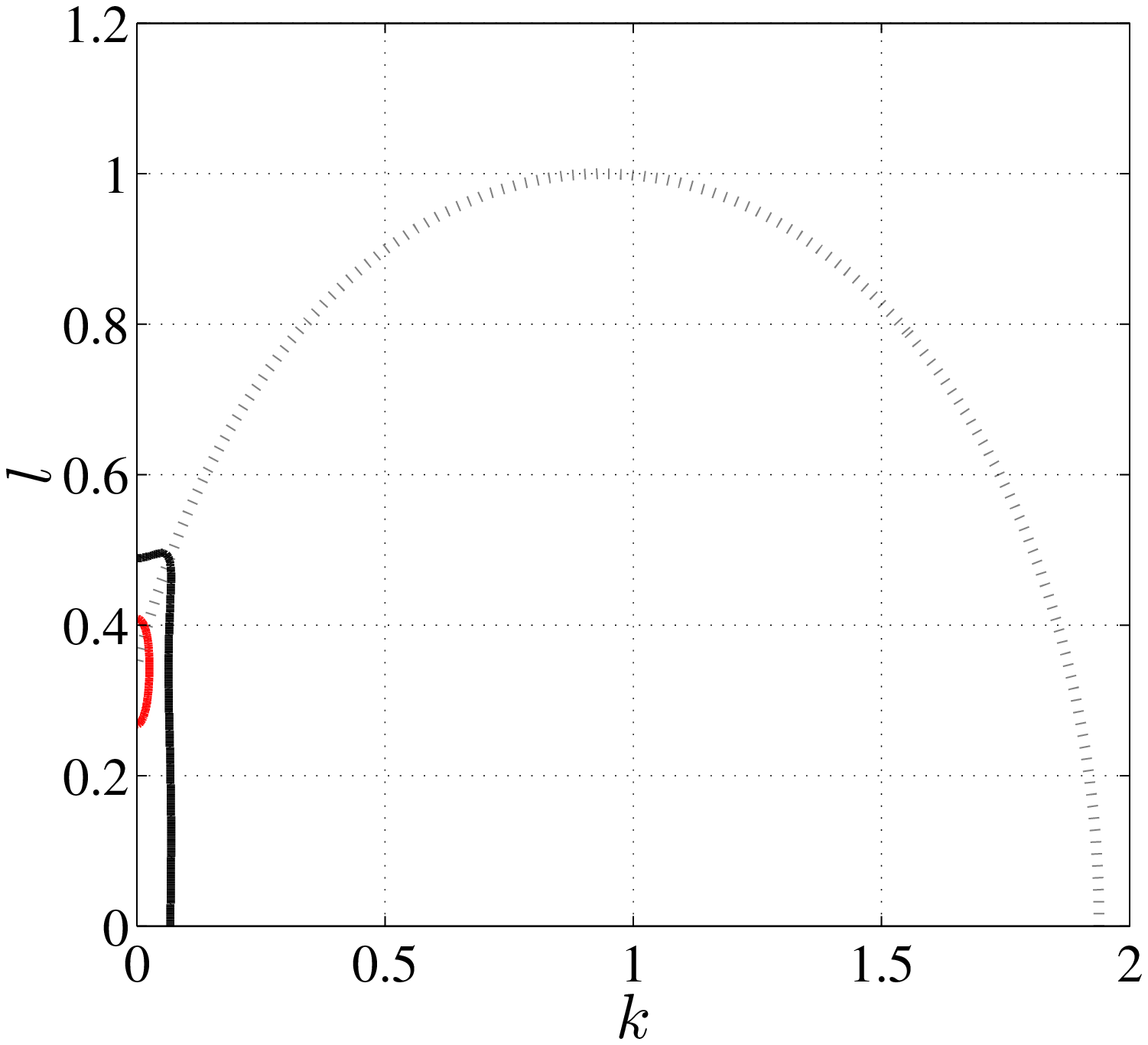}}
\label{fig:17st2}}
}
\mbox{
\vspace{-0.7in}
\hspace{-0.2in}
\subfigure[  ]{
\scalebox{0.27}{\includegraphics*[viewport=0.4in 0.7in 7in 6.6in,width=6in,  height=7in]{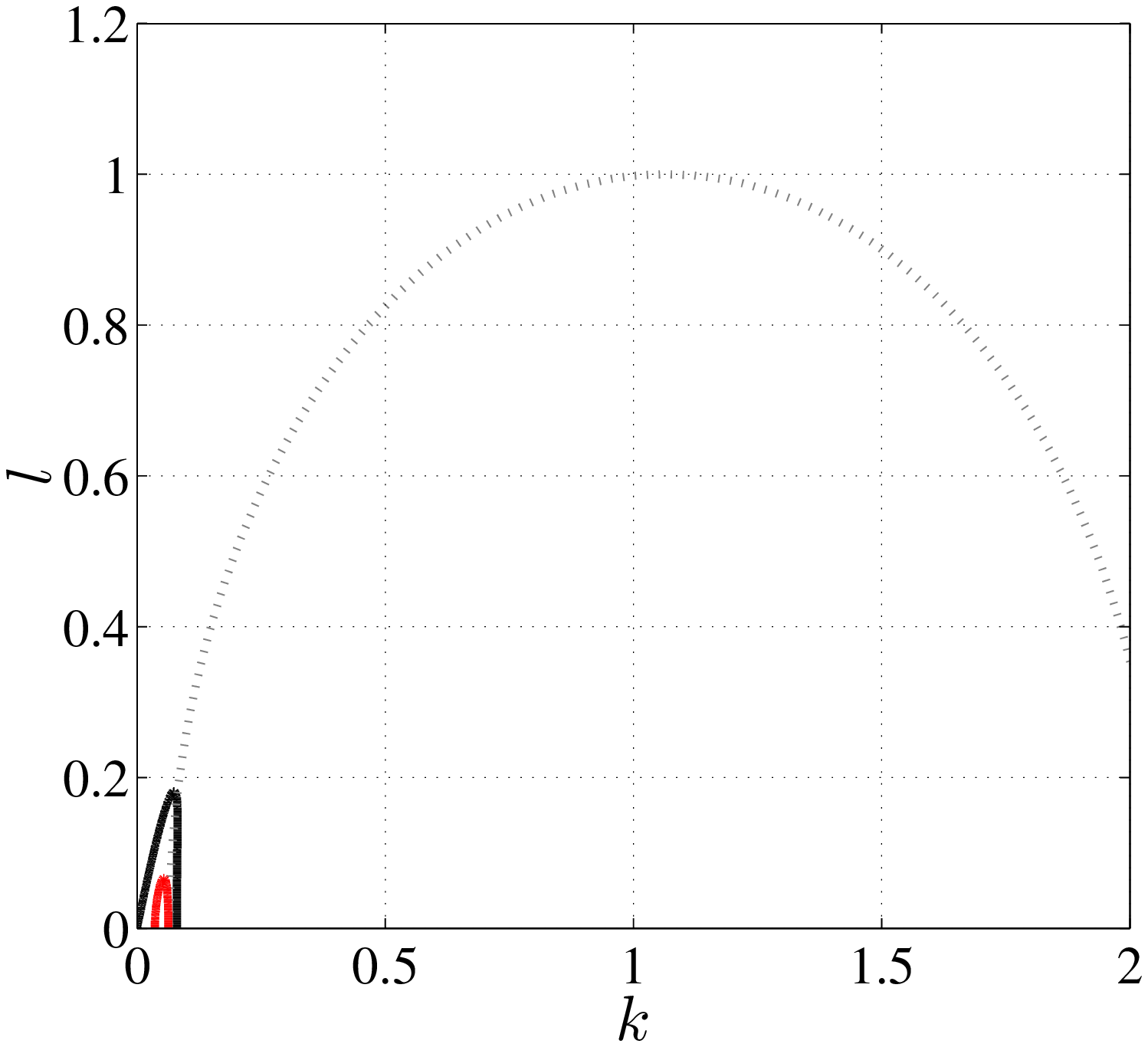}}
\label{fig:17st3}}
\subfigure[  ]{
\hspace{-0.15in}
\scalebox{0.27}{\includegraphics*[viewport=0.4in 0.7in 7in 6.6in,width=6in, height=7in]{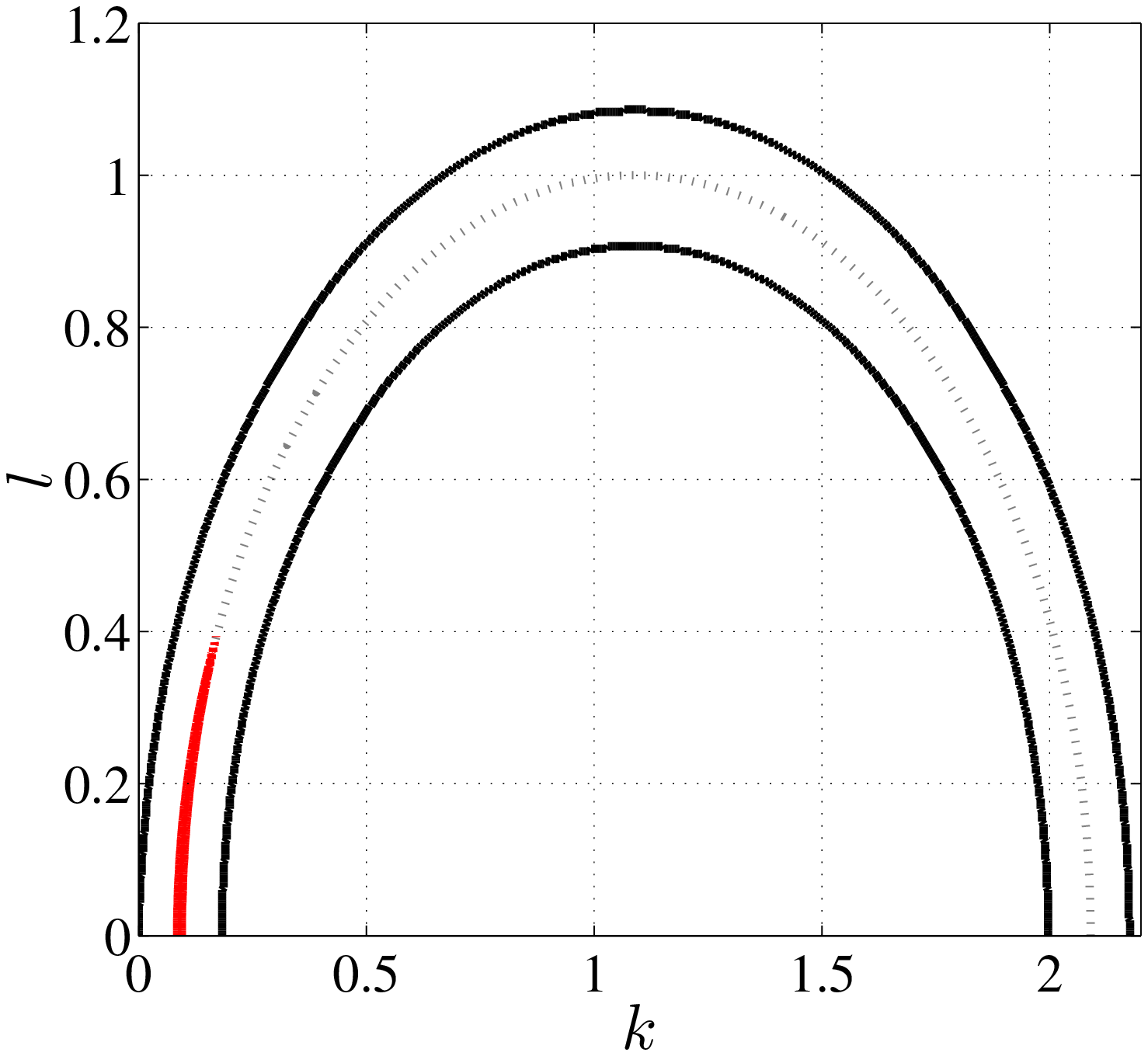}}
\label{fig:17st4}}
}
\caption{(Color online) Growth rates of perturbations as a function of $(k ,l)$ for the Swift--Hohenberg equation~(\ref{eq:she}), with $\psi^3$ replaced by $P_{\alpha}(\psi^3$). Figures (a) and (b) correspond to points in the left Eckhaus band, $q<0$, whereas  (c) and (d) correspond to points in the right Eckhaus band, $q>0$. (a) Stripes are zigzag and Eckhaus unstable and the unstable modes lie approximately in a annulus of unit radius centered at $(k,l)=(1+q,0)$ have positive real eigenvalues. (b) Stripes are zigzag and Eckhaus unstable. (c) Stripes are Eckhaus unstable.  (d) Stripes are Eckhaus unstable and in addition, an annulus of unit radius centered at $(k,l)=(1+q,0)$ has positive real eigenvalues. Thick black: zero contour. thick Red (inner in gray): positive contour. Dotted curve: circle of unit radius centered at $(k,l)=(1+q,0)$. }
 \label{fig:strangec}
\end{figure}
\begin{figure}[t]
\vspace{-0.1in}
\mbox{
\vspace{-1.5in}
\subfigure[  ]{
\scalebox{0.36}{\includegraphics*[viewport=0.5in 0.8in 9.6in 6.3in]{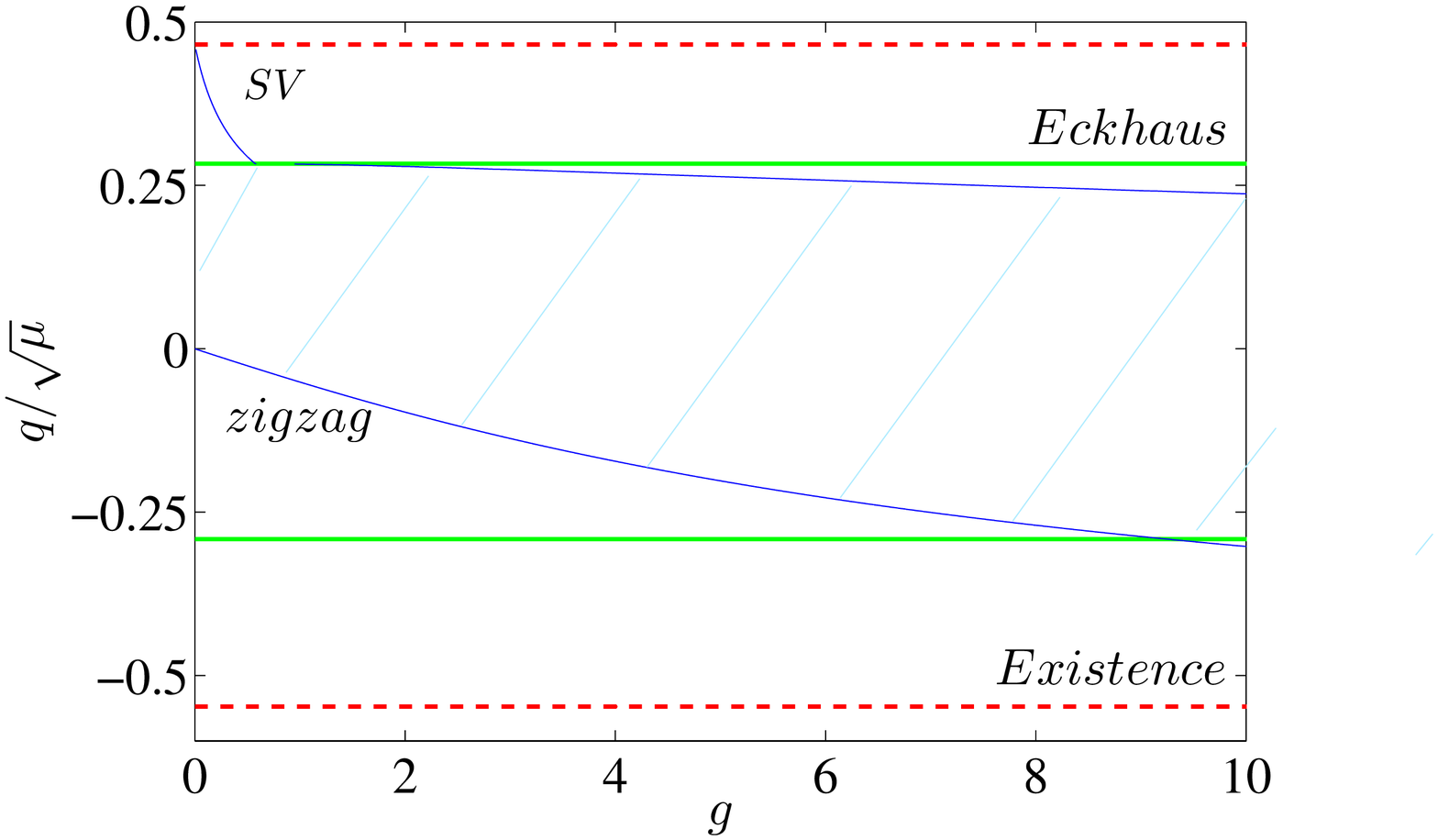}}
\label{fig:g1}}}
\mbox{
\subfigure[  ]{
\scalebox{0.36}{\includegraphics*[viewport=0.5in 0.4in 9.6in 6.3in]{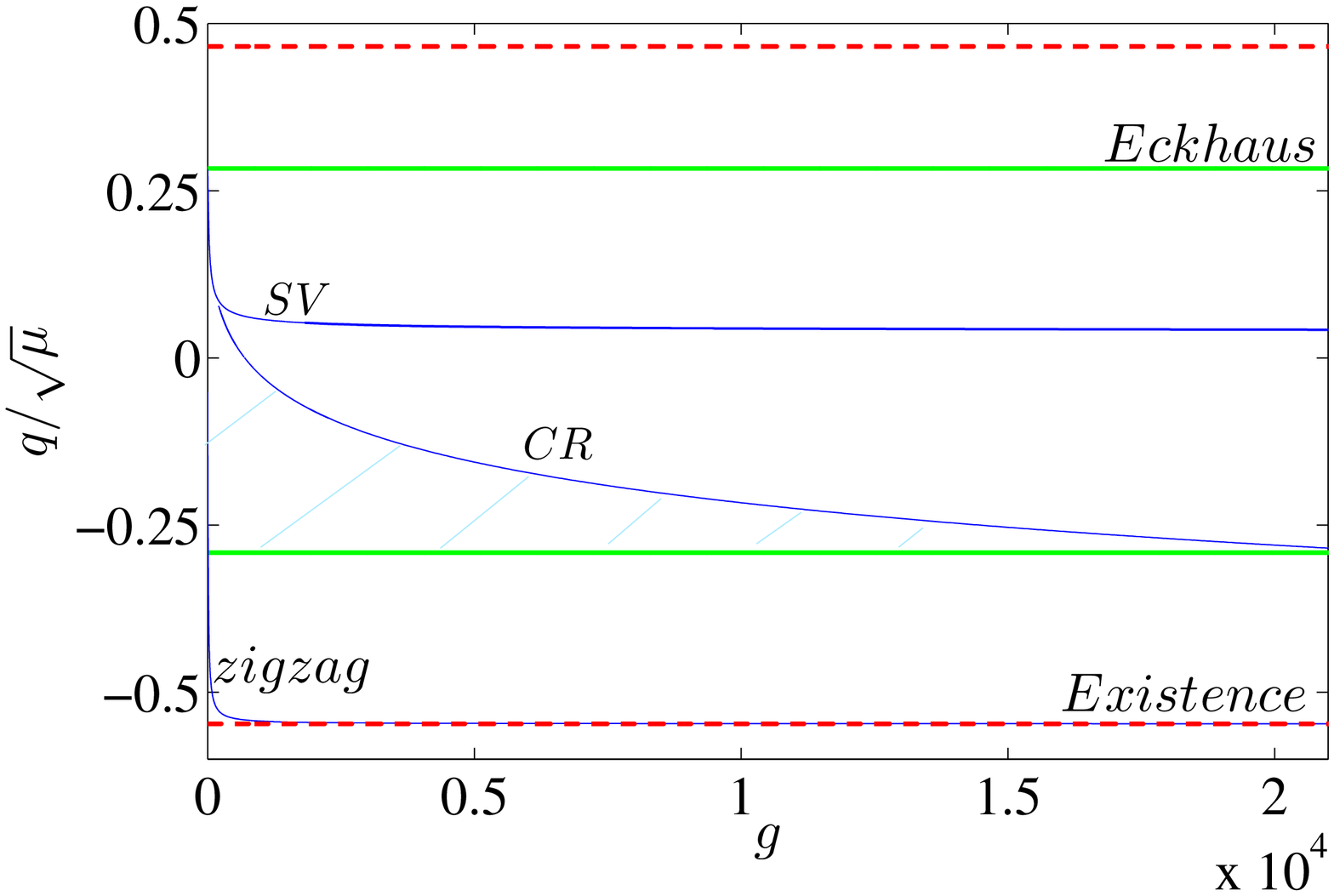}}
\label{fig:g2}}
}
\caption{(Color online) Stability diagrams in $(g,q/\sqrt{\mu})$ plane for model 2 with  $\mu=0.1$ and $\gamma=2.5$. (a) For small $g$: the region of stable stripes is mainly bounded by zigzag and SV instabilities.  (b) For large $g$, the region of stable stripes is bounded by the Eckhaus (thick green)  and CR instabilities. These instability boundaries cross around $g=2\times10^{4}$, eliminating the region of stable stripes. The region of stable stripes is hatched.}
 \label{fig:gq}
\end{figure}
\begin{figure*}[t]
\mbox{
\vspace{-0.8in}
\subfigure[  ]{
\hspace{-0.4in}
\scalebox{0.28}{\includegraphics*[viewport=0.5in 0.37in 9.6in 6.3in]{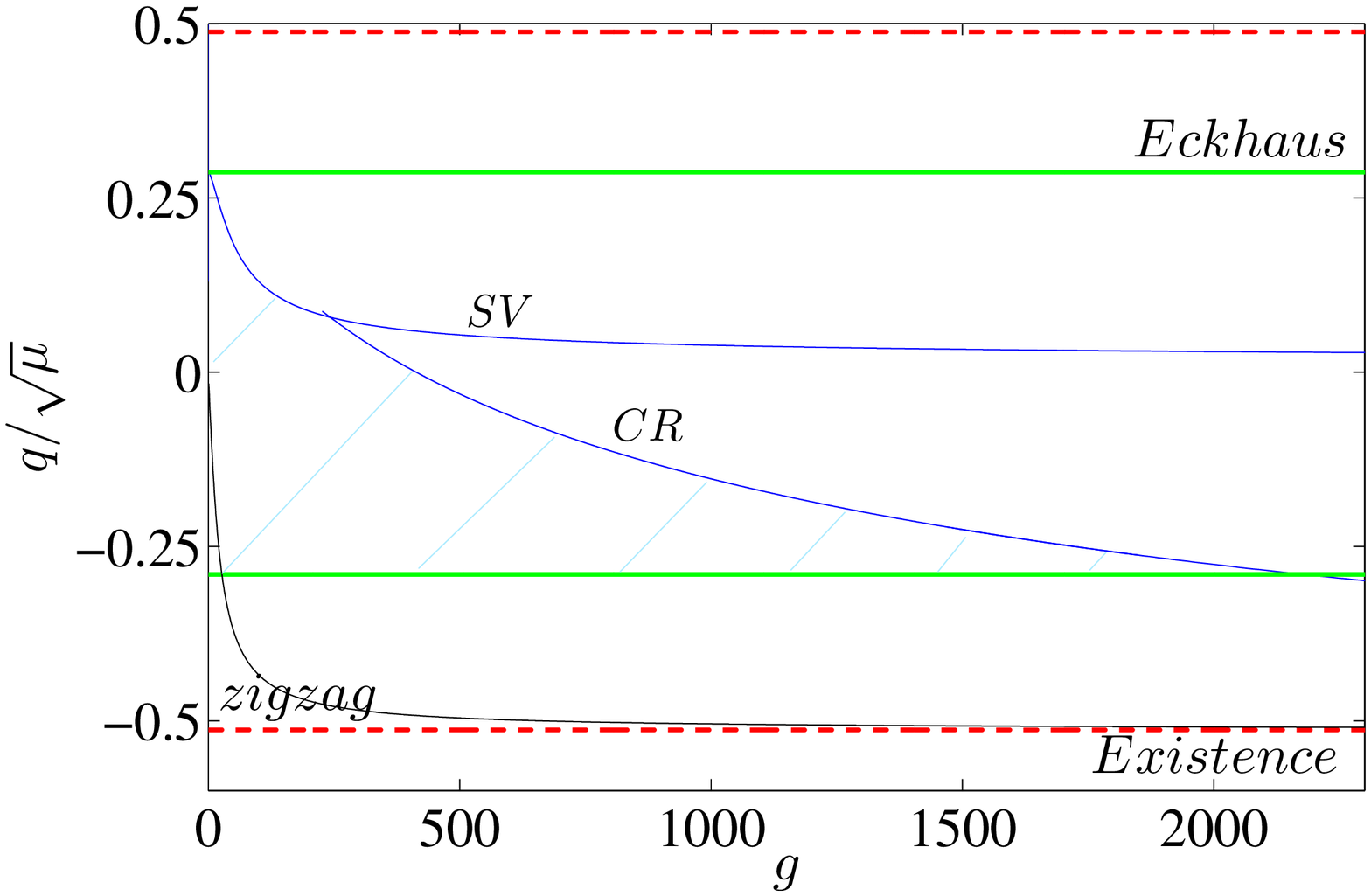}}
\label{fig:st1}}
\subfigure[  ]{
\hspace{-0.2in}
\scalebox{0.28}{\includegraphics*[viewport=0.5in 0.37in 9.6in 6.3in]{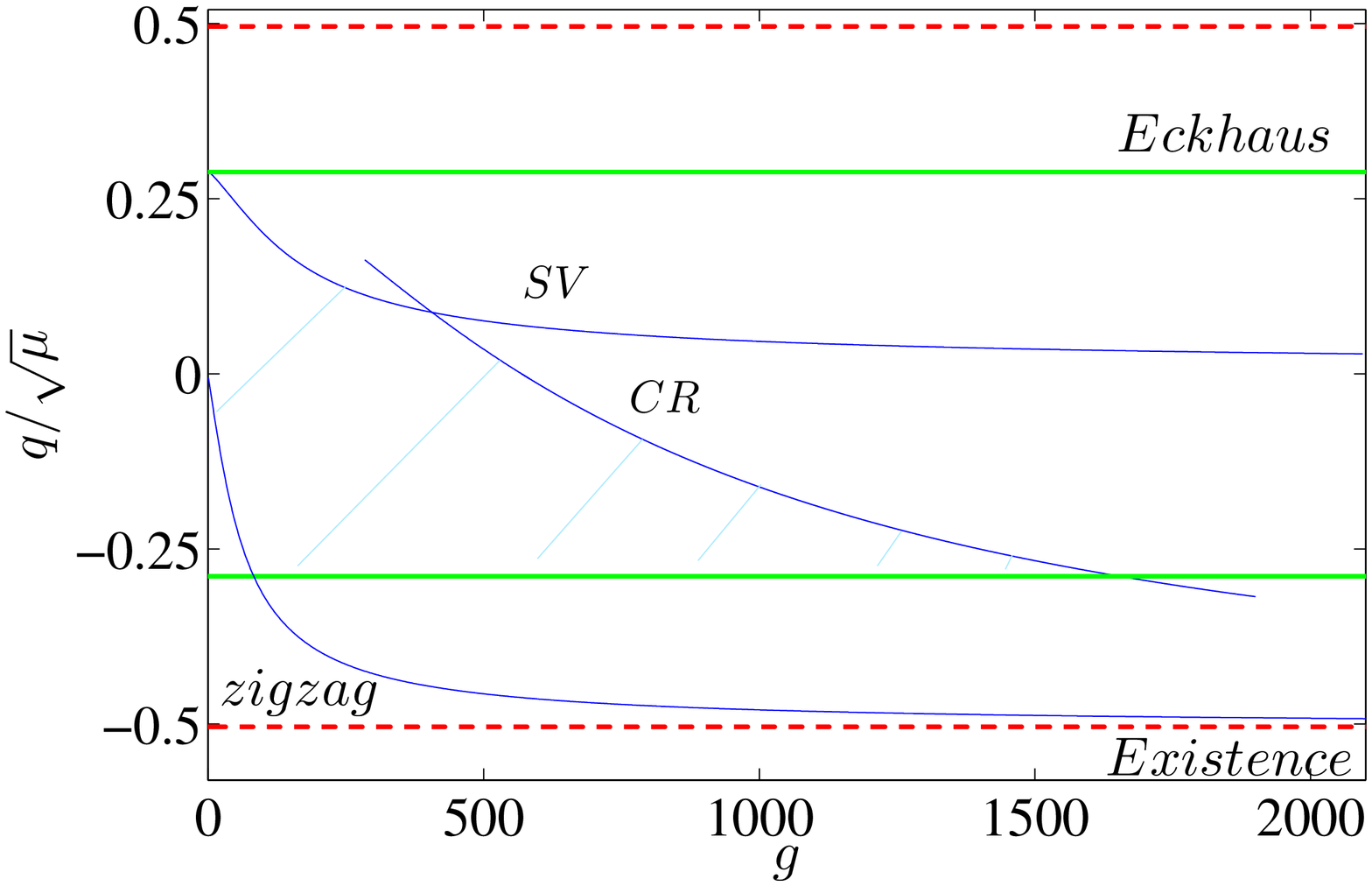}}
\label{fig:st2}}
\subfigure[  ]{
\hspace{-0.3in}
\scalebox{0.28}{\includegraphics*[viewport=0.5in 0.4in 9.6in 6.3in]{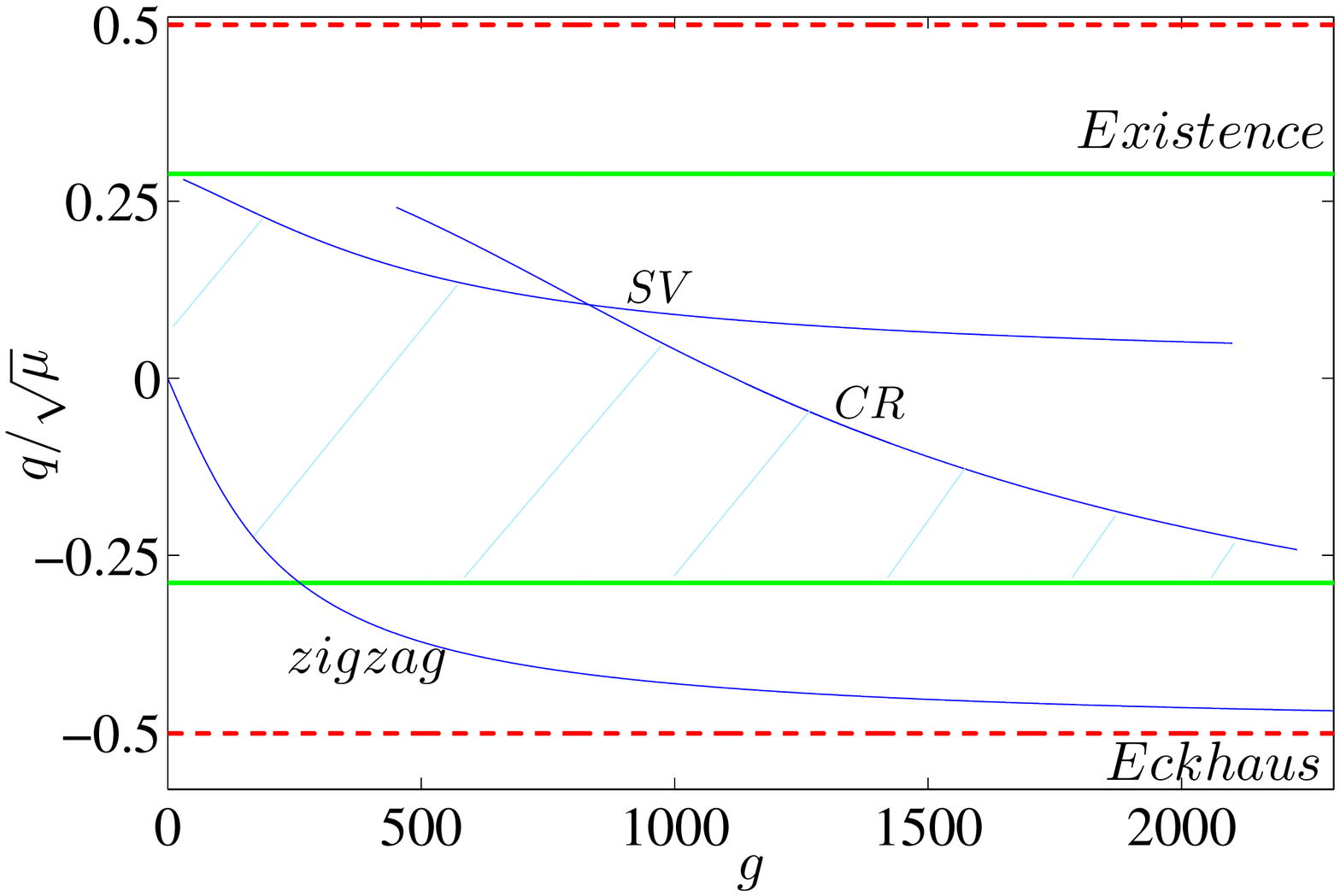}}
\label{fig:st3}}
}
\caption{(Color online) Stability diagrams in $(g,q/\sqrt{\mu})$ plane for model 2 with $\gamma=2.5$. (a) $\mu=0.01$ (b) $\mu=0.001$ (c) $\mu=0.0001$. In all three cases CR instability reduces the region of stable stripes which is  hatched. }
 \label{fig:19}
\end{figure*}
For the same parameter values, the stability diagram for a larger range of $q$ and $\mu$ is shown in figure~\ref{fig:rwc}. The region of stable stripes is bounded by the Eckhaus instability from below and the CR instability from above. The zigzag instability boundary lies close to the existence curve and is of less interest for this  value of $g_{m}$.

Figure~\ref{fig:13} shows the eigenvalue behavior at selected points in the $(q,\mu)$ space from figure 11, as a function of $(k,l)$, showing how stripes can be unstable to one or both of the SV and CR instabilities. We note that the CR instability occurs for reasonably large values of $k\approx 0.04$  and $l\approx 0.2$. In contrast, in the SVI, contours of positive growth rate emerge from $(k,l)=(0,0)$. When both instabilities exist, two separate peaks of growth  rates appear. For large $q$, these contours can join to form  one large contour.

We have computed the stability diagrams for model 2 with  $g=500$ and $g=25$, and these are qualitatively the same as~\ref{fig:sv50} and \ref{fig:zoom}, consistent with the relation  $g=g_{m}/Pr\,c^2$. 

Figure~\ref{fig:14} and~\ref{fig:16} similarly show the instability boundaries for model 1, with parameters $c=0$ (corresponding to stress-free boundary conditions), $Pr=1$, $\gamma=2.5$ and $g_m=50$ and $g_m=1000$. The  results  agree qualitatively with earlier calculations by Bernoff \cite{Bernoff:1994}, who found a similar linear relation between  $\mu$ and $q$ for the SV and OSV instabilities in convection with stress-free boundary conditions; he did not consider the CR instability. 

Disregarding the CR instability, stripes would be stable between the OSV instability and SVI boundaries. However, as seen in figure~\ref{fig:14}, the SVI is  always preempted by the CR instability; these boundaries appear to be parallel for  larger $\mu$. For $\mu<0.32$, there are no stable stripes.  For higher $\mu$, the stable region is  bounded by the CR instability from above and by the OSV instability from below. 

Figure~\ref{fig:sc} shows the change of structure of  the eigenvalues when moving from left to right in the stability diagram shown in figure~\ref{fig:14}. At $\mu=0.1$, we selected four different wavenumbers: $q=-0.05$, where  stripes are OSV unstable, $q=-0.02$, where stripes are CR and OSV unstable, $q=0.0058$, where stripes are CR unstable but OSV stable and $q=0.05$,  where stripes are CR and skew-varicose unstable, though the distinction between these two instabilities has become blurred.

Figure~\ref{fig:16} presents instability boundaries for stress-free boundary conditions with $Pr=1$, $g_m=50$ and $\gamma=2.5$. This provides an illustration of the change of the CR instability boundary with $g_{m}$. For $g_{m}=50$, the CR instability boundary crosses the SVI boundary and the effect of the CR instability is reduced.

\subsection{Curious behavior of growth rates }
In previous sections, we considered the growth rate behavior of $(k,l)$ close to  zero. However, when a larger range of $k$ and $l$ is considered, we noticed additional instabilities in the regime where stripes are already unstable. These instabilities occur even in the SHE and so are not related to presence of the mean-flow. They have not been studied before, though they are of less interest since they do not bound the region of stable stripes. We consider them briefly here.

Contours of the growth rates of perturbations for selected parameter values are shown in figure~\ref{fig:strangec}. As shown in figure~\ref{fig:17st1}, for a fixed small $\mu$, when $q$ is increased from the left existence boundary $(q<0)$, modes approximately in an annulus of unit radius centered at $(k,l)=(1+q,0)$ have positive real eigenvalues in addition to unstable modes for small $k$ and $l$, corresponding to the Eckhaus and zigzag instabilities. This annulus disappears for larger $q$  leaving only the Eckhaus and zigzag unstable modes close to $(k,l)=(0,0)$, as  shown in figure~\ref{fig:17st2}. This process is in part reversed for $q>0$ going towards the right existence boundary: here modes are Eckhaus unstable only, but for large enough $q$, the annulus of unstable modes reappears. The behavior of growth rates in the right Eckhaus band is illustrated in figure~\ref{fig:17st3}  and \ref{fig:17st4}.

We have confirmed that this curious behaviour is not a result of using the projection operator $P_{\alpha}$ in the Swift--Hohenberg equation.

\section{Analysis of the role of mean-flows}

We now illustrate the stability diagrams in a different manner. For a fixed $\mu$, we show how the coupling to the mean flow affects the region of stable stripes. This choice of presentation provides useful information for numerical simulations of the PDEs in large domains.
Figure~\ref{fig:gq} represents the region of stable stripes for model 2 in the $(g,q/\sqrt{\mu})$ plane for $\mu=0.1$. We choose $q/\sqrt{\mu}$ as the coordinate for ease of comparison between different values of $\mu$. Figure \ref{fig:g1} shows the stability diagram for small $g$, where the region of stable  stripes is bounded from above by the Eckhaus instability  for $g<0.574$  and by the SVI  for $g>0.574$, and by the   
zigzag instability from below. Figure \ref{fig:g2} shows how for large $g$, the region of stable stripes is bounded by the CR instability  from above and by the Eckhaus instability from below, and the region is eliminated when $g\gtrsim2\times10^{4}$. Figure~\ref{fig:19} shows the location of the CR instability  for $\mu=0.01,0.001$ and $0.0001$. The upper bound on $g$ beyond which there are no stable stripes initially decreases with $\mu$ and then increases. The behavior of model 1 is qualitatively the same.

\begin{figure*}[t]
\mbox{
\vspace{-0.8in}
\subfigure[  ]{
\hspace{-0.4in}
\scalebox{0.28}{\includegraphics*[viewport=0.5in 0.37in 9.6in 6.3in]{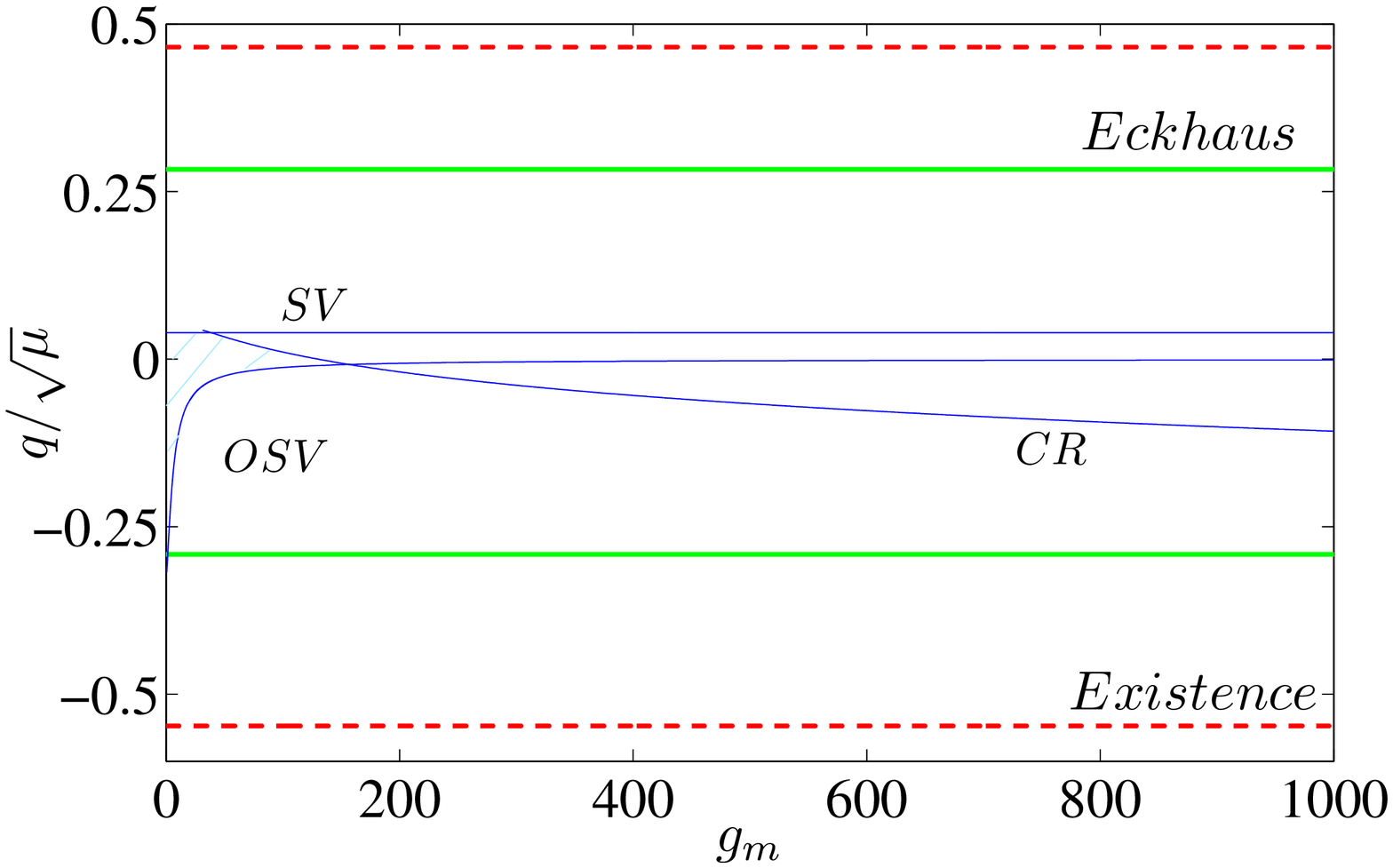}}
\label{fig:20st1}}
\subfigure[  ]{
\hspace{-0.2in}
\scalebox{0.28}{\includegraphics*[viewport=0.5in 0.37in 9.7in 6.3in]{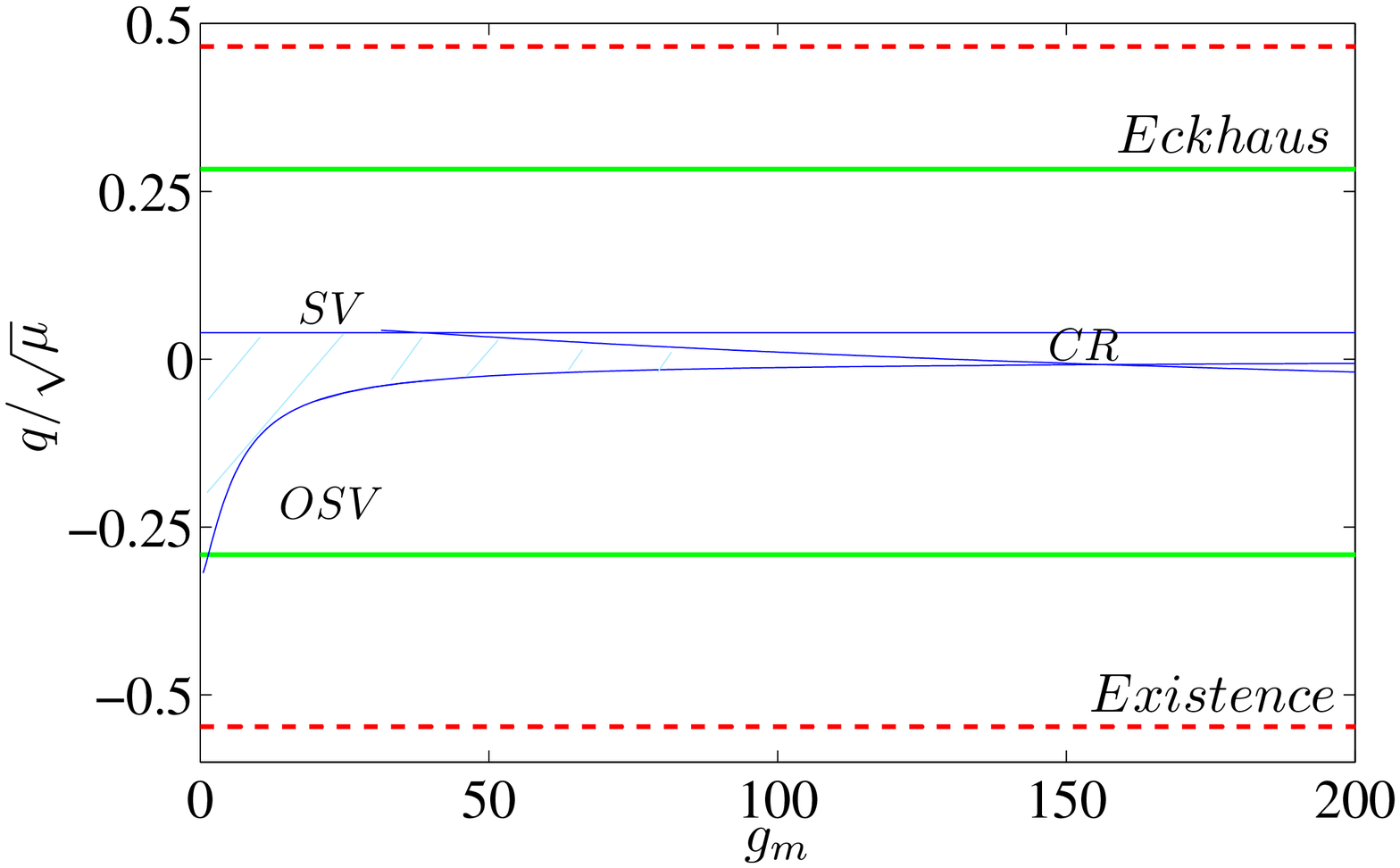}}
\label{fig:20st2}}
\subfigure[  ]{
\hspace{-0.3in}
\scalebox{0.28}{\includegraphics*[viewport=0.5in 0.4in 9.6in 6.3in]{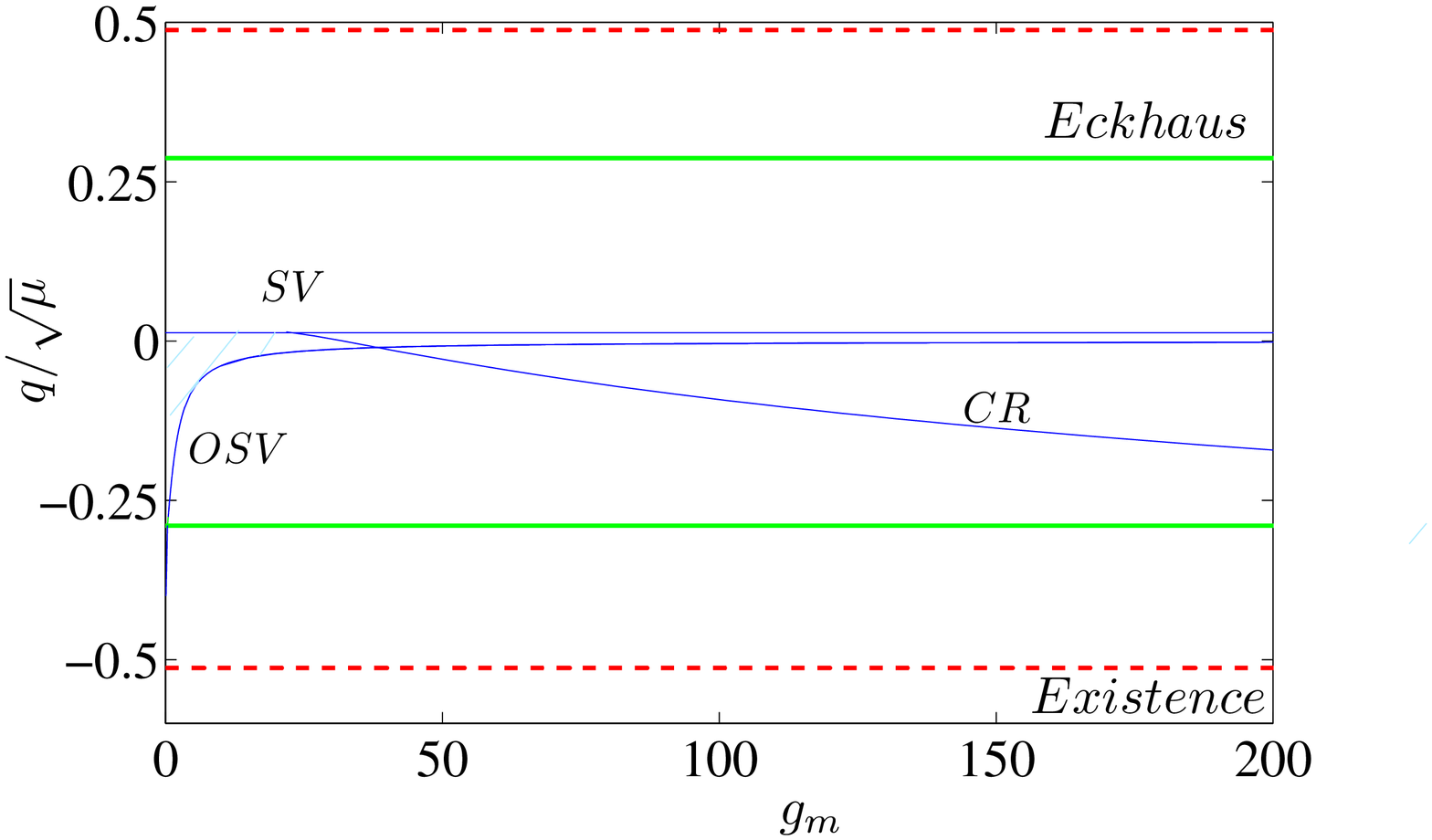}}
\label{fig:20st1}}
}
\caption{(Color online) Stability diagrams in $(g_{m},q/\sqrt{\mu})$ plane with $c=0$ (stress-free boundary conditions), $Pr=1$  and $\gamma=2.5$. (a) $\mu=0.1$ and for large $g_{m}$: (b) $\mu=0.1$ and for small $g_{m}$: stable stripes are completely  eliminated  when $g_{m}\gtrsim130$. (c) $\mu=0.01$ and for small $g_{m}$. In all three cases, the region of stable stripes (hatched ) is mainly bounded by the OSV, SV and CR instabilities, and the CR instability makes the upper bound in $g_{m}$ and reduces the region of stable stripes with $\mu$. }
 \label{fig:gqsf}
\end{figure*}

Figure~\ref{fig:gqsf} presents stability diagrams for model 1 with $c=0$ (stress-free boundary conditions), for $\mu=0.1$ and $\mu=0.01$. The SV and the OSV instabilities bound the region of stable rolls from above and below and the CR instability makes the upper bound on $g_{m}$. Stable stripes exists only for small $g_{m}$ and the upper bound on $g_{m}$ further reduces with $\mu$.

\section {conclusion}
In this work we considered the consequences of including the mean-flow on the generalized Swift--Hohenberg models. We analysed two models:  
in the first model vorticity has  its own independent dynamics~\cite{Xi:1993}. In the second, vorticity is directly slaved to the order parameter~\cite{Green:1985}. These two models  are related to each other through $g=g_{m}/(Pr\ c^{2})$, where $g_{m}$  and $g$ are the couplings to mean flow in model 1 and model 2 respectively.
Two boundary conditions were considered in this work: stress-free ($c=0$) and  no-slip ($c^2=2$). 

In order to explore long-wavelength instabilities, we carried out a complete  linear stability analysis of stripes. We expressed the relevant determinants as power series in $k^{2}$ and $l^{2}$, where $(k,l)$ is the perturbation wavevector. We were able to derive explicit expressions for the largest growth rates in most cases. This has led to an improved understanding of the instabilities of stripes. Unlike in previous work \cite{Hoyle:2006,ponty:1997,Bernoff:1994}, we have not had to make assumptions on the relation between $k$, $l$ and the amplitude of the basic stripe solution. This approach has been made possible through the use of the projection operator, $P_{\alpha}$, which allows the exact stripe solution to be written down easily~\cite{Ian}.

The skew-varicose instability, which was our main concern, has two different behaviors:    
if stripes are stable to the Eckhaus instability, in the limit of $\mu=0$, the SVI goes as $\mu\sim12q^2$, provided $g>0.75$. The most unstable wavevector satisfies $k^{2}/l^{2}=\mathcal{O}(1)$.  For $g<0.75$, the SVI boundary crosses the Eckhaus curve, and  in the limit of $\mu=0$, it goes as  $\mu\sim aq^2$ with $4<a<12$. In model 1, the critical $g_{m}$ is $0.75Pr\,c^{2}$.
In the large $g$ limit (that is, for very low $Pr$, or for stress-free boundary conditions), there is a transition of the SVI boundary from $\mu=12q^2$ to $\mu=8q$ at a wavenumber satisfying $q\propto1/g$. 

An additional instability, the oscillatory skew-varicose (OSV) instability, is encountered  for stress-free boundary conditions in model 1.
The OSV instability boundary is approximately $\mu=\left(\frac{-3+\sqrt{5}}{3}\right)qg_{m}$, for small $\mu$.

We confirmed our  analytical results by numerical computations of the eigenvalues of the stability matrices. These eigenvalues also allow us to explore short-wavelength instabilities: cross-roll and the oscillatory instability. 
Finally, we have shown how the region of stability of stripes is  eliminated for small $\mu$ and large enough $g$.

The use of the projection operator $P_{\alpha}$, which is equivalent to a truncation to selected wave numbers, made this analysis straightforward and allowed the complete understanding of the skew-varicose instability in our models. Numerical simulations of these projected models for small $\mu$  have qualitatively the same solutions as the unprojected PDEs (for example, they have SDC solutions); this is our justifcation for using these projected models in the stability analysis for small $\mu$. The projected and unprojected models will of course differ for large $\mu$. It is of interest to address the question whether a similar projection could be used in the analysis of the Navier–-Stokes equations.

We finally comment on the implications of our analysis on the direct numerical simulations of PDE's in large domains. The most striking signature of the inclusion of a mean-flow is the existence of the skew-varicose instability, which can play an important role in the formation of the spiral defect chaos or defect chaos~\cite{Bodenschatz:2000,Busse:1978}.  Hence the improved understanding of the stability of stripes in this work provides a foundation for numerical simulations of spiral defect chaos and  defect chaos. The SDC state exists inside the Busse balloon, where convection rolls are stable~\cite{Morris:1993}. Thus we intend in a future paper to relate the SDC and defect chaotic states present in these generalized SH models to calculations carried out using Rayleigh--B{\'e}nard convection  with stress-free~\cite{Morris:1993} and no-slip~\cite{Chiam2003b} boundary conditions, aiming to improve the understanding of why SDC occurs in convection. The results of this work provide useful information for the choice of parameters for different instability regimes in  model 1 and model 2. In particular we will be interested in numerical simulations with small $\mu$ and in large domains. This work also justifies using  model 2 with large $g$, since this has been shown to have similar stability properties as to model 1 with small $Pr$.

\subsection*{Appendix I}
   We compute the stability boundaries by defining conditions on the maximum eigenvalue, $\sigma_{max}$ and using the continuation package MATCONT~\cite{Dhooge:2004}. 

We calculate eigenvalues of the matrices $J_{1}$ and $J_{2}$  and determine numerically how these depend on $k$ and $l$ in order to compute the derivatives that are needed.

\begin{itemize}
\item  Stripes are unstable to the Eckhaus instability if $\sigma_{max}$ is positive for $l=0$ and $k=k_{max}>0$. The Eckhaus stability boundary occurs when  $k_{max}\rightarrow0$, in which case $\frac{\partial^{2}\sigma_{max}}{\partial k^{2}}=0$ at $(k,l)=(0,0)$.

\item  Stripes are unstable to the zigzag instability if $\sigma_{max}$ is positive for $k=0$ and $l=l_{max}>0$. The zigzag stability boundary occurs when  $l_{max}\rightarrow0$, in which case $\frac{\partial^{2}\sigma_{max}}{\partial l^{2}}=0$ at $(k,l)=(0,0)$.

\item There are two cases of the SVI, as discussed in the text. First if the SVI precedes the Eckhaus instability, we consider $\sigma_{max}$ as a function of $k$ and $l$ in polar coordinate, so $k=\epsilon \cos(\theta)$ and $l=\epsilon \sin(\theta)$. Stripes are unstable to the SVI, if $\sigma_{max}$ is positive for some $\theta$ in the limit of $\epsilon=0$.  The SV stability boundary occurs when $\epsilon\rightarrow0$, in which case $\sigma_{max}=0$ and $\frac{\partial\sigma_{max}}{\partial \theta}=0$. 
Second if the SVI follows the Eckhaus instability, stripes are unstable to the SVI, if $\sigma_{max}$ is positive for  sufficiently small $l=l_{max}$ and $k=k_{max}$.  The SV stability boundary occurs when $l_{max}\rightarrow0$, in which case $\frac{\partial\sigma_{max}}{\partial k^{2}}=0$ and $\frac{\partial\sigma_{max}}{\partial l^{2}}=0$ at  $(k,l)=(k_{max},0)$.
\item In the case of CR instability, $\sigma_{max}=0$ occurs at non-zero $(k,l)$. Hence the three conditions are  $\frac{\partial\sigma_{max}}{\partial k}=0$, $\frac{\partial\sigma_{max}}{\partial l}=0$ and ${\sigma_{max}}=0$ at  $k=k_{max}\neq0$ and $l=l_{max}\neq0$.
\end{itemize}

The SV and CR instabilities can be  oscillatory. In these cases, we consider  ${\sigma_{max}}$ to be the real part of the eigenvalue and use the same relevant conditions stated above.
\subsection*{Acknowledgements}
We are grateful to  Prof. Ian Melbourne for suggesting the  method of using the projection operator $P_{\alpha}$ to allow an exact stripe solution, which enables the calculations to proceed without   making assumptions on asymptotic  relations between the perturbation wavevector and the amplitude of the stripes. JAW is grateful for the support from the Overseas Research Scholarship Scheme and from the School of Mathematics, University of Leeds.

\bibliographystyle{unsrt}
\bibliography{1}
\end{document}